\documentclass[final]{siamart250211}
\usepackage{amsfonts}
\usepackage{todonotes}
\usepackage{amsfonts,amsmath,amssymb,stackengine}
\usepackage{mathrsfs,mathtools,stmaryrd,wasysym}
\usepackage{enumerate}
\usepackage{xspace,mydef}
\usepackage{esint}
\usepackage{graphicx,psfrag}
\usepackage{hyperref}
\usepackage{pgf,tikz,pgfplots}
\usepackage{pstricks-add}
\usepackage{enumitem}
\usepackage{yfonts}
\usepackage{rotating}

\usetikzlibrary{arrows}
\DeclareMathAlphabet{\mathpzc}{OT1}{pzc}{m}{it}

\usepackage[notcite,notref]{showkeys}

\newsiamremark{remark}{Remark}
\newsiamremark{assumption}{Assumption}

\newcommand{\norm}[1]{{\left\vert\kern-0.25ex\left\vert\kern-0.25ex\left\vert #1
\right\vert\kern-0.25ex\right\vert\kern-0.25ex\right\vert}}

\newcommand{\TheTitle}{A Class of Monotone Quasi-Newtonian Flows Coupled to Nonlinear Equations:
	Analysis and Finite Element Approximation}
\newcommand{\ShortTitle}{Quasi-Newtonian Flows Coupled to Nonlinear Equations}
\newcommand{\TheAuthors}{Gilberto Campa\~na, Victor Vallejo and Jesus Vellojin}

\headers{\ShortTitle}{\TheAuthors}

\title{{\TheTitle}\thanks{GC is partially supported by ANID through FONDECYT project 11260204.}}

\author{Gilberto Campa\~na\thanks{Departamento de Ciencias, Universidad T\'ecnica Federico Santa Mar\'ia, Valpara\'iso, Chile. (\email{gilberto.campana@usm.cl}).}
\and
Victor Vallejo\thanks{Departamento de Ciencias, Universidad T\'ecnica Federico Santa Mar\'ia, Valpara\'iso, Chile. (\email{victor.vallejo@usm.cl}).}
\and
Jesus Vellojin\thanks{Departamento de Ciencias, Universidad T\'ecnica Federico Santa Mar\'ia, Valpara\'iso, Chile. (\email{jesus.vellojinm@usm.cl}).}
}

\ifpdf
\hypersetup{
  pdftitle={\TheTitle},
  pdfauthor={\TheAuthors}
}
\fi

\date{Draft version of \today.}

\begin{document}

\maketitle

\begin{abstract}
We analyze a monotone quasi-Newtonian flow coupled with nonlinear heat and concentration equations. In the model, the viscosity depends on the temperature, the concentration, and the velocity gradient, and the diffusivities depend on the temperature and the concentration. We show the existence of weak solutions without smallness assumptions on the data, and their uniqueness under suitable smallness and regularity assumptions. For the numerical approximation, we consider the Taylor--Hood and mini elements. We prove that the discrete solutions converge, up to a subsequence, to a weak solution, and we derive a quasi-best approximation result and convergence rates. We present numerical experiments in two and three dimensions that illustrate the theory.
\end{abstract}

\begin{keywords}
Quasi-Newtonian flow, double diffusion, monotone operators, finite element approximation.
\end{keywords}

\begin{AMS}
35A01,
35A02,
35Q35,
65N15,
65N30,
76A05,
76R50
\end{AMS}

\section{Introduction}
\label{sec:intro}
In this work, we focus on the analysis and finite element approximation of a class of monotone quasi-Newtonian incompressible flows coupled with nonlinear heat and concentration equations. To do this, we consider $\Omega \subset \mathbb{R}^d$, with $d \in \{2,3\}$, to be an open and bounded domain with Lipschitz boundary $\partial\Omega$. To be precise, the model under study is described by the following problem:
\begin{equation}\label{eq:model}
	\left\{
\begin{array}{c}
{-}\mathsf{div}(\nu(T,C)\,\varepsilon(x,|\nabla \mathbf{u}|)\nabla\mathbf{u})+\nabla \mathsf{p}
  =  \mathbf{f}(T,C),
\quad
\mathsf{div}\,\mathbf{u}=0 \textrm{ in } \Omega,
\\
{-}\mathsf{div}(\kappa(T,C)\nabla T)+\mathbf{u}\cdot\nabla T=\mathsf{g}_1,\hspace{0.25cm} -\mathsf{div}(\mathcal{D}(T,C)\nabla C)+\mathbf{u}\cdot\nabla C=\mathsf{g}_2 \textrm{ in } \Omega,
\end{array}
\right.
\end{equation}
supplemented with $\mathbf{u}=\mathbf{0}$, $T=0$, $C=0$ on $\partial\Omega$. The unknowns of system~\eqref{eq:model} are the velocity field $\mathbf{u}$, the pressure $\mathsf{p}$, the temperature $T$, and the concentration $C$. The data of the model are the thermal and solutal diffusivities $\kappa$ and $\mathcal{D}$, the external density force $\mathbf{f}$, the external heat and solute sources $\mathsf{g}_1$ and $\mathsf{g}_2$, and the viscosity $\nu(T,C)\,\varepsilon(x,|\nabla\mathbf{u}|)$, which depends on the temperature and the concentration through $\nu$ and on the velocity gradient through $\varepsilon$. We refer the reader to \cite{MR3423229,MR3959546,MR5061536} for related flow and transport models with variable viscosity and diffusivities.

In these works, the fluid is Newtonian, i.e., its viscosity may depend on the temperature or the concentration, but not on the velocity gradient. We also refer to \cite{MR4344862,MR4659442,MR4404048} for mixed finite element methods for double-diffusion problems with Newtonian fluids. In \cite{MR4344862,MR4404048}, the well-posedness of the continuous and discrete problems is established by fixed-point arguments under smallness assumptions on the data. However, many fluids of practical interest, such as blood and polymer solutions, have a viscosity that depends on the shear rate, and their flow is often accompanied by heat and mass transfer. Examples include blood oxygenators, in which the blood is cooled and dissolved carbon dioxide is removed, and cooling and separation processes in chemical engineering.

The finite element approximation of quasi-Newtonian Stokes flows whose viscosity follows a power law or a Carreau law was analyzed in \cite{MR1069652,MR1211613,MR1213411,MR1301740}. See also \cite{MR3119721} for discontinuous Galerkin methods. For quasi-Newtonian flows coupled with a heat equation, the existence of stationary solutions with temperature-dependent viscosity and viscous heating was studied in \cite{MR1348583}, and the convergence of a decoupled fixed-point iteration was proved in \cite{MR1636360} under smallness assumptions on the data. Finite element approximations of thermally coupled quasi-Newtonian flows were analyzed in \cite{MR2320257,MR2794183,MR3070180}. More recently, in \cite{MR4610451} the authors studied the existence, regularity, and conditional uniqueness of weak solutions for stationary non-isothermal power-law fluids with temperature-dependent viscosity, and derived a priori error estimates for a finite element discretization of the Stokes--Carreau case under smallness and regularity assumptions. We also mention \cite{MR4379972} for the finite element approximation and preconditioning of anisothermal flows of implicitly constituted fluids, and \cite{MR2515781} for an existence result for large data in the unsteady case.

On the other hand, flows of generalized Newtonian fluids with a power-law index that depends on a concentration were analyzed in \cite{MR3023245,MR3262601}, where the existence of weak solutions is proved without smallness assumptions on the data. For this model, the convergence of a subsequence of finite element approximations to a weak solution was proved in \cite{MR3834434} in two dimensions, with skew-symmetric convective terms, and in \cite{MR3904139} in three dimensions, by means of a regularization of the momentum equation.

As far as we know, the existence of weak solutions and the convergence of finite element approximations without smallness assumptions on the data have not been established for problem \eqref{eq:model}. In this work, we aim to prove these results for a viscosity that depends on the temperature, the concentration, and the velocity gradient, and for diffusivities that depend on both scalar fields.

Throughout this manuscript, we will consider the following assumptions.

$\checkmark$ \emph{Viscosity:} We assume that, for $s\geq 0$, the function $x\mapsto\varepsilon(x,s)$ is Lebesgue measurable in $\Omega$. Inspired by \cite[Section 2]{MR3119721}, we assume that there exist constants $\varepsilon_{-}, \varepsilon_{+}>0$ such that
\begin{align}\label{assumption_1}
|\varepsilon(x,|\boldsymbol{\tau}|)\boldsymbol{\tau}-\varepsilon(x,|\boldsymbol{\omega}|)\boldsymbol{\omega}|
\leq \varepsilon_{+} |\boldsymbol{\tau}-\boldsymbol{\omega}|,\qquad \\ \label{assumption_2}
\varepsilon_{-}|\boldsymbol{\tau}-\boldsymbol{\omega}|^{2} \leq
(\varepsilon(x,|\boldsymbol{\tau}|)\boldsymbol{\tau}-\varepsilon(x,|\boldsymbol{\omega}|)\boldsymbol{\omega}):(\boldsymbol{\tau}-\boldsymbol{\omega}),
\end{align}
for all $\boldsymbol{\tau},\boldsymbol{\omega}\in\mathbb{R}^{d\times d}$ and a.e.\ $x\in\Omega$. The function $\nu\in C^{0,1}(\mathbb{R}^2)$ is strictly positive and bounded: there exist positive constants $\nu_{-}$ and $\nu_{+}$ such that $\nu_{-}\leq\nu(r,s)\leq\nu_{+}$ for all $(r,s)\in\mathbb{R}^2$. Moreover, $\nu$ is Lipschitz continuous with constant $\mathcal{L}_{\nu}$, i.e.,
\[
|\nu(r_1,s_1)-\nu(r_2,s_2)|\leq\mathcal{L}_{\nu}\left(|r_1-r_2|+|s_1-s_2|\right),
\]
for all $(r_1,s_1),(r_2,s_2)\in\mathbb{R}^2$.

$\checkmark$ \emph{Diffusivities:} The coefficients $\kappa,\mathcal{D}\in C^{0,1}(\mathbb{R}^2)$ are strictly positive and uniformly bounded: there exist positive constants $\kappa_{-}$, $\kappa_{+}$, $\mathcal{D}_{-}$, and $\mathcal{D}_{+}$ such that $\kappa_{-}\leq\kappa(r,s)\leq\kappa_{+}$ and $\mathcal{D}_{-}\leq\mathcal{D}(r,s)\leq\mathcal{D}_{+}$ for all $(r,s)\in\mathbb{R}^2$. Moreover, $\kappa$ and $\mathcal{D}$ are Lipschitz continuous with constants $\mathcal{L}_{\kappa}$ and $\mathcal{L}_{\mathcal{D}}$, respectively, i.e., for all $(r_1,s_1),(r_2,s_2)\in\mathbb{R}^2$ we have that
\[
\begin{aligned}
|\kappa(r_1,s_1)-\kappa(r_2,s_2)|&\leq\mathcal{L}_{\kappa}\left(|r_1-r_2|+|s_1-s_2|\right),
\\
|\mathcal{D}(r_1,s_1)-\mathcal{D}(r_2,s_2)|&\leq\mathcal{L}_{\mathcal{D}}\left(|r_1-r_2|+|s_1-s_2|\right).
\end{aligned}
\]

$\checkmark$ \emph{External forces:} As in \cite{MR4337455,MR4659334}, we consider Lipschitz continuous functions $\mathbf{f}_1,\mathbf{f}_2$ with constants $C_{\mathbf{f}_1}$ and $C_{\mathbf{f}_2}$, respectively, and $\mathbf{f}_1(0)=\mathbf{f}_2(0)=\mathbf{0}$. In particular, $|\mathbf{f}_i(s)|\leq C_{\mathbf{f}_i}|s|$ for every $s\in\mathbb{R}$ and $i\in\{1,2\}$. On the other hand, we set $\mathbf{f}_0\in \mathbf{H}^{-1}(\Omega)$. Thus,
\[
\mathbf{f}(T,C):= \mathbf{f}_0+\mathbf{f}_1(T)+\mathbf{f}_2(C).
\]
In addition, we assume that the external sources $\mathsf{g}_1$ and $\mathsf{g}_2$ belong to $H^{-1}(\Omega)$.

Although classical non-Newtonian fluid models involve the symmetric gradient, we work with the full gradient, as in \cite{MR1370148,MR3232446,MR3959546}. This choice is consistent with the open boundary conditions used in the numerical applications, following the distinction between open and traction boundary conditions discussed in \cite[Section 1]{LIU}.

We summarize the main contributions of this work as follows:

$\checkmark$ \emph{Existence and uniqueness:} We prove that problem \eqref{eq:model} admits weak solutions without smallness assumptions on the data. The proof combines a Galerkin method with Brouwer's fixed-point theorem.
We also prove uniqueness in two and three dimensions under smallness conditions.%

$\checkmark$ \emph{Finite element approximation:} We discretize problem \eqref{eq:model} with the Taylor--Hood and mini elements for the velocity and the pressure, continuous piecewise polynomials for the temperature and the concentration, and a skew-symmetric form of the convective terms. We prove the existence of discrete solutions and the strong convergence of a subsequence to a weak solution, without smallness assumptions on the data and without additional regularity of the solution.

$\checkmark$ \emph{Error estimates:} Under the assumptions of our uniqueness results, we derive a quasi-best approximation estimate. For sufficiently regular solutions, this estimate yields convergence rates of order $h^2$ for the Taylor--Hood element and of order $h$ for the mini element.

$\checkmark$ \emph{A convergent fixed-point algorithm:} We propose a fixed-point algorithm that decouples the flow equations and the scalar equations.
Inspired by \cite{MR1636360}, we prove its linear convergence to the discrete solution under explicit smallness conditions.%

$\checkmark$ \emph{Numerical experiments:} We confirm the theoretical convergence rates in two and three dimensions, study the effect of the rheological contrast and of small coefficients on the accuracy and on the number of fixed-point iterations, and present two applications: staged cooling in a hollow-fiber oxygenator and cooling with solute removal in two connected tanks.

We organize this work as follows. In Section~\ref{sec:notation}, we present the notation and some preliminary results. Sections~\ref{sec:coupled} and~\ref{sec:fem} are devoted to the analysis of the continuous and discrete problems. Finally, in Section~\ref{sec:numerics}, we present a fixed-point algorithm, its convergence analysis, and the numerical experiments that illustrate our results.

\section{Notation and preliminary remarks}
\label{sec:notation}
We now present the standard notation and a classical framework within which we develop our analysis.

\subsection{Notation}
Throughout the text, we use the classical notation for Lebesgue and Sobolev spaces and denote the vector-valued counterparts of these spaces with bold letters. In particular, we consider the space
$
\mathbf{V}(\Omega):=\{\mathbf{v}\in\mathbf{H}_0^1(\Omega):\text{div }\mathbf{v}=0\}.
$

Let $\mathcal{X}$ and $\mathcal{Y}$ be Banach function spaces. We write $\mathcal{X} \hookrightarrow \mathcal{Y}$ to indicate that $\mathcal{X}$ is continuously embedded in $\mathcal{Y}$. We denote by $\mathcal{X}'$ and $\|\cdot\|_{\mathcal{X}}$ the dual space and the norm of $\mathcal{X}$, respectively. On the other hand, given $\mathtt{p}\in(1, \infty)$, we denote by $\mathtt{p}\prime$ its H\"older conjugate, which is such that $1/\mathtt{p} + 1/\mathtt{p}\prime = 1$. The relation $a\lesssim b$ states that $a\leq Cb$ with a positive constant $C$ that does not depend on $a$ and $b$, respectively. Note that the value of the constant $C$ might change at each occurrence.

\subsection{A class of quasi-Newtonian models}
In this subsection, we analyze the well-posedness of a class of quasi-Newtonian systems. To do this, we consider an external density force $\mathbf{f}\in\mathbf{H}^{-1}(\Omega)$ and a function $\nu\in L^\infty(\Omega)$ that satisfies $\nu_{-}\leq\nu(x)\leq\nu_{+}$ for a.e.\ $x\in\Omega$, and we are interested in finding the pair $(\mathbf{u},\mathsf{p})\in \mathbf{H}_0^1(\Omega)\times L_0^2(\Omega)$ such that
\begin{equation}\label{eq:fluid_problem}
\int_\Omega\left(\nu\,\varepsilon(x,|\nabla\mathbf{u}|)\nabla\mathbf{u}:\nabla\mathbf{v}-\mathsf{p}\,\mathsf{div}\,\mathbf{v} \right)\mathrm{d}x
=
\displaystyle \langle\mathbf{f},\mathbf{v}\rangle,\quad
\displaystyle
\int_{\Omega}\mathsf{q}\,\mathsf{div}\,\mathbf{u}\,\mathrm{d}x=0,
\end{equation}
for all $\mathbf{v}\in \mathbf{H}_0^{1}(\Omega)$ and $\mathsf{q}\in L_0^{2}(\Omega)$. Here, $\langle\cdot,\cdot\rangle$ represents the duality pairing between $\mathbf{H}^{-1}(\Omega)$ and $\mathbf{H}_0^{1}(\Omega)$. We note that owing to our assumptions on the data and the conditions \eqref{assumption_1}--\eqref{assumption_2}, all terms in the weak formulation \eqref{eq:fluid_problem} are well-defined.

Since the divergence operator is surjective from $\mathbf{H}_0^1(\Omega)$ to $L_0^2(\Omega)$, there exists $\beta>0$ such that the inf-sup condition holds (see for more details \cite[Chapter I, Section 5.1]{MR851383}, \cite[Corollary B.71]{Guermond-Ern}):
\begin{equation}\label{eq:infsup_cont}
\sup_{\boldsymbol{0}\neq\mathbf{v}\in\mathbf{H}_0^1(\Omega)}\dfrac{(\mathsf{q},\text{div }\mathbf{v})_{L^2(\Omega)}}{\|\nabla \mathbf{v}\|_{\mathbf{L}^2(\Omega)}}\geq \beta \|\mathsf{q}\|_{L^2(\Omega)}\quad\forall \mathsf{q}\in L_0^2(\Omega).
\end{equation}

In order to study problem \eqref{eq:fluid_problem}, it is convenient to introduce the map
\begin{equation}\label{eq:operator_Aa}
\mathcal{A}:\mathbf{H}_0^1(\Omega)\to\mathbf{H}^{-1}(\Omega)\quad
\langle\mathcal{A}(\mathbf{v}),\mathbf{w}\rangle
:=\int_\Omega\nu\,\varepsilon(x,|\nabla\mathbf{v}|)\nabla\mathbf{v}:\nabla\mathbf{w}\,\mathrm{d}x.
\end{equation}

\begin{proposition}[Properties of $\mathcal{A}$]\label{proposition_the_map_A} The map $\mathcal{A}$ satisfies the following properties:

(i) For every $\mathbf{u},\mathbf{v},\mathbf{w}\in\mathbf{H}_0^1(\Omega)$, we have that
\begin{align}\label{eq:prop_01}
\langle\mathcal{A}(\mathbf{u})-\mathcal{A}(\mathbf{v}),\mathbf{w}\rangle\leq \nu_{+}\varepsilon_{+}\|\nabla(\mathbf{u}-\mathbf{v})\|_{\mathbf{L}^2(\Omega)}\|\nabla \mathbf{w}\|_{\mathbf{L}^2(\Omega)}.
\end{align}

(ii) $\mathcal{A}$ is strongly monotone in $\mathbf{H}_0^1(\Omega)$: If $\mathbf{v},\mathbf{w}\in\mathbf{H}_0^1(\Omega)$, then
\begin{align}\label{eq:prop_02}
\langle\mathcal{A}(\mathbf{v})-\mathcal{A}(\mathbf{w}),\mathbf{v}-\mathbf{w}\rangle\geq \nu_{-}\varepsilon_{-}\|\nabla(\mathbf{v}-\mathbf{w})\|_{\mathbf{L}^2(\Omega)}^2.
\end{align}

(iii) $\mathcal{A}$ is coercive in $\mathbf{H}_0^1(\Omega)$:
\begin{align}\label{eq:prop_03}
\langle\mathcal{A}(\mathbf{v}),\mathbf{v}\rangle /\|\nabla \mathbf{v}\|_{\mathbf{L}^2(\Omega)}\to +\infty ,\quad\text{as}\quad \|\nabla\mathbf{v}\|_{\mathbf{L}^2(\Omega)}\to +\infty.
\end{align}

(iv) $\mathcal{A}$ is hemicontinuous in $\mathbf{H}_0^1(\Omega)$: If $\mathbf{u},\mathbf{v},\mathbf{w}\in\mathbf{H}_0^1(\Omega)$, then the mapping
\begin{align}\label{eq:prop_05}
t\mapsto \langle\mathcal{A}(\mathbf{u}+t\mathbf{v}),\mathbf{w}\rangle
\end{align}
is continuous from $\mathbb{R}$ into $\mathbb{R}$.
\end{proposition}

\begin{proof}
The properties \eqref{eq:prop_01} and \eqref{eq:prop_02} are a direct consequence of assumptions \eqref{assumption_1}--\eqref{assumption_2}, the estimates $\nu_{-}\leq\nu\leq\nu_{+}$, and H\"older's inequality. Property \eqref{eq:prop_03} follows from \eqref{eq:prop_02} with $\mathbf{w}=\mathbf{0}$. Finally, to derive \eqref{eq:prop_05} we invoke \eqref{eq:prop_01} to establish that, if $t_n\to t$ in $\mathbb{R}$, then $|\langle\mathcal{A}(\mathbf{u}+t_n\mathbf{v}),\mathbf{w}\rangle-\langle\mathcal{A}(\mathbf{u}+t\mathbf{v}),\mathbf{w}\rangle|\lesssim|t_n-t|\|\nabla\mathbf{v}\|_{\mathbf{L}^2(\Omega)}\|\nabla\mathbf{w}\|_{\mathbf{L}^2(\Omega)}\to 0$. This concludes the proof.
\end{proof}

Having introduced the operator $\mathcal{A}$ and the space $\mathbf{V}(\Omega)$, we invoke de Rham's theorem \cite[Theorem B.73]{Guermond-Ern} to conclude that problem \eqref{eq:fluid_problem} is equivalent to the reduced formulation: Find $\mathbf{u}\in \mathbf{V}(\Omega)$ such that
\begin{align}\label{eq:problem_reduced}
 \langle\mathcal{A}(\mathbf{u}),\mathbf{v}\rangle
 =  \langle \mathbf{f}, \mathbf{v} \rangle \quad \forall\mathbf{v}\in \mathbf{V}(\Omega).
\end{align}

\begin{theorem}[well-posedness]\label{th:well_posedness_darcy_forchheimer} Let $\mathbf{f}\in\mathbf{H}^{-1}(\Omega)$ and let $\varepsilon$ be such that it satisfies the assumptions \eqref{assumption_1} and \eqref{assumption_2}. Then, there exists a unique solution $(\mathbf{u},\mathsf{p})\in \mathbf{H}_0^{1}(\Omega)\times L_0^{2}(\Omega)$ of problem \eqref{eq:fluid_problem}, which satisfies the following estimate
\begin{align}\label{eq:estimate_quasi}
\|\nabla\mathbf{u}\|_{\mathbf{L}^{2}(\Omega)} +\| \mathsf{p}\|_{L^{2}(\Omega)}&\lesssim \|\mathbf{f}\|_{\mathbf{H}^{-1}(\Omega)}.
\end{align}
Here, the hidden constant depends only on the constants $\beta$, $\nu_{-}$, $\nu_{+}$, $\varepsilon_{-}$, and $\varepsilon_{+}$.
\end{theorem}

\begin{proof}
First, we observe that \eqref{eq:fluid_problem} and \eqref{eq:problem_reduced} are equivalent. Hence, we study problem \eqref{eq:problem_reduced}. Since $\mathbf{V}(\Omega)$ is a closed subspace of $\mathbf{H}_0^1(\Omega)$, it is a separable Hilbert space, and the properties in Proposition~\ref{proposition_the_map_A} also hold for the restriction of $\mathcal{A}$ to $\mathbf{V}(\Omega)$, which is understood as a map from $\mathbf{V}(\Omega)$ into $\mathbf{V}(\Omega)'$. Note that \eqref{eq:prop_02} shows that this restriction is strictly monotone. Thus, we apply \cite[Theorem 9.14-1]{MR3136903} to guarantee the existence of a unique $\mathbf{u}\in\mathbf{V}(\Omega)$ that satisfies \eqref{eq:problem_reduced}. Second, the existence of a unique $\mathsf{p}\in L_0^2(\Omega)$ such that $(\mathbf{u},\mathsf{p})$ solves \eqref{eq:fluid_problem} follows from de Rham's theorem \cite[Theorem B.73]{Guermond-Ern} and the inf-sup condition \eqref{eq:infsup_cont}. Finally, to derive the estimate in \eqref{eq:estimate_quasi}, we set $\mathbf{v}=\mathbf{u}$ in the first equation of problem \eqref{eq:fluid_problem} and invoke \eqref{eq:prop_02} with $\mathbf{w}=\mathbf{0}$. This yields $\nu_{-}\varepsilon_{-}\|\nabla\mathbf{u}\|_{\mathbf{L}^2(\Omega)}\leq\|\mathbf{f}\|_{\mathbf{H}^{-1}(\Omega)}$. The bound for $\mathsf{p}$ follows from the inf-sup condition \eqref{eq:infsup_cont}, the first equation of \eqref{eq:fluid_problem}, and \eqref{eq:prop_01} with $\mathbf{v}=\mathbf{0}$, which produce $\beta\|\mathsf{p}\|_{L^2(\Omega)}\leq\nu_{+}\varepsilon_{+}\|\nabla\mathbf{u}\|_{\mathbf{L}^2(\Omega)}+\|\mathbf{f}\|_{\mathbf{H}^{-1}(\Omega)}$. This concludes the proof.
\end{proof}

\section{The coupled problem}
\label{sec:coupled}
In the following, we show the existence of suitable weak solutions for the coupled problem \eqref{eq:model} and derive a uniqueness result under appropriate assumptions. To continue, we introduce the concept of a weak solution.

\subsection{Weak solution}
\label{sec:weak_solution}
The problem \eqref{eq:model} is analyzed in the following context.

\begin{definition}[weak solution]\label{def:weak_solution}
Let $\mathbf{f}_0\in\mathbf{H}^{-1}(\Omega)$, $\mathsf{g}_1\in H^{-1}(\Omega)$, and $\mathsf{g}_2\in H^{-1}(\Omega)$. We say that $(\mathbf{u},\mathsf{p},T,C)\in\mathbf{H}_0^1(\Omega)\times L_0^2(\Omega)\times H_0^1(\Omega)\times H_0^1(\Omega)$ is a weak solution to \eqref{eq:model} if
\begin{align}\nonumber
\int_\Omega\left(\nu(T,C)\,\varepsilon(x,|\nabla\mathbf{u}|)\nabla\mathbf{u}:\nabla\mathbf{v}-\mathsf{p}\operatorname{div}\mathbf{v}\right)\mathrm{d}x=\langle\mathbf{f}(T,C),\mathbf{v}\rangle,\hspace{1cm}
\\ \label{eq:weak_coupled}
\int_\Omega\mathsf{q}\operatorname{div}\mathbf{u}\,\mathrm{d}x=0,
\quad
\int_\Omega\left(\kappa(T,C)\nabla T\cdot\nabla S+(\mathbf{u}\cdot\nabla T)S\right)\mathrm{d}x=\langle\mathsf{g}_1,S\rangle,
\\ \nonumber
\int_\Omega\left(\mathcal{D}(T,C)\nabla C\cdot\nabla R+(\mathbf{u}\cdot\nabla C)R\right)\mathrm{d}x=\langle\mathsf{g}_2,R\rangle,\hspace{1cm}
\end{align}
for all $(\mathbf{v},\mathsf{q},S,R)\in\mathbf{H}_0^1(\Omega)\times L_0^2(\Omega)\times H_0^1(\Omega)\times H_0^1(\Omega)$.
\end{definition}
Under the assumptions presented in Section \ref{sec:intro}, every term in system \eqref{eq:weak_coupled} is well-defined.
To simplify the presentation of the analysis, given $T,C\in H_0^1(\Omega)$, we introduce the map $\mathcal{A}_{T,C}:\mathbf{H}_0^1(\Omega)\to\mathbf{H}^{-1}(\Omega)$ defined by
\begin{equation}\label{eq:operator_ATC}
\langle\mathcal{A}_{T,C}(\mathbf{v}),\mathbf{w}\rangle
:=\int_\Omega\nu(T,C)\,\varepsilon(x,|\nabla\mathbf{v}|)\nabla\mathbf{v}:\nabla\mathbf{w}\,\mathrm{d}x,
\end{equation}
which is the map defined in \eqref{eq:operator_Aa} with $\nu$ replaced by $\nu(T,C)$. Since $\nu(T,C)\in L^\infty(\Omega)$ and $\nu_{-}\leq\nu(T,C)\leq\nu_{+}$ a.e.\ in $\Omega$, Proposition~\ref{proposition_the_map_A} applies to $\mathcal{A}_{T,C}$, with constants that do not depend on $T$ and $C$. With this notation, the first equation in \eqref{eq:weak_coupled} reads as follows: $\langle\mathcal{A}_{T,C}(\mathbf{u}),\mathbf{v}\rangle-\int_\Omega\mathsf{p}\operatorname{div}\mathbf{v}\,\mathrm{d}x=\langle\mathbf{f}(T,C),\mathbf{v}\rangle$.

\subsection{Existence of solutions}
\label{sec:existence}

In what follows, we derive an existence result for the coupled problem \eqref{eq:weak_coupled}. To do this, we introduce the spaces $\mathbb{X}:=\mathbf{H}_0^1(\Omega)\times H_0^1(\Omega)\times H_0^1(\Omega)$ and $\mathbb{X}':=\mathbf{H}^{-1}(\Omega)\times H^{-1}(\Omega)\times H^{-1}(\Omega)$, and we use the following notation:
\begin{align*}
& \forall\,(\mathbf{v},S,R)\in\mathbb{X}:
\quad
\|(\mathbf{v},S,R)\|_{\mathbb{X}}
:=
\left[
\|\nabla\mathbf{v}\|_{\mathbf{L}^2(\Omega)}^2+\|\nabla S\|_{\mathbf{L}^2(\Omega)}^2+\|\nabla R\|_{\mathbf{L}^2(\Omega)}^2
\right]^{1/2},
\\
& \forall\,(\mathbf{f},\mathsf{g}_1,\mathsf{g}_2)\in\mathbb{X}':
\quad
\|(\mathbf{f},\mathsf{g}_1,\mathsf{g}_2)\|_{\mathbb{X}'}
:=
\left[
\|\mathbf{f}\|_{\mathbf{H}^{-1}(\Omega)}^2+\|\mathsf{g}_1\|_{H^{-1}(\Omega)}^2+\|\mathsf{g}_2\|_{H^{-1}(\Omega)}^2
\right]^{1/2}.
\end{align*}

\begin{theorem}[existence of solutions]\label{th:existence_coupled}
Let $d\in\{2,3\}$, and let $\Omega\subset\mathbb{R}^d$ be an open and bounded domain with Lipschitz boundary $\partial\Omega$. Let $\varepsilon$, $\nu$, $\kappa$, $\mathcal{D}$, $\mathbf{f}_1$, and $\mathbf{f}_2$ be such that the assumptions of Section~\ref{sec:intro} hold. If $\mathbf{f}_0\in\mathbf{H}^{-1}(\Omega)$ and $\mathsf{g}_1,\mathsf{g}_2\in H^{-1}(\Omega)$, then problem \eqref{eq:weak_coupled} has at least one solution $(\mathbf{u},\mathsf{p},T,C)\in\mathbf{H}_0^1(\Omega)\times L_0^2(\Omega)\times H_0^1(\Omega)\times H_0^1(\Omega)$. Moreover, we have
\begin{align}
\label{eq:estimate_TC_coupled}
\|\nabla T\|_{\mathbf{L}^2(\Omega)}
& \leq \kappa_{-}^{-1}\|\mathsf{g}_1\|_{H^{-1}(\Omega)},
\quad
\|\nabla C\|_{\mathbf{L}^2(\Omega)}
\leq \mathcal{D}_{-}^{-1}\|\mathsf{g}_2\|_{H^{-1}(\Omega)},
\\
\label{eq:estimate_up_coupled}
\|\nabla\mathbf{u}\|_{\mathbf{L}^2(\Omega)}
& \leq (\nu_{-}\varepsilon_{-})^{-1}\Lambda,
\quad
\|\mathsf{p}\|_{L^2(\Omega)}
\leq \beta^{-1}\left(1+\nu_{+}\varepsilon_{+}(\nu_{-}\varepsilon_{-})^{-1}\right)\Lambda.
\end{align}
Here, $\mathcal{C}_{2\hookrightarrow2}$ denotes the best constant in the Sobolev embedding $H_0^1(\Omega)\hookrightarrow L^2(\Omega)$, and $\Lambda:=\|\mathbf{f}_0\|_{\mathbf{H}^{-1}(\Omega)}
+
\mathcal{C}_{2\hookrightarrow2}^2
\left(
C_{\mathbf{f}_1}\kappa_{-}^{-1}\|\mathsf{g}_1\|_{H^{-1}(\Omega)}
+
C_{\mathbf{f}_2}\mathcal{D}_{-}^{-1}\|\mathsf{g}_2\|_{H^{-1}(\Omega)}
\right).$
\end{theorem}

\begin{proof}
We adapt the proof of \cite[Theorem 8]{MR4911804} to our case and divide the proof into five steps. Throughout the proof, we use the map $\mathcal{A}_{T,C}$ defined in \eqref{eq:operator_ATC} and its properties stated in Proposition~\ref{proposition_the_map_A}. The same notation is used when $T$ and $C$ are replaced by other functions in $H_0^1(\Omega)$; for instance, $\mathcal{A}_{T_k,C_k}$ denotes the map defined in \eqref{eq:operator_ATC} with $T$ and $C$ replaced by $T_k$ and $C_k$, respectively.

\emph{Step 1.} \emph{(A mapping $\Phi$).} Let $\mathbf{u},\mathbf{v}\in\mathbf{V}(\Omega)$, and let $T,C,S,R\in H_0^1(\Omega)$ be arbitrary. To simplify the presentation of the proof, we define $\boldsymbol{\mathcal{U}}:=(\mathbf{u},T,C)$ and $\boldsymbol{\mathcal{V}}:=(\mathbf{v},S,R)$. On the other hand, given $\lambda>0$, we introduce the mapping $\Phi$ from the space $\mathbf{V}(\Omega)\times H_0^1(\Omega)\times H_0^1(\Omega)$ into its dual space as follows
\begin{multline}\label{eq:definition_Phi}
\langle\Phi(\boldsymbol{\mathcal{U}}),\boldsymbol{\mathcal{V}}\rangle
:=
\langle\mathcal{A}_{T,C}(\mathbf{u})-\mathbf{f}(T,C),\mathbf{v}\rangle
+\lambda\Big(-\langle\mathsf{g}_1,S\rangle-\langle\mathsf{g}_2,R\rangle
\\
+\int_\Omega\left(\kappa(T,C)\nabla T\cdot\nabla S+(\mathbf{u}\cdot\nabla T)S\right)\mathrm{d}x
+\int_\Omega\left(\mathcal{D}(T,C)\nabla C\cdot\nabla R+(\mathbf{u}\cdot\nabla C)R\right)\mathrm{d}x\Big).
\end{multline}

$\checkmark$ \emph{Continuity of $\Phi$.} Let $\boldsymbol{\mathcal{U}}_n=(\mathbf{u}_n,T_n,C_n)\to\boldsymbol{\mathcal{U}}=(\mathbf{u},T,C)$ in $\mathbf{V}(\Omega)\times H_0^1(\Omega)\times H_0^1(\Omega)$ as $n\uparrow\infty$. We now invoke \eqref{eq:prop_01} for $\mathcal{A}_{T_n,C_n}$, the condition \eqref{assumption_1}, the upper bounds of $\nu$, $\kappa$, and $\mathcal{D}$, the assumption regarding $\mathbf{f}_1$ and $\mathbf{f}_2$, H\"older's inequality, the Sobolev embeddings $H_0^1(\Omega)\hookrightarrow L^2(\Omega)$ and $H_0^1(\Omega)\hookrightarrow L^4(\Omega)$, and the Lebesgue dominated convergence theorem to produce that
\begin{multline*}
|\langle\Phi(\boldsymbol{\mathcal{U}}_n)-\Phi(\boldsymbol{\mathcal{U}}),\boldsymbol{\mathcal{V}}\rangle|/\|\boldsymbol{\mathcal{V}}\|_{\mathbb{X}}
\lesssim \|\nabla(\mathbf{u}_n-\mathbf{u})\|_{\mathbf{L}^2(\Omega)}
+\|(\nu(T_n,C_n)-\nu(T,C))\nabla\mathbf{u}\|_{\mathbf{L}^2(\Omega)}
\\+\|\nabla(T_n-T)\|_{\mathbf{L}^2(\Omega)}+\|\nabla(C_n-C)\|_{\mathbf{L}^2(\Omega)}
+\|(\kappa(T_n,C_n)-\kappa(T,C))\nabla T\|_{\mathbf{L}^2(\Omega)}
\\
+\|(\mathcal{D}(T_n,C_n)-\mathcal{D}(T,C))\nabla C\|_{\mathbf{L}^2(\Omega)}
+\|\nabla(\mathbf{u}_n-\mathbf{u})\|_{\mathbf{L}^2(\Omega)}\left(\|\nabla T_n\|_{\mathbf{L}^2(\Omega)}+\|\nabla C_n\|_{\mathbf{L}^2(\Omega)}\right)
\\
+\|\nabla\mathbf{u}\|_{\mathbf{L}^2(\Omega)}\left(\|\nabla(T_n-T)\|_{\mathbf{L}^2(\Omega)}+\|\nabla(C_n-C)\|_{\mathbf{L}^2(\Omega)}\right)\to 0,\quad\text{as}\quad n\uparrow \infty.
\end{multline*}
$\checkmark$ \emph{The map $\Phi$ satisfies $\langle\Phi(\boldsymbol{\mathcal{V}}),\boldsymbol{\mathcal{V}}\rangle\geq0$ for all $\boldsymbol{\mathcal{V}}\in\mathbf{V}(\Omega)\times H_0^1(\Omega)\times H_0^1(\Omega)$ such that $\|\boldsymbol{\mathcal{V}}\|_{\mathbb{X}}=\delta$, with $\delta$ defined in \eqref{eq:definition_delta}}. If we replace $\boldsymbol{\mathcal{U}}$ by $\boldsymbol{\mathcal{V}}$ in \eqref{eq:definition_Phi} and use that $\int_\Omega(\mathbf{v}\cdot\nabla S)S\,\mathrm{d}x=\int_\Omega(\mathbf{v}\cdot\nabla R)R\,\mathrm{d}x=0$, together with \eqref{eq:prop_02} with $\mathbf{w}=\mathbf{0}$ for $\mathcal{A}_{S,R}$, the lower bounds of $\kappa$ and $\mathcal{D}$, H\"older's inequality, the Sobolev embedding $H_0^1(\Omega)\hookrightarrow L^2(\Omega)$, and the estimate
\begin{eqnarray}\label{eq:bound_f_TC}
\|\mathbf{f}(T,C)\|_{\mathbf{H}^{-1}(\Omega)}
\leq
\|\mathbf{f}_0\|_{\mathbf{H}^{-1}(\Omega)}
+\mathcal{C}_{2\hookrightarrow2}^2
\left(C_{\mathbf{f}_1}\|\nabla T\|_{\mathbf{L}^2(\Omega)}+C_{\mathbf{f}_2}\|\nabla C\|_{\mathbf{L}^2(\Omega)}\right),
\end{eqnarray}
we can conclude that
\begin{multline*}
\langle\Phi(\boldsymbol{\mathcal{V}}),\boldsymbol{\mathcal{V}}\rangle
\geq
\nu_{-}\varepsilon_{-}\|\nabla\mathbf{v}\|_{\mathbf{L}^2(\Omega)}^2
+\lambda\kappa_{-}\|\nabla S\|_{\mathbf{L}^2(\Omega)}^2
+\lambda\mathcal{D}_{-}\|\nabla R\|_{\mathbf{L}^2(\Omega)}^2
-\mathcal{C}_{2\hookrightarrow2}^2\|\nabla\mathbf{v}\|_{\mathbf{L}^2(\Omega)}\\
\cdot\big(C_{\mathbf{f}_1}\|\nabla S\|_{\mathbf{L}^2(\Omega)}+C_{\mathbf{f}_2}\|\nabla R\|_{\mathbf{L}^2(\Omega)}\big)
-\|(\mathbf{f}_0,\lambda\mathsf{g}_1,\lambda\mathsf{g}_2)\|_{\mathbb{X}'}
\|(\mathbf{v},S,R)\|_{\mathbb{X}}.
\end{multline*}
Young's inequality, combined with the choice of $\lambda:=\max\left\{
1,\,
\frac{2\,\mathcal{C}_{2\hookrightarrow2}^4C_{\mathbf{f}_1}^2}{\nu_{-}\varepsilon_{-}\kappa_{-}},\,
\frac{2\,\mathcal{C}_{2\hookrightarrow2}^4C_{\mathbf{f}_2}^2}{\nu_{-}\varepsilon_{-}\mathcal{D}_{-}}
\right\}$, allows us to conclude that
$
\langle\Phi(\boldsymbol{\mathcal{V}}),\boldsymbol{\mathcal{V}}\rangle
\geq
\gamma_{\lambda}\|(\mathbf{v},S,R)\|_{\mathbb{X}}^2
-\|(\mathbf{f}_0,\lambda\mathsf{g}_1,\lambda\mathsf{g}_2)\|_{\mathbb{X}'}\|(\mathbf{v},S,R)\|_{\mathbb{X}},
$
where $\gamma_{\lambda}:=\frac{1}{2}\min\{\nu_{-}\varepsilon_{-},\lambda\kappa_{-},\lambda\mathcal{D}_{-}\}$. From this, we can conclude that $\langle\Phi(\boldsymbol{\mathcal{V}}),\boldsymbol{\mathcal{V}}\rangle\geq0$ for all $\boldsymbol{\mathcal{V}}\in\mathbf{V}(\Omega)\times H_0^1(\Omega)\times H_0^1(\Omega)$ such that
\begin{align}\label{eq:definition_delta}
\|(\mathbf{v},S,R)\|_{\mathbb{X}}
=
2\left(\min\{\nu_{-}\varepsilon_{-},\lambda\kappa_{-},\lambda\mathcal{D}_{-}\}\right)^{-1}\|(\mathbf{f}_0,\lambda\mathsf{g}_1,\lambda\mathsf{g}_2)\|_{\mathbb{X}'}
=:\delta.
\end{align}

\emph{Step 2.} \emph{(Galerkin approximation).} We use arguments presented in \cite[Theorem 2.6]{MR3014456}. We observe that, since the spaces $\mathbf{V}(\Omega)$ and $H_0^1(\Omega)$ are separable, there exist sequences $\{\mathbb{V}_k\}_{k\in\mathbb{N}}$ and $\{\mathbb{W}_k\}_{k\in\mathbb{N}}$ of finite-dimensional subspaces of $\mathbf{V}(\Omega)$ and $H_0^1(\Omega)$, respectively, such that for every $k\in\mathbb{N}$,
\begin{equation}\label{eq:density}
\mathbb{V}_k\subset\mathbb{V}_{k+1},
\quad
\mathbb{W}_k\subset\mathbb{W}_{k+1},
\qquad
\overline{\cup_{k\in\mathbb{N}}\mathbb{V}_k}=\mathbf{V}(\Omega),
\quad
\overline{\cup_{k\in\mathbb{N}}\mathbb{W}_k}=H_0^1(\Omega).
\end{equation}
Therefore, setting $\mathbb{X}_k:=\mathbb{V}_k\times\mathbb{W}_k\times\mathbb{W}_k$, we define a Galerkin approximation $\boldsymbol{\mathcal{U}}_k=(\mathbf{u}_k,T_k,C_k)\in\mathbb{X}_k$ as follows $
\langle\Phi(\boldsymbol{\mathcal{U}}_k),\boldsymbol{\mathcal{V}}_k\rangle=0$
for all $\boldsymbol{\mathcal{V}}_k\in\mathbb{X}_k$. Moreover, due to the continuity of $\Phi$, the definition \eqref{eq:definition_delta}, and the fact that $\langle\Phi(\boldsymbol{\mathcal{V}}_k),\boldsymbol{\mathcal{V}_k}\rangle\geq0$ for all $\boldsymbol{\mathcal{V}_k}\in\mathbb{X}_k$ such that $\|\boldsymbol{\mathcal{V}}_k\|_{\mathbb{X}}=\delta$, and by invoking a consequence of Brouwer's fixed-point theorem \cite[Chapter IV, Corollary 1.1]{MR851383}, we obtain the existence of $\boldsymbol{\mathcal{U}}_k$ such that
\begin{equation}\label{eq:Galerkin_approximation}
\langle\Phi(\boldsymbol{\mathcal{U}}_k),\boldsymbol{\mathcal{V}}_k\rangle=0
\quad
\forall\boldsymbol{\mathcal{V}}_k\in\mathbb{X}_k,
\qquad
\|(\mathbf{u}_k,T_k,C_k)\|_{\mathbb{X}}\leq\delta.
\end{equation}

\emph{Step 3.} \emph{(Limit passage).} First, we observe that the sequences $\{\mathbf{u}_k\}_{k\in\mathbb{N}}$, $\{T_k\}_{k\in\mathbb{N}}$, and $\{C_k\}_{k\in\mathbb{N}}$ are uniformly bounded in $\mathbf{H}_0^1(\Omega)$, $H_0^1(\Omega)$, and $H_0^1(\Omega)$, respectively. Thus, we deduce the existence of nonrelabeled subsequences such that $\mathbf{u}_k\rightharpoonup\mathbf{u}$ in $\mathbf{H}_0^1(\Omega)$, and $T_k\rightharpoonup T$ and $C_k\rightharpoonup C$ in $H_0^1(\Omega)$, as $k\uparrow\infty$. Second, we invoke \cite[Theorem 6.3, Part I]{MR2424078} to guarantee that $\mathbf{u}_k\to\mathbf{u}$ in $\mathbf{L}^q(\Omega)$, and $T_k\to T$ and $C_k\to C$ in $L^q(\Omega)$, as $k\uparrow\infty$, for $q<\infty$ in two dimensions and $q<6$ in three dimensions. In what follows, we write $\boldsymbol{\mathcal{V}}_\ell=(\mathbf{v}_\ell,S_\ell,R_\ell)$. Therefore, for every $\ell\leq k$ we have that
\begin{equation}\label{eq:Galerkin_approximation_ell}
\boldsymbol{\mathcal{U}}_k=(\mathbf{u}_k,T_k,C_k)\in\mathbb{X}_k:
\quad
\langle\Phi(\boldsymbol{\mathcal{U}}_k),\boldsymbol{\mathcal{V}}_\ell\rangle=0
\quad
\forall\boldsymbol{\mathcal{V}}_\ell\in\mathbb{X}_\ell.
\end{equation}
In what follows, we show that $\boldsymbol{\mathcal{U}}=(\mathbf{u},T,C)$ satisfies $\langle\Phi(\boldsymbol{\mathcal{U}}),\boldsymbol{\mathcal{V}}\rangle=0$ for all $\boldsymbol{\mathcal{V}}\in\mathbf{V}(\Omega)\times H_0^1(\Omega)\times H_0^1(\Omega)$. To do so, we consider:

 (I) \emph{Preliminary convergences:} We observe that, since $\nu$ is continuous and uniformly bounded, and given the fact that $|\nu(T_k,C_k)\boldsymbol{\tau}|\leq\nu_{+}|\boldsymbol{\tau}|$, we can apply the Lebesgue dominated convergence theorem to obtain $\nu(T_k,C_k)\mathbf{w}\to\nu(T,C)\mathbf{w}$ in $\mathbf{L}^2(\Omega)$ as $k\uparrow\infty$ for every $\mathbf{w}\in\mathbf{L}^2(\Omega)$. On the other hand, using similar arguments, we can obtain the same results for $\kappa(T_k,C_k)$ and $\mathcal{D}(T_k,C_k)$. Moreover, the Lipschitz continuity of $\mathbf{f}_1$ and $\mathbf{f}_2$, and Sobolev embeddings yield $\|\mathbf{f}_1(T_k)-\mathbf{f}_1(T)\|_{\mathbf{L}^2(\Omega)}\leq C_{\mathbf{f}_1}\|T_k-T\|_{L^2(\Omega)}\to0$ and $\|\mathbf{f}_2(C_k)-\mathbf{f}_2(C)\|_{\mathbf{L}^2(\Omega)}\leq C_{\mathbf{f}_2}\|C_k-C\|_{L^2(\Omega)}\to0$ as $k\uparrow\infty$.

(II) \emph{Strong convergence $\mathbf{u}_k\to\mathbf{u}$ in $\mathbf{H}_0^1(\Omega)$ as $k\uparrow\infty$:} In view of \eqref{eq:density}, there exists a sequence $\{\mathbf{w}_k\}_{k\in\mathbb{N}}$, with $\mathbf{w}_k\in\mathbb{V}_k$, such that $\mathbf{w}_k\to\mathbf{u}$ in $\mathbf{H}_0^1(\Omega)$ as $k\uparrow\infty$. Define $\mathbf{e}_k:=\mathbf{u}_k-\mathbf{w}_k\in\mathbb{V}_k$. We set $\boldsymbol{\mathcal{V}}_k=(\mathbf{e}_k,0,0)$ in \eqref{eq:Galerkin_approximation}. This, combined with \eqref{eq:prop_02} for $\mathcal{A}_{T_k,C_k}$, yields
$\nu_{-}\varepsilon_{-}\|\nabla\mathbf{e}_k\|_{\mathbf{L}^2(\Omega)}^2
\leq
\langle\mathcal{A}_{T_k,C_k}(\mathbf{u}_k)-\mathcal{A}_{T_k,C_k}(\mathbf{w}_k),\mathbf{e}_k\rangle
=
\langle\mathbf{f}(T_k,C_k),\mathbf{e}_k\rangle
-\langle\mathcal{A}_{T_k,C_k}(\mathbf{w}_k),\mathbf{e}_k\rangle=:\mathrm{I}_k+\mathrm{II}_k$.
Since $\mathbf{e}_k\rightharpoonup\mathbf{0}$ in $\mathbf{H}_0^1(\Omega)$, $\mathbf{e}_k\to\mathbf{0}$ in $\mathbf{L}^2(\Omega)$, and the sequence $\{\mathbf{f}_1(T_k)+\mathbf{f}_2(C_k)\}_{k\in\mathbb{N}}$ is bounded in $\mathbf{L}^2(\Omega)$, we have that the term $\mathrm{I}_k\to0$ as $k\uparrow\infty$. On the other hand, we invoke assumption \eqref{assumption_1} and the arguments presented in (I) to obtain that
\begin{multline*}
\|\nu(T_k,C_k)\,\varepsilon(x,|\nabla\mathbf{w}_k|)\nabla\mathbf{w}_k-\nu(T,C)\,\varepsilon(x,|\nabla\mathbf{u}|)\nabla\mathbf{u}\|_{\mathbf{L}^2(\Omega)}
\\
\leq
\nu_{+}\varepsilon_{+}\|\nabla(\mathbf{w}_k-\mathbf{u})\|_{\mathbf{L}^2(\Omega)}
+\|(\nu(T_k,C_k)-\nu(T,C))\,\varepsilon(x,|\nabla\mathbf{u}|)\nabla\mathbf{u}\|_{\mathbf{L}^2(\Omega)}
\to0
\end{multline*}
as $k\uparrow\infty$. This and the weak convergence $\nabla\mathbf{e}_k\rightharpoonup\mathbf{0}$ in $\mathbf{L}^2(\Omega)$ show that the term $\mathrm{II}_k\to0$ as $k\uparrow\infty$. Therefore, $\nabla\mathbf{e}_k\to\mathbf{0}$ in $\mathbf{L}^2(\Omega)$ and, since $\mathbf{w}_k\to\mathbf{u}$ in $\mathbf{H}_0^1(\Omega)$, we conclude that $\mathbf{u}_k\to\mathbf{u}$ in $\mathbf{H}_0^1(\Omega)$ as $k\uparrow\infty$.

(III) \emph{The operator $\mathcal{A}$ satisfies $\langle\mathcal{A}_{T_k,C_k}(\mathbf{u}_k),\mathbf{v}\rangle\to\langle\mathcal{A}_{T,C}(\mathbf{u}),\mathbf{v}\rangle$ as $k\uparrow\infty$:} We invoke \eqref{eq:prop_01} for $\mathcal{A}_{T_k,C_k}$, \eqref{assumption_1}, and arguments presented in (I) and (II). Thus,
\begin{multline*}
|\langle\mathcal{A}_{T_k,C_k}(\mathbf{u}_k)-\mathcal{A}_{T,C}(\mathbf{u}),\mathbf{v}\rangle|
\lesssim \|\nabla(\mathbf{u}_k-\mathbf{u})\|_{\mathbf{L}^2(\Omega)}\|\nabla\mathbf{v}\|_{\mathbf{L}^2(\Omega)}
\\
+\|(\nu(T_k,C_k)-\nu(T,C))\,\varepsilon(x,|\nabla\mathbf{u}|)\nabla\mathbf{u}\|_{\mathbf{L}^2(\Omega)}
\|\nabla\mathbf{v}\|_{\mathbf{L}^2(\Omega)}\to 0,\quad\text{as}\quad k\uparrow\infty.
\end{multline*}

(IV) \emph{The right-hand side satisfies $\langle\mathbf{f}(T_k,C_k),\mathbf{v}\rangle\to\langle\mathbf{f}(T,C),\mathbf{v}\rangle$ as $k\uparrow\infty$:} This result is an immediate consequence of item (I).

(V) \emph{The diffusion terms:} In view of item (I), the fact that $\kappa(T_k,C_k)\nabla S\to\kappa(T,C)\nabla S$ in $\mathbf{L}^2(\Omega)$, together with the weak convergence $\nabla T_k\rightharpoonup\nabla T$ in $\mathbf{L}^2(\Omega)$ yields $\int_\Omega\kappa(T_k,C_k)\nabla T_k\cdot\nabla S\,\mathrm{d}x\to\int_\Omega\kappa(T,C)\nabla T\cdot\nabla S\,\mathrm{d}x$ as $k\uparrow\infty$. The same arguments show that $\int_\Omega\mathcal{D}(T_k,C_k)\nabla C_k\cdot\nabla R\,\mathrm{d}x\to\int_\Omega\mathcal{D}(T,C)\nabla C\cdot\nabla R\,\mathrm{d}x$ as $k\uparrow\infty$.

(VI) \emph{The convective term satisfies $\int_\Omega(\mathbf{u}_k\cdot\nabla T_k)S\,\mathrm{d}x\to\int_\Omega(\mathbf{u}\cdot\nabla T)S\,\mathrm{d}x$ as $k\uparrow\infty$:} H\"older's inequality gives $\|(\mathbf{u}_k-\mathbf{u})S\|_{\mathbf{L}^2(\Omega)}\leq\|\mathbf{u}_k-\mathbf{u}\|_{\mathbf{L}^4(\Omega)}\|S\|_{L^4(\Omega)}$, and thus $\mathbf{u}_kS\to\mathbf{u}S$ in $\mathbf{L}^2(\Omega)$. These arguments and the weak convergence $\nabla T_k\rightharpoonup\nabla T$ in $\mathbf{L}^2(\Omega)$ yield the desired result. Similarly, $\int_\Omega(\mathbf{u}_k\cdot\nabla C_k)R\,\mathrm{d}x\to\int_\Omega(\mathbf{u}\cdot\nabla C)R\,\mathrm{d}x$ as $k\uparrow\infty$.

The results obtained in (III)--(VI) allow us to pass to the limit in \eqref{eq:Galerkin_approximation_ell} and to conclude that $\langle\Phi(\boldsymbol{\mathcal{U}}),\boldsymbol{\mathcal{V}}_\ell\rangle=0$ for all $\boldsymbol{\mathcal{V}}_\ell\in\mathbb{X}_\ell$. Then, a density argument shows that $\langle\Phi(\boldsymbol{\mathcal{U}}),\boldsymbol{\mathcal{V}}\rangle=0$ for all $\boldsymbol{\mathcal{V}}\in\mathbf{V}(\Omega)\times H_0^1(\Omega)\times H_0^1(\Omega)$. On the other hand, we can conclude that $(\mathbf{u},T,C)$ satisfies the second, third, and fourth equations of \eqref{eq:weak_coupled} and
\[
\langle\mathcal{H},\mathbf{v}\rangle
:=
\langle\mathcal{A}_{T,C}(\mathbf{u})-\mathbf{f}(T,C),\mathbf{v}\rangle
=0
\quad
\forall\mathbf{v}\in\mathbf{V}(\Omega).
\]

\emph{Step 4.} \emph{The pressure.} We observe that the functional $\mathcal{H}$ is linear and continuous on $\mathbf{H}_0^1(\Omega)$. In addition, $\langle\mathcal{H},\mathbf{v}\rangle=0$ for all $\mathbf{v}\in\mathbf{V}(\Omega)$. With these elements, we can apply de Rham's theorem \cite[Theorem B.73]{Guermond-Ern} which guarantees the existence of $\mathsf{p}\in L_0^2(\Omega)$ such that $\langle\mathcal{H},\mathbf{v}\rangle=\int_\Omega\mathsf{p}\operatorname{div}\mathbf{v}\,\mathrm{d}x$ for all $\mathbf{v}\in\mathbf{H}_0^1(\Omega)$. Thus, we have proved the existence of a solution $(\mathbf{u},\mathsf{p},T,C)\in\mathbf{H}_0^1(\Omega)\times L_0^2(\Omega)\times H_0^1(\Omega)\times H_0^1(\Omega)$ to the coupled problem \eqref{eq:weak_coupled}.

\emph{Step 5.} \emph{Stability bounds:} Let $(\mathbf{u},\mathsf{p},T,C)\in\mathbf{H}_0^1(\Omega)\times L_0^2(\Omega)\times H_0^1(\Omega)\times H_0^1(\Omega)$ be a solution of problem \eqref{eq:weak_coupled}. Taking $S=T$ in the temperature equation of \eqref{eq:weak_coupled}, and using that $\int_\Omega(\mathbf{u}\cdot\nabla T)T\,\mathrm{d}x=0$, since $\mathbf{u}\in\mathbf{V}(\Omega)$, we immediately obtain the estimate $\kappa_{-}\|\nabla T\|_{\mathbf{L}^2(\Omega)}^2\leq\|\mathsf{g}_1\|_{H^{-1}(\Omega)}\|\nabla T\|_{\mathbf{L}^2(\Omega)}$.
The bound for $C$ in \eqref{eq:estimate_TC_coupled} follows by similar arguments. We now set $\mathbf{v}=\mathbf{u}$ in the first equation of \eqref{eq:weak_coupled} and invoke \eqref{eq:prop_02} with $\mathbf{w}=\mathbf{0}$ for $\mathcal{A}_{T,C}$, and the estimate $|\mathbf{f}_i(r)|\leq C_{\mathbf{f}_i}|r|$ with $i=1,2$ in order to obtain
\[
\nu_{-}\varepsilon_{-}\|\nabla\mathbf{u}\|_{\mathbf{L}^2(\Omega)}
\leq
\|\mathbf{f}_0\|_{\mathbf{H}^{-1}(\Omega)}
+\mathcal{C}_{2\hookrightarrow2}^2
\left(C_{\mathbf{f}_1}\|\nabla T\|_{\mathbf{L}^2(\Omega)}+C_{\mathbf{f}_2}\|\nabla C\|_{\mathbf{L}^2(\Omega)}\right).
\]
This bound and \eqref{eq:estimate_TC_coupled} show the estimate for $\mathbf{u}$ in \eqref{eq:estimate_up_coupled}. Finally, the inf-sup condition \eqref{eq:infsup_cont}, the first equation of \eqref{eq:weak_coupled}, assumption \eqref{assumption_1} with $\boldsymbol{\omega}=\mathbf{0}$, \eqref{eq:bound_f_TC} and \eqref{eq:estimate_TC_coupled} produce
$\beta\|\mathsf{p}\|_{L^2(\Omega)}
\leq
\nu_{+}\varepsilon_{+}\|\nabla\mathbf{u}\|_{\mathbf{L}^2(\Omega)}+\|\mathbf{f}_0\|_{\mathbf{H}^{-1}(\Omega)}
+
\mathcal{C}_{2\hookrightarrow2}^2
\big(
C_{\mathbf{f}_1}\kappa_{-}^{-1}\|\mathsf{g}_1\|_{H^{-1}(\Omega)}
+
C_{\mathbf{f}_2}\mathcal{D}_{-}^{-1}\|\mathsf{g}_2\|_{H^{-1}(\Omega)}
\big)$. Thus, we have the estimate for the pressure $\mathsf{p}$.
\end{proof}

\subsection{Uniqueness of solutions} In the following, we will introduce regularity assumptions regarding the solutions in order to establish the uniqueness of the solution to the nonlinear system \eqref{eq:weak_coupled}.

\begin{theorem}[uniqueness of solutions]\label{th:uniqueness_coupled_3D}
Assume that $d=3$ and that the data satisfy the assumptions of Theorem \ref{th:existence_coupled}. Let $(\mathbf{u}_1,\mathsf{p}_1,T_1,C_1)\in[\mathbf{W}^{1,3}(\Omega)\cap\mathbf{H}_0^1(\Omega)]\times L_0^2(\Omega)\times[W^{1,3}(\Omega)\cap H_0^1(\Omega)]\times[W^{1,3}(\Omega)\cap H_0^1(\Omega)]$ solve \eqref{eq:weak_coupled} and satisfy
\begin{align}
\label{eq:small_data_coupled}
\frac{C_{\mathbf{f}}\,\mathcal{C}_{2\hookrightarrow2}^2\mathcal{C}_{4\hookrightarrow2}^2\,\Gamma(\mathsf{g}_1,\mathsf{g}_2)}{\nu_{-}\varepsilon_{-}(1-\theta)}
+\frac{\varepsilon_{+}\mathcal{L}_{\nu}\mathcal{C}_{6\hookrightarrow2}\mathcal{C}_{4\hookrightarrow2}^2\,\Gamma(\mathsf{g}_1,\mathsf{g}_2)\,\|\nabla\mathbf{u}_1\|_{\mathbf{L}^3(\Omega)}}{\nu_{-}\varepsilon_{-}(1-\theta)}
&<1,
\\
\label{eq:small_data_coupled_2}
\theta:=\frac{\mathcal{L}_{\kappa}\mathcal{C}_{6\hookrightarrow2}\|\nabla T_1\|_{\mathbf{L}^3(\Omega)}}{\kappa_{-}}
+\frac{\mathcal{L}_{\mathcal{D}}\mathcal{C}_{6\hookrightarrow2}\|\nabla C_1\|_{\mathbf{L}^3(\Omega)}}{\mathcal{D}_{-}}
&<1,
\end{align}
where $\Gamma(\mathsf{g}_1,\mathsf{g}_2):=\kappa_{-}^{-2}\|\mathsf{g}_1\|_{H^{-1}(\Omega)}+\mathcal{D}_{-}^{-2}\|\mathsf{g}_2\|_{H^{-1}(\Omega)}$ and $C_{\mathbf{f}}:=\max\{C_{\mathbf{f}_1},C_{\mathbf{f}_2}\}$. Then, problem \eqref{eq:weak_coupled} admits exactly one solution in $\mathbf{H}_0^1(\Omega)\times L_0^2(\Omega)\times H_0^1(\Omega)\times H_0^1(\Omega)$, namely $(\mathbf{u}_1,\mathsf{p}_1,T_1,C_1)$. Here, $\mathcal{C}_{4\hookrightarrow2}$, $\mathcal{C}_{2\hookrightarrow2}$, and $\mathcal{C}_{6\hookrightarrow2}$ are the best constants in the Sobolev embeddings $H_0^1(\Omega)\hookrightarrow L^4(\Omega)$, $H_0^1(\Omega)\hookrightarrow L^2(\Omega)$, and $H_0^1(\Omega)\hookrightarrow L^6(\Omega)$, respectively.
\end{theorem}

\begin{proof}
Assume that $(\mathbf{u}_2,\mathsf{p}_2,T_2,C_2)\in\mathbf{H}_0^1(\Omega)\times L_0^2(\Omega)\times H_0^1(\Omega)\times H_0^1(\Omega)$ also solves \eqref{eq:weak_coupled}. In what follows, we introduce the notation $\mathbf{e}_{\mathbf{u}}:=\mathbf{u}_1-\mathbf{u}_2$, $e_{\mathsf{p}}:=\mathsf{p}_1-\mathsf{p}_2$, $e_T:=T_1-T_2$, and $e_C:=C_1-C_2$. We subtract the temperature equations satisfied by both solutions to obtain, for all $S\in H_0^1(\Omega)$ the identity
$\int_\Omega\kappa(T_2,C_2)\nabla e_T\cdot\nabla S\,\mathrm{d}x
+\int_\Omega(\mathbf{u}_1\cdot\nabla T_1-\mathbf{u}_2\cdot\nabla T_2)S\,\mathrm{d}x
=\int_\Omega(\kappa(T_2,C_2)-\kappa(T_1,C_1))\nabla T_1\cdot\nabla S\,\mathrm{d}x$.
If we set $S=e_T$, use the fact that $\mathbf{u}_1\cdot\nabla T_1-\mathbf{u}_2\cdot\nabla T_2=\mathbf{e}_{\mathbf{u}}\cdot\nabla T_1+\mathbf{u}_2\cdot\nabla e_T$, and note that $\int_{\Omega}(\mathbf{u}_2\cdot\nabla e_T)e_T\, \mathrm{d}x=0$, since $\mathbf{u}_2\in\mathbf{V}(\Omega)$, we have that
\[
\int_\Omega\kappa(T_2,C_2)|\nabla e_T|^2\,\mathrm{d}x
=-\int_\Omega(\mathbf{e}_{\mathbf{u}}\cdot\nabla T_1)e_T\,\mathrm{d}x
+\int_\Omega(\kappa(T_2,C_2)-\kappa(T_1,C_1))\nabla T_1\cdot\nabla e_T\,\mathrm{d}x.
\]
We now apply the bound $\kappa(T_2,C_2)\geq\kappa_{-}$, the Lipschitz continuity of $\kappa$, suitable H\"older's inequalities, the Sobolev embeddings $H_0^1(\Omega)\hookrightarrow L^4(\Omega)$ and $H_0^1(\Omega)\hookrightarrow L^6(\Omega)$, and the bound $\|\nabla T_1\|_{\mathbf{L}^2(\Omega)}\leq\kappa_{-}^{-1}\|\mathsf{g}_1\|_{H^{-1}(\Omega)}$ from \eqref{eq:estimate_TC_coupled} to arrive at
\begin{multline*}
\kappa_{-}\|\nabla e_T\|_{\mathbf{L}^2(\Omega)}^2
\leq
\kappa_{-}^{-1}\mathcal{C}_{4\hookrightarrow2}^2\|\mathsf{g}_1\|_{H^{-1}(\Omega)}\|\nabla\mathbf{e}_{\mathbf{u}}\|_{\mathbf{L}^2(\Omega)}\|\nabla e_T\|_{\mathbf{L}^2(\Omega)}
\\
+\mathcal{L}_{\kappa}\mathcal{C}_{6\hookrightarrow2}\|\nabla T_1\|_{\mathbf{L}^3(\Omega)}\left(\|\nabla e_T\|_{\mathbf{L}^2(\Omega)}+\|\nabla e_C\|_{\mathbf{L}^2(\Omega)}\right)\|\nabla e_T\|_{\mathbf{L}^2(\Omega)}.
\end{multline*}
This estimate immediately leads to the following bound:
\begin{multline*}
\|\nabla e_T\|_{\mathbf{L}^2(\Omega)}
\leq
\kappa_{-}^{-2}\mathcal{C}_{4\hookrightarrow2}^2\|\mathsf{g}_1\|_{H^{-1}(\Omega)}\|\nabla\mathbf{e}_{\mathbf{u}}\|_{\mathbf{L}^2(\Omega)}
\\
+\kappa_{-}^{-1}\mathcal{L}_{\kappa}\mathcal{C}_{6\hookrightarrow2}\|\nabla T_1\|_{\mathbf{L}^3(\Omega)}\left(\|\nabla e_T\|_{\mathbf{L}^2(\Omega)}+\|\nabla e_C\|_{\mathbf{L}^2(\Omega)}\right).
\end{multline*}
Similar arguments, applied to the concentration equation with $R=e_C$, show that
\begin{multline*}
\|\nabla e_C\|_{\mathbf{L}^2(\Omega)}
\leq
\mathcal{D}_{-}^{-2}\mathcal{C}_{4\hookrightarrow2}^2\|\mathsf{g}_2\|_{H^{-1}(\Omega)}\|\nabla\mathbf{e}_{\mathbf{u}}\|_{\mathbf{L}^2(\Omega)}
\\
+\mathcal{D}_{-}^{-1}\mathcal{L}_{\mathcal{D}}\mathcal{C}_{6\hookrightarrow2}\|\nabla C_1\|_{\mathbf{L}^3(\Omega)}\left(\|\nabla e_T\|_{\mathbf{L}^2(\Omega)}+\|\nabla e_C\|_{\mathbf{L}^2(\Omega)}\right).
\end{multline*}
Adding the last two inequalities and using the definitions of $\theta$ in \eqref{eq:small_data_coupled_2} and of $\Gamma(\mathsf{g}_1,\mathsf{g}_2)$, we obtain that
\begin{equation}\label{eq:bound_TC_uniq}
(1-\theta)\left(\|\nabla e_T\|_{\mathbf{L}^2(\Omega)}+\|\nabla e_C\|_{\mathbf{L}^2(\Omega)}\right)
\leq
\mathcal{C}_{4\hookrightarrow2}^2\,\Gamma(\mathsf{g}_1,\mathsf{g}_2)\,\|\nabla\mathbf{e}_{\mathbf{u}}\|_{\mathbf{L}^2(\Omega)}.
\end{equation}

We now estimate $\|\nabla\mathbf{e}_{\mathbf{u}}\|_{\mathbf{L}^2(\Omega)}$. To this end, we subtract the first equations in \eqref{eq:weak_coupled} with test functions $\mathbf{v}\mathrel{\in}\mathbf{V}(\Omega)$. On the other hand, we use the map $\mathcal{A}_{T,C}$ defined in \eqref{eq:operator_ATC}, and we write $\langle\mathcal{A}_{T_1,C_1}(\mathbf{u}_1),\mathbf{v}\rangle=\langle\mathcal{A}_{T_2,C_2}(\mathbf{u}_1)+(\mathcal{A}_{T_1,C_1}(\mathbf{u}_1)-\mathcal{A}_{T_2,C_2}(\mathbf{u}_1)),\mathbf{v}\rangle$. Thus, we obtain that
\begin{multline*}
\langle\mathcal{A}_{T_2,C_2}(\mathbf{u}_1)-\mathcal{A}_{T_2,C_2}(\mathbf{u}_2),\mathbf{v}\rangle
=
\langle\mathbf{f}_1(T_1)-\mathbf{f}_1(T_2)+\mathbf{f}_2(C_1)-\mathbf{f}_2(C_2),\mathbf{v}\rangle
\\
+\int_\Omega(\nu(T_2,C_2)-\nu(T_1,C_1))\,\varepsilon(x,|\nabla\mathbf{u}_1|)\nabla\mathbf{u}_1:\nabla\mathbf{v}\,\mathrm{d}x
\end{multline*}
for all $\mathbf{v}\in\mathbf{V}(\Omega)$. Let $\mathbf{v}=\mathbf{e}_{\mathbf{u}}$. In view of \eqref{eq:prop_02} for $\mathcal{A}_{T_2,C_2}$, the left-hand side is bounded below by $\nu_{-}\varepsilon_{-}\|\nabla\mathbf{e}_{\mathbf{u}}\|_{\mathbf{L}^2(\Omega)}^2$. To bound the right-hand side, we use the Lipschitz continuity of $\mathbf{f}_1$, $\mathbf{f}_2$, and $\nu$, the estimate $|\varepsilon(x,|\boldsymbol{\tau}|)\boldsymbol{\tau}|\leq\varepsilon_{+}|\boldsymbol{\tau}|$, H\"older's inequality, and the embeddings $H_0^1(\Omega)\hookrightarrow L^2(\Omega)$ and $H_0^1(\Omega)\hookrightarrow L^6(\Omega)$ to produce that
\begin{multline*}
\nu_{-}\varepsilon_{-}\|\nabla\mathbf{e}_{\mathbf{u}}\|_{\mathbf{L}^2(\Omega)}^2
\leq
\left(C_{\mathbf{f}}\,\mathcal{C}_{2\hookrightarrow2}^2+\varepsilon_{+}\mathcal{L}_{\nu}\mathcal{C}_{6\hookrightarrow2}\|\nabla\mathbf{u}_1\|_{\mathbf{L}^3(\Omega)}\right)
\\
\cdot\left(\|\nabla e_T\|_{\mathbf{L}^2(\Omega)}+\|\nabla e_C\|_{\mathbf{L}^2(\Omega)}\right)\|\nabla\mathbf{e}_{\mathbf{u}}\|_{\mathbf{L}^2(\Omega)}.
\end{multline*}
Since $\theta<1$, see \eqref{eq:small_data_coupled_2}, we can use \eqref{eq:bound_TC_uniq} in the last inequality. This yields
\[
\|\nabla\mathbf{e}_{\mathbf{u}}\|_{\mathbf{L}^2(\Omega)}^2
\!\left(1
\!-\!\frac{C_{\mathbf{f}}\,\mathcal{C}_{2\hookrightarrow2}^2\mathcal{C}_{4\hookrightarrow2}^2\,\Gamma(\mathsf{g}_1,\mathsf{g}_2)}{\nu_{-}\varepsilon_{-}(1-\theta)}
\!-\!\frac{\varepsilon_{+}\mathcal{L}_{\nu}\mathcal{C}_{6\hookrightarrow2}\mathcal{C}_{4\hookrightarrow2}^2\,\Gamma(\mathsf{g}_1,\mathsf{g}_2)\,\|\nabla\mathbf{u}_1\|_{\mathbf{L}^3(\Omega)}}{\nu_{-}\varepsilon_{-}(1-\theta)}\right)
\!\leq \!0.
\]
In view of \eqref{eq:small_data_coupled}, the term that multiplies $\|\nabla\mathbf{e}_{\mathbf{u}}\|_{\mathbf{L}^2(\Omega)}^2$ is positive, and therefore $\mathbf{e}_{\mathbf{u}}=\mathbf{0}$, i.e., $\mathbf{u}_1=\mathbf{u}_2$. Then, \eqref{eq:bound_TC_uniq} with $\mathbf{e}_{\mathbf{u}}=\mathbf{0}$ and $\theta<1$ imply that $e_T=e_C=0$. Finally, since $(\mathbf{u}_1,T_1,C_1)=(\mathbf{u}_2,T_2,C_2)$, the first equation in \eqref{eq:weak_coupled} yields $\int_\Omega e_{\mathsf{p}}\operatorname{div}\mathbf{v}\,\mathrm{d}x=0$ for every $\mathbf{v}\in\mathbf{H}_0^1(\Omega)$. The inf-sup condition \eqref{eq:infsup_cont} thus guarantees that $e_{\mathsf{p}}=0$.
\end{proof}

Since $H_0^1(\Omega)\hookrightarrow L^{\sigma}(\Omega)$ for every $\sigma<\infty$ when $d=2$, the regularity assumed in Theorem~\ref{th:uniqueness_coupled_3D} can be reduced in two dimensions.

\begin{theorem}[uniqueness of solutions for $d=2$]\label{th:uniqueness_coupled_2D}
Assume that $d=2$ and that the data satisfy the assumptions of Theorem~\ref{th:existence_coupled}. Let $\ell>2$, and let $\sigma$ be such that $1/\sigma+1/\ell=1/2$. Let $(\mathbf{u}_1,\mathsf{p}_1,T_1,C_1)\in[\mathbf{W}^{1,\ell}(\Omega)\cap\mathbf{H}_0^1(\Omega)]\times L_0^2(\Omega)\times[W^{1,\ell}(\Omega)\cap H_0^1(\Omega)]\times[W^{1,\ell}(\Omega)\cap H_0^1(\Omega)]$ solve \eqref{eq:weak_coupled} and satisfy
\begin{align}
\label{eq:small_data_coupled_2D}
\frac{C_{\mathbf{f}}\,\mathcal{C}_{2\hookrightarrow2}^2\mathcal{C}_{4\hookrightarrow2}^2\,\Gamma(\mathsf{g}_1,\mathsf{g}_2)}{\nu_{-}\varepsilon_{-}(1-\theta_{\ell})}
+\frac{\varepsilon_{+}\mathcal{L}_{\nu}\mathcal{C}_{\sigma\hookrightarrow2}\mathcal{C}_{4\hookrightarrow2}^2\,\Gamma(\mathsf{g}_1,\mathsf{g}_2)\,\|\nabla\mathbf{u}_1\|_{\mathbf{L}^{\ell}(\Omega)}}{\nu_{-}\varepsilon_{-}(1-\theta_{\ell})}
&<1,
\\
\label{eq:small_data_coupled_2D_2}
\theta_{\ell}:=\frac{\mathcal{L}_{\kappa}\mathcal{C}_{\sigma\hookrightarrow2}\|\nabla T_1\|_{\mathbf{L}^{\ell}(\Omega)}}{\kappa_{-}}
+\frac{\mathcal{L}_{\mathcal{D}}\mathcal{C}_{\sigma\hookrightarrow2}\|\nabla C_1\|_{\mathbf{L}^{\ell}(\Omega)}}{\mathcal{D}_{-}}
&<1,
\end{align}
where $\Gamma(\mathsf{g}_1,\mathsf{g}_2)$ and $C_{\mathbf{f}}$ are as in Theorem~\ref{th:uniqueness_coupled_3D}. Then, problem \eqref{eq:weak_coupled} admits exactly one solution in $\mathbf{H}_0^1(\Omega)\times L_0^2(\Omega)\times H_0^1(\Omega)\times H_0^1(\Omega)$, namely $(\mathbf{u}_1,\mathsf{p}_1,T_1,C_1)$. Here, $\mathcal{C}_{4\hookrightarrow2}$, $\mathcal{C}_{2\hookrightarrow2}$, and $\mathcal{C}_{\sigma\hookrightarrow2}$ are the best constants in the Sobolev embeddings $H_0^1(\Omega)\hookrightarrow L^4(\Omega)$, $H_0^1(\Omega)\hookrightarrow L^2(\Omega)$, and $H_0^1(\Omega)\hookrightarrow L^{\sigma}(\Omega)$, respectively.
\end{theorem}

\begin{proof}
The proof is a suitable adaptation of the proof of Theorem~\ref{th:uniqueness_coupled_3D}.
\end{proof}

\section{Finite element analysis} In this section, we present a finite element discretization of problem \eqref{eq:model}. In particular, we study the convergence of the discrete scheme and we derive a priori error bounds. We begin with some standard elements of the finite element method; see \cite{MR3097958,MR2373954,CiarletBook,Guermond-Ern}.
\label{sec:fem}

\subsection{Framework}
\label{sec:fem_notation}
Let $\Omega$ be a Lipschitz polytope. We denote by $\mathscr{T}_h=\{K\}$ a conforming partition of $\bar\Omega$ into closed simplices $K$ with size $h_K=\operatorname{diam}(K)$, and $h:=\max\{h_K:K\in\mathscr{T}_h\}$. We consider a family $\mathbb{T}=\{\mathscr{T}_h\}_{h>0}$ of such triangulations, obtained by refining an initial mesh $\mathscr{T}_0$, and we assume that $\mathbb{T}$ is shape regular.

\subsection{Finite element spaces}
\label{sec:fem_spaces}
Given $\mathscr{T}_h\in \mathbb{T}$, we approximate the velocity field $\mathbf{u}$ and the pressure $\mathsf{p}$ of the fluid with a pair $(\mathbf{X}_h,M_h)$ of finite element spaces satisfying the classical uniform discrete inf-sup condition, i.e., there is $\tilde\beta>0$, which does not depend on $h$, such that
\begin{equation}\label{eq:infsup_disc}
\inf_{0\neq\mathsf{q}_h\in M_h}\sup_{\mathbf{0}\neq\mathbf{v}_h\in\mathbf{X}_h}\frac{\int_\Omega\mathsf{q}_h\operatorname{div}\mathbf{v}_h\,\mathrm{d}x}{\|\nabla\mathbf{v}_h\|_{\mathbf{L}^{2}(\Omega)}\|\mathsf{q}_h\|_{L^2(\Omega)}}\geq\tilde\beta>0.
\end{equation}
In this work, we consider two classical pairs:

$\checkmark$ The \emph{Taylor--Hood element} of lowest order, proposed in \cite{TaylorHood1973} in two dimensions (see also \cite[Chapter II, Section 4.2]{MR851383}, \cite[Section 4.2.5]{Guermond-Ern}, and \cite[Section 8.8.2]{MR3097958}):
\begin{align}\label{eq:TH_X}
\mathbf{X}_h&=\{\mathbf{v}_h\in\mathbf{C}(\bar\Omega):\mathbf{v}_h|_K\in[\mathbb{P}_2(K)]^d\ \forall K\in\mathscr{T}_h\}\cap\mathbf{H}_0^1(\Omega),
\\
\label{eq:TH_M}
M_h&=\{\mathsf{q}_h\in L_0^2(\Omega)\cap C(\bar\Omega):\mathsf{q}_h|_K\in\mathbb{P}_1(K)\ \forall K\in\mathscr{T}_h\}.
\end{align}
$\checkmark$ The \emph{mini element}, proposed in \cite{MR0799997} in two dimensions (see also \cite[Chapter II, Section 4.1]{MR851383}, \cite[Section 4.2.4]{Guermond-Ern}, and \cite[Section 8.4.2]{MR3097958}):
\begin{align}\label{eq:mini_X}
\mathbf{X}_h&=\{\mathbf{v}_h\in\mathbf{C}(\bar\Omega):\mathbf{v}_h|_K\in[\mathbb{W}(K)]^d\ \forall K\in\mathscr{T}_h\}\cap\mathbf{H}_0^1(\Omega),
\\
\label{eq:mini_M}
M_h&=\{\mathsf{q}_h\in L_0^2(\Omega)\cap C(\bar\Omega):\mathsf{q}_h|_K\in\mathbb{P}_1(K)\ \forall K\in\mathscr{T}_h\},
\end{align}
where $\mathbb{W}(K):=\mathbb{P}_1(K)\oplus\mathbb{B}(K)$ and $\mathbb{B}(K)$ is the space generated by the bubble function associated with $K$.
Both pairs satisfy \eqref{eq:infsup_disc}. For the mini element, we refer to \cite[Lemma 4.20]{Guermond-Ern}. For the Taylor--Hood element, we refer to \cite[Theorems 8.8.1 and 8.8.2]{MR3097958}.

To approximate the temperature and the concentration, we use the space
\begin{equation}\label{eq:Y_h}
Y_h:=\{S_h\in C(\bar\Omega):S_h|_K\in\mathbb{P}_{\mathfrak{r}}(K)\ \forall K\in\mathscr{T}_h\}\cap H_0^1(\Omega),
\end{equation}
with $\mathfrak{r}=1$ when $(\mathbf{u},\mathsf{p})$ is approximated with the mini element and $\mathfrak{r}=2$ when the Taylor--Hood element is used.

Moreover, the spaces $\mathbf{X}_h$, $M_h$, and $Y_h$ approximate $\mathbf{H}_0^1(\Omega)$, $L_0^2(\Omega)$, and $H_0^1(\Omega)$, respectively, in the following sense: For all $\mathbf{v}\in\mathbf{H}_0^1(\Omega)$, $\mathsf{q}\in L_0^2(\Omega)$, and $S\in H_0^1(\Omega)$,%
\begin{equation}\label{eq:approx_props}
\begin{gathered}
\lim_{h\to0}\inf_{\mathbf{v}_h\in\mathbf{X}_h}\|\nabla(\mathbf{v}-\mathbf{v}_h)\|_{\mathbf{L}^{2}(\Omega)}=0, \quad
\lim_{h\to0}\inf_{\mathsf{q}_h\in M_h}\|\mathsf{q}-\mathsf{q}_h\|_{L^2(\Omega)}=0,
\\
\lim_{h\to0}\inf_{S_h\in Y_h}\|\nabla(S-S_h)\|_{\mathbf{L}^{2}(\Omega)}=0;
\end{gathered}
\end{equation}
see \cite[Corollaries 1.109 and 1.110]{Guermond-Ern}.
In addition, there are interpolation operators $\boldsymbol{\Pi}_h:\mathbf{H}_0^1(\Omega)\to\mathbf{X}_h$ and $\Pi_h:H_0^1(\Omega)\to Y_h$, for instance of Scott--Zhang or Cl\'ement type \cite{MR1011446,Clement1975}, such that, for every $\mathbf{v}\in\mathbf{H}_0^1(\Omega)$ and $S\in H_0^1(\Omega)$,
\begin{equation}\label{eq:interp_conv}
\lim_{h\to0}\|\nabla(\mathbf{v}-\boldsymbol{\Pi}_h\mathbf{v})\|_{\mathbf{L}^{2}(\Omega)}=0,\qquad\lim_{h\to0}\|\nabla(S-\Pi_hS)\|_{\mathbf{L}^{2}(\Omega)}=0.
\end{equation}

We also need the space of discretely divergence-free functions $\mathbf{V}_h:=\{\mathbf{v}_h\in\mathbf{X}_h:\int_\Omega\mathsf{q}_h\operatorname{div}\mathbf{v}_h\,\mathrm{d}x=0\ \forall\mathsf{q}_h\in M_h\}$. In general, a function in $\mathbf{V}_h$ does not have zero divergence, and thus $\mathbf{V}_h\not\subset\mathbf{V}(\Omega)$. Nevertheless, the discrete inf-sup condition \eqref{eq:infsup_disc} guarantees that $\mathbf{V}(\Omega)$ is well approximated by $\mathbf{V}_h$ \cite[Chapter II, Section 1.1]{MR851383}: For every $\mathbf{v}\in\mathbf{V}(\Omega)$, we have
\begin{equation}\label{eq:approx_Vh}
\inf_{\mathbf{w}_h\in\mathbf{V}_h}\|\nabla(\mathbf{v}-\mathbf{w}_h)\|_{\mathbf{L}^{2}(\Omega)}\leq\left(1+\tilde\beta^{-1}\right)\inf_{\mathbf{v}_h\in\mathbf{X}_h}\|\nabla(\mathbf{v}-\mathbf{v}_h)\|_{\mathbf{L}^{2}(\Omega)}.
\end{equation}

\subsection{The discrete coupled problem}
\label{sec:discrete_problem}
Since $\mathbf{V}_h\not\subset\mathbf{V}(\Omega)$, the property $\int_\Omega(\mathbf{v}\cdot\nabla S)S\,\mathrm{d}x=0$, which holds for all $\mathbf{v}\in\mathbf{V}(\Omega)$ and $S\in H_0^1(\Omega)$, does not hold, in general, for discrete velocities. Therefore, as in \cite[Chapter II, Section 3.2]{MR0603444} and \cite[Chapter IV, Section 3]{MR0548867}, we write the convective terms of the temperature and concentration equations in skew-symmetric form. To be precise, we introduce the term $\frac{1}{2}\int_\Omega((\mathbf{w}_h\cdot\nabla\varphi_h)\psi_h-(\mathbf{w}_h\cdot\nabla\psi_h)\varphi_h)\,\mathrm{d}x$. This term is zero when $\psi_h=\varphi_h$, for every $\mathbf{w}_h\in\mathbf{X}_h$ and $\varphi_h\in Y_h$. Thus, the discrete coupled problem reads as follows: Find $(\mathbf{u}_h,\mathsf{p}_h,T_h,C_h)\in\mathbf{X}_h\times M_h\times Y_h\times Y_h$ such that
\begin{gather}
\langle\mathcal{A}_{T_h,C_h}(\mathbf{u}_h),\mathbf{v}_h\rangle-\int_\Omega\mathsf{p}_h\operatorname{div}\mathbf{v}_h\,\mathrm{d}x=\langle\mathbf{f}(T_h,C_h),\mathbf{v}_h\rangle,\quad
\int_\Omega\mathsf{q}_h\operatorname{div}\mathbf{u}_h\,\mathrm{d}x=0, \nonumber
\\
\int_\Omega\left(\kappa(T_h,C_h)\nabla T_h\cdot\nabla S_h+\tfrac{1}{2}\left[(\mathbf{u}_h\cdot\nabla T_h)S_h-(\mathbf{u}_h\cdot\nabla S_h)T_h\right]\right)\mathrm{d}x=\langle\mathsf{g}_1,S_h\rangle,\label{eq:discrete_coupled}
\\
\int_\Omega\left(\mathcal{D}(T_h,C_h)\nabla C_h\cdot\nabla R_h+\tfrac{1}{2}\left[(\mathbf{u}_h\cdot\nabla C_h)R_h-(\mathbf{u}_h\cdot\nabla R_h)C_h\right]\right)\mathrm{d}x=\langle\mathsf{g}_2,R_h\rangle,\nonumber
\end{gather}
for all $(\mathbf{v}_h,\mathsf{q}_h,S_h,R_h)\in\mathbf{X}_h\times M_h\times Y_h\times Y_h$.
Here, $\mathcal{A}_{T_h,C_h}$ is the discrete version of the map defined in \eqref{eq:operator_ATC}, i.e., $\langle\mathcal{A}_{T_h,C_h}(\mathbf{v}),\mathbf{w}\rangle:=\int_\Omega\nu(T_h,C_h)\,\varepsilon(x,|\nabla\mathbf{v}|)\nabla\mathbf{v}:\nabla\mathbf{w}\,\mathrm{d}x$ for all $\mathbf{v},\mathbf{w}\in\mathbf{H}_0^1(\Omega)$.

\begin{theorem}[existence of discrete solutions]\label{th:existence_discrete}
In the framework of Theorem \ref{th:existence_coupled}, for every $\mathscr{T}_h\in\mathbb{T}$, problem \eqref{eq:discrete_coupled} admits at least one solution $(\mathbf{u}_h,\mathsf{p}_h,T_h,C_h)\in\mathbf{X}_h\times M_h\times Y_h\times Y_h$. Moreover, we have
\begin{align}
\label{eq:estimate_TC_discrete}
\|\nabla T_h\|_{\mathbf{L}^2(\Omega)}
& \leq \kappa_{-}^{-1}\|\mathsf{g}_1\|_{H^{-1}(\Omega)},
\quad
\|\nabla C_h\|_{\mathbf{L}^2(\Omega)}
\leq \mathcal{D}_{-}^{-1}\|\mathsf{g}_2\|_{H^{-1}(\Omega)},
\\
\label{eq:estimate_up_discrete}
\|\nabla\mathbf{u}_h\|_{\mathbf{L}^2(\Omega)}
& \leq (\nu_{-}\varepsilon_{-})^{-1}\Lambda,
\quad
\|\mathsf{p}_h\|_{L^2(\Omega)}
\leq \tilde\beta^{-1}\left(1+\nu_{+}\varepsilon_{+}(\nu_{-}\varepsilon_{-})^{-1}\right)\Lambda.
\end{align}
Here, $\mathcal{C}_{2\hookrightarrow2}$ and $\Lambda$ are as in Theorem~\ref{th:existence_coupled}.
\end{theorem}

\begin{proof}
We proceed as in the proof of Theorem~\ref{th:existence_coupled}, choosing a discrete analog of the map $\Phi$ on $\mathbf{V}_h\times Y_h\times Y_h$ and incorporating the skew-symmetric convective terms. These terms are zero when $S_h=T_h$ and $R_h=C_h$. The existence of a discrete solution follows immediately from \cite[Chapter IV, Corollary 1.1]{MR851383} and the condition \eqref{eq:infsup_disc}. On the other hand, the bounds  \eqref{eq:estimate_TC_discrete} and \eqref{eq:estimate_up_discrete} follow as in Step 5 of the proof of Theorem \ref{th:existence_coupled}, with $\tilde\beta$ instead of $\beta$. This completes the proof.
\end{proof}

\begin{theorem}[uniqueness of discrete solutions]\label{th:uniqueness_discrete_3D}
In the framework of Theorem~\ref{th:existence_discrete}, let $d=3$, and let $(\mathbf{u}_{h,1},\mathsf{p}_{h,1},T_{h,1},C_{h,1})\in\mathbf{X}_h\times M_h\times Y_h\times Y_h$ be a solution of \eqref{eq:discrete_coupled} that satisfies
\begin{align}
\label{eq:small_data_discrete}
\frac{C_{\mathbf{f}}\,\mathcal{C}_{2\hookrightarrow2}^2\mathcal{C}_{4\hookrightarrow2}^2\,\Gamma(\mathsf{g}_1,\mathsf{g}_2)}{\nu_{-}\varepsilon_{-}(1-\theta_h)}
+\frac{\varepsilon_{+}\mathcal{L}_{\nu}\mathcal{C}_{6\hookrightarrow2}\mathcal{C}_{4\hookrightarrow2}^2\,\Gamma(\mathsf{g}_1,\mathsf{g}_2)\,\|\nabla\mathbf{u}_{h,1}\|_{\mathbf{L}^3(\Omega)}}{\nu_{-}\varepsilon_{-}(1-\theta_h)}
&<1,
\\
\label{eq:small_data_discrete_2}
\theta_h:=\frac{\mathcal{L}_{\kappa}\mathcal{C}_{6\hookrightarrow2}\|\nabla T_{h,1}\|_{\mathbf{L}^3(\Omega)}}{\kappa_{-}}
+\frac{\mathcal{L}_{\mathcal{D}}\mathcal{C}_{6\hookrightarrow2}\|\nabla C_{h,1}\|_{\mathbf{L}^3(\Omega)}}{\mathcal{D}_{-}}
&<1,
\end{align}
where $\Gamma(\mathsf{g}_1,\mathsf{g}_2)$, $C_{\mathbf{f}}$, $\mathcal{C}_{4\hookrightarrow2}$, $\mathcal{C}_{2\hookrightarrow2}$, and $\mathcal{C}_{6\hookrightarrow2}$ are as in Theorem~\ref{th:uniqueness_coupled_3D}. Then, problem \eqref{eq:discrete_coupled} has no other solution.
\end{theorem}

\begin{proof}
The proof is a suitable adaptation of the proof of Theorem \ref{th:uniqueness_coupled_3D}.
\end{proof}

Similar to the continuous case, the regularity requirements for the solution in two dimensions can be relaxed.

\begin{theorem}[uniqueness of discrete solutions for $d=2$]\label{th:uniqueness_discrete_2D}
In the framework of Theorem~\ref{th:existence_discrete}, let $d=2$, let $\ell>2$, and let $\sigma$ be such that $1/\sigma+1/\ell=1/2$. Let $(\mathbf{u}_{h,1},\mathsf{p}_{h,1},T_{h,1},C_{h,1})\in\mathbf{X}_h\times M_h\times Y_h\times Y_h$ be a solution of \eqref{eq:discrete_coupled} that satisfies
\begin{align}
\label{eq:small_data_discrete_2D}
\frac{C_{\mathbf{f}}\,\mathcal{C}_{2\hookrightarrow2}^2\mathcal{C}_{4\hookrightarrow2}^2\,\Gamma(\mathsf{g}_1,\mathsf{g}_2)}{\nu_{-}\varepsilon_{-}(1-\theta_{h,\ell})}
+\frac{\varepsilon_{+}\mathcal{L}_{\nu}\mathcal{C}_{\sigma\hookrightarrow2}\mathcal{C}_{4\hookrightarrow2}^2\,\Gamma(\mathsf{g}_1,\mathsf{g}_2)\,\|\nabla\mathbf{u}_{h,1}\|_{\mathbf{L}^{\ell}(\Omega)}}{\nu_{-}\varepsilon_{-}(1-\theta_{h,\ell})}
&<1,
\\
\label{eq:small_data_discrete_2D_2}
\theta_{h,\ell}:=\frac{\mathcal{L}_{\kappa}\mathcal{C}_{\sigma\hookrightarrow2}\|\nabla T_{h,1}\|_{\mathbf{L}^{\ell}(\Omega)}}{\kappa_{-}}
+\frac{\mathcal{L}_{\mathcal{D}}\mathcal{C}_{\sigma\hookrightarrow2}\|\nabla C_{h,1}\|_{\mathbf{L}^{\ell}(\Omega)}}{\mathcal{D}_{-}}
&<1,
\end{align}
where $\Gamma(\mathsf{g}_1,\mathsf{g}_2)$, $C_{\mathbf{f}}$, $\mathcal{C}_{4\hookrightarrow2}$, $\mathcal{C}_{2\hookrightarrow2}$, and $\mathcal{C}_{\sigma\hookrightarrow2}$ are as in Theorem~\ref{th:uniqueness_coupled_2D}. Then, problem \eqref{eq:discrete_coupled} has no other solution.
\end{theorem}

\begin{proof}
The proof is a suitable adaptation of the proof of Theorem~\ref{th:uniqueness_discrete_3D}.
\end{proof}

\subsection{Convergence}

We note that the following result requires no smallness assumptions on the data and no additional regularity of the solutions.

\begin{theorem}
In the framework of Theorem~\ref{th:existence_discrete}, let $(\mathbf{u}_h,\mathsf{p}_h,T_h,C_h)\in\mathbf{X}_h\times M_h\times Y_h\times Y_h$ be a solution of \eqref{eq:discrete_coupled}. Then, there exist a nonrelabeled subsequence of $\{(\mathbf{u}_h,\mathsf{p}_h,T_h,C_h)\}_{h>0}$ and $(\mathbf{u},\mathsf{p},T,C)\in\mathbf{H}_0^1(\Omega)\times L_0^2(\Omega)\times H_0^1(\Omega)\times H_0^1(\Omega)$ such that
\[
\mathbf{u}_h\to\mathbf{u}\ \text{in}\ \mathbf{H}_0^1(\Omega),
\quad
\mathsf{p}_h\to\mathsf{p}\ \text{in}\ L^2(\Omega),
\quad
T_h\to T\ \text{in}\ H_0^1(\Omega),
\quad
C_h\to C\ \text{in}\ H_0^1(\Omega),
\]
as $h\to0$. Moreover, $(\mathbf{u},\mathsf{p},T,C)$ solves the coupled problem \eqref{eq:weak_coupled}.
\end{theorem}

\begin{proof}
To simplify the presentation, we introduce the notation $\nu_h:=\nu(T_h,C_h)$, $\kappa_h:=\kappa(T_h,C_h)$, and $\mathcal{D}_h:=\mathcal{D}(T_h,C_h)$. In view of \eqref{eq:estimate_TC_discrete}--\eqref{eq:estimate_up_discrete}, the discrete solutions are uniformly bounded in $\mathbf{H}_0^1(\Omega)\times L_0^2(\Omega)\times H_0^1(\Omega)\times H_0^1(\Omega)$. Thus, we deduce the existence of a nonrelabeled subsequence $\{(\mathbf{u}_h,\mathsf{p}_h,T_h,C_h)\}_{h>0}$ such that
\begin{equation}\label{eq:conv_weak}
\mathbf{u}_h\rightharpoonup\mathbf{u}\ \text{in}\ \mathbf{H}_0^1(\Omega),
\quad
\mathsf{p}_h\rightharpoonup\mathsf{p}\ \text{in}\ L^2(\Omega),
\quad
T_h\rightharpoonup T,\ C_h\rightharpoonup C\ \text{in}\ H_0^1(\Omega),
\quad h\to0.
\end{equation}
In the following, we show that $(\mathbf{u},\mathsf{p},T,C)$ solves the coupled problem. To accomplish this task, we divide the proof into several steps.

\emph{Step 1.} \emph{Strong convergence $\mathbf{u}_h\to\mathbf{u}$ in $\mathbf{H}_0^1(\Omega)$.} We prove that $\mathbf{u}\in\mathbf{V}(\Omega)$. Let $\mathsf{q}\in L_0^2(\Omega)$, and let $\mathsf{q}_h\in M_h$. Using that $\int_\Omega\mathsf{q}_h\operatorname{div}\mathbf{u}_h\,\mathrm{d}x=0$, \eqref{eq:conv_weak}, \eqref{eq:approx_props}, and the bound $\|\operatorname{div}\mathbf{w}\|_{L^2(\Omega)}\lesssim\|\nabla\mathbf{w}\|_{\mathbf{L}^2(\Omega)}$ for $\mathbf{w}\in\mathbf{H}_0^1(\Omega)$, we obtain that $\int_\Omega\mathsf{q}\operatorname{div}\mathbf{u}\,\mathrm{d}x=0$ for all $\mathsf{q}\in L_0^2(\Omega)$. Thus, $\mathbf{u}\in\mathbf{V}(\Omega)$.

Since $\mathbf{u}\in\mathbf{V}(\Omega)$, from \eqref{eq:approx_Vh} and \eqref{eq:approx_props} it follows that there exists $\mathbf{w}_h\in\mathbf{V}_h$ such that $\|\nabla(\mathbf{u}-\mathbf{w}_h)\|_{\mathbf{L}^2(\Omega)}\to0$ as $h\to0$. We set $\mathbf{v}_h=\mathbf{e}_h:=\mathbf{u}_h-\mathbf{w}_h\in\mathbf{V}_h$ in the first equation of \eqref{eq:discrete_coupled}, and note that $\mathbf{e}_h\rightharpoonup\mathbf{0}$ in $\mathbf{H}_0^1(\Omega)$. Since $\int_\Omega\mathsf{p}_h\operatorname{div}\mathbf{e}_h\,\mathrm{d}x=0$, the property \eqref{eq:prop_02} for $\mathcal{A}_{T_h,C_h}$ yields
\begin{align*}
\nu_{-}\varepsilon_{-}\|\nabla\mathbf{e}_h\|_{\mathbf{L}^2(\Omega)}^2
\leq
\langle\mathbf{f}(T_h,C_h),\mathbf{e}_h\rangle
-\langle\mathcal{A}_{T_h,C_h}(\mathbf{w}_h),\mathbf{e}_h\rangle
=:\mathrm{I}_h+\mathrm{II}_h.
\end{align*}

We note that the term $\mathrm{I}_h\to0$, due to the Sobolev embedding, the bounds $|\mathbf{f}_i(s)|\leq C_{\mathbf{f}_i}|s|$ and \eqref{eq:estimate_TC_discrete}, and $\mathbf{e}_h\to\mathbf{0}$ in $\mathbf{L}^2(\Omega)$. To control the term $\mathrm{II}_h$, we note the identity $
\mathrm{II}_h=\langle\mathcal{A}_{T_h,C_h}(\mathbf{u})-\mathcal{A}_{T_h,C_h}(\mathbf{w}_h),\mathbf{e}_h\rangle-\langle\mathcal{A}_{T_h,C_h}(\mathbf{u})-\mathcal{A}_{T,C}(\mathbf{u}),\mathbf{e}_h\rangle-\langle\mathcal{A}_{T,C}(\mathbf{u}),\mathbf{e}_h\rangle$. Since $\nu$ is continuous and bounded, and according to \eqref{assumption_1} with $\boldsymbol{\omega}=\mathbf{0}$, and by the Lebesgue dominated convergence theorem, we can obtain
\begin{equation}\label{eq:conv_A_u}
\|\mathcal{A}_{T_h,C_h}(\mathbf{u})-\mathcal{A}_{T,C}(\mathbf{u})\|_{\mathbf{H}^{-1}(\Omega)}
\leq
\|(\nu_h-\nu(T,C))\,\varepsilon(x,|\nabla\mathbf{u}|)\nabla\mathbf{u}\|_{\mathbf{L}^2(\Omega)}
\to0,
\end{equation}
as $h\to 0$. Since $\{\mathbf{e}_h\}_{h>0}$ is bounded in $\mathbf{H}_0^1(\Omega)$, by \eqref{eq:estimate_up_discrete} and $\|\nabla(\mathbf{u}-\mathbf{w}_h)\|_{\mathbf{L}^2(\Omega)}\to0$, the first two terms in the identity of $\mathrm{II}_h$ converge to zero, due to \eqref{eq:prop_01} and \eqref{eq:conv_A_u}, respectively. The last term converges to zero because $\mathbf{e}_h\rightharpoonup\mathbf{0}$ in $\mathbf{H}_0^1(\Omega)$. Therefore, the term $\mathrm{II}_h\to 0$ as $h\to 0$. Thus, we can conclude that $\|\nabla\mathbf{e}_h\|_{\mathbf{L}^2(\Omega)}\to0$ and by the triangle inequality $\|\nabla(\mathbf{u}_h-\mathbf{u})\|_{\mathbf{L}^2(\Omega)}\to 0$. The strong convergence of $\mathbf{u}_h$, \eqref{eq:prop_01}, and \eqref{eq:conv_A_u} also produce the following convergence as $h\to 0$
\begin{multline}\label{eq:conv_A}
\|\mathcal{A}_{T_h,C_h}(\mathbf{u}_h)-\mathcal{A}_{T,C}(\mathbf{u})\|_{\mathbf{H}^{-1}(\Omega)}
\\
\lesssim
\|\nabla(\mathbf{u}_h-\mathbf{u})\|_{\mathbf{L}^2(\Omega)}
+\|\mathcal{A}_{T_h,C_h}(\mathbf{u})-\mathcal{A}_{T,C}(\mathbf{u})\|_{\mathbf{H}^{-1}(\Omega)}
\to0.%
\end{multline}

\emph{Step 2.} \emph{The first equation in \eqref{eq:weak_coupled}.} Let $\mathbf{v}\in\mathbf{H}_0^1(\Omega)$, and set $\mathbf{v}_h:=\boldsymbol{\Pi}_h\mathbf{v}\in\mathbf{X}_h$, with $\boldsymbol{\Pi}_h$ as in \eqref{eq:interp_conv}. Since $(\mathbf{u}_h,\mathsf{p}_h,T_h,C_h)$ satisfies the first equation of \eqref{eq:discrete_coupled} with $\mathbf{v}_h$, we can write $\langle\mathcal{A}_{T,C}(\mathbf{u}),\mathbf{v}\rangle
-\int_\Omega\mathsf{p}\operatorname{div}\mathbf{v}\,\mathrm{d}x
-\langle\mathbf{f}(T,C),\mathbf{v}\rangle
=\mathrm{I}_h+\mathrm{II}_h+\mathrm{III}_h$, with
\begin{equation*}
\begin{gathered}
\mathrm{I}_h:=\langle\mathcal{A}_{T,C}(\mathbf{u}),\mathbf{v}\rangle-\langle\mathcal{A}_{T_h,C_h}(\mathbf{u}_h),\mathbf{v}_h\rangle,
\qquad
\mathrm{II}_h:=-\int_\Omega\left(\mathsf{p}\operatorname{div}\mathbf{v}-\mathsf{p}_h\operatorname{div}\mathbf{v}_h\right)\mathrm{d}x,
\\
\mathrm{III}_h:=-\langle\mathbf{f}(T,C),\mathbf{v}\rangle+\langle\mathbf{f}(T_h,C_h),\mathbf{v}_h\rangle.
\end{gathered}
\end{equation*}
We now show that each term converges to zero. In view of \eqref{eq:conv_A}, \eqref{eq:interp_conv}, \eqref{eq:estimate_up_discrete}, and \eqref{eq:prop_01}, we have
\[
|\mathrm{I}_h|
\lesssim
\|\mathcal{A}_{T,C}(\mathbf{u})-\mathcal{A}_{T_h,C_h}(\mathbf{u}_h)\|_{\mathbf{H}^{-1}(\Omega)}\|\nabla\mathbf{v}\|_{\mathbf{L}^2(\Omega)}
+\|\nabla\mathbf{u}_h\|_{\mathbf{L}^2(\Omega)}\|\nabla(\mathbf{v}-\mathbf{v}_h)\|_{\mathbf{L}^2(\Omega)}
\to0.
\]
For the term $\mathrm{II}_h$, we use \eqref{eq:conv_weak}, \eqref{eq:estimate_up_discrete}, \eqref{eq:interp_conv}, and the estimate $\|\operatorname{div}\mathbf{w}\|_{L^2(\Omega)}\lesssim\|\nabla\mathbf{w}\|_{\mathbf{L}^2(\Omega)}$ to conclude that the term $\mathrm{II}_h$ satisfies $|\mathrm{II}_h|
\leq
\left|\int_\Omega(\mathsf{p}-\mathsf{p}_h)\operatorname{div}\mathbf{v}\,\mathrm{d}x\right|
+\|\mathsf{p}_h\|_{L^2(\Omega)}\|\operatorname{div}(\mathbf{v}-\mathbf{v}_h)\|_{L^2(\Omega)}
\to0$ as $h\to 0$. Finally, the Lipschitz condition of $\mathbf{f}_1$ and $\mathbf{f}_2$, the bound $|\mathbf{f}_i(s)|\leq C_{\mathbf{f}_i}|s|$, the convergences $T_h\to T$ and $C_h\to C$ in $L^2(\Omega)$, \eqref{eq:interp_conv}, and the embedding $H_0^1(\Omega)\hookrightarrow L^2(\Omega)$ show that
\begin{multline*}
|\mathrm{III}_h|
\leq
\|\mathbf{f}_0\|_{\mathbf{H}^{-1}(\Omega)}\|\nabla(\mathbf{v}-\mathbf{v}_h)\|_{\mathbf{L}^2(\Omega)}
+\left(C_{\mathbf{f}_1}\|T_h\|_{L^2(\Omega)}+C_{\mathbf{f}_2}\|C_h\|_{L^2(\Omega)}\right)\|\mathbf{v}-\mathbf{v}_h\|_{\mathbf{L}^2(\Omega)}
\\
+\left(C_{\mathbf{f}_1}\|T-T_h\|_{L^2(\Omega)}+C_{\mathbf{f}_2}\|C-C_h\|_{L^2(\Omega)}\right)\|\mathbf{v}\|_{\mathbf{L}^2(\Omega)}
\to0.
\end{multline*}
Thus, $(\mathbf{u},\mathsf{p},T,C)$ satisfies the first equation in \eqref{eq:weak_coupled}.

\emph{Step 3.} \emph{The temperature and concentration equations in \eqref{eq:weak_coupled}.} Let $S\in H_0^1(\Omega)$ be arbitrary, and set $S_h:=\Pi_hS\in Y_h$, with $\Pi_h$ as in \eqref{eq:interp_conv}. As in Step 2, the temperature equation in \eqref{eq:discrete_coupled} allows us to write
\begin{equation*}
\begin{gathered}
\int_\Omega\kappa(T,C)\nabla T\cdot\nabla S\,\mathrm{d}x
+\int_\Omega(\mathbf{u}\cdot\nabla T)S\,\mathrm{d}x
-\langle\mathsf{g}_1,S\rangle
=\mathrm{I}_h+\mathrm{II}_h-\langle\mathsf{g}_1,S-S_h\rangle,
\\
\mathrm{I}_h:=\int_\Omega\left(\kappa(T,C)\nabla T\cdot\nabla S-\kappa_h\nabla T_h\cdot\nabla S_h\right)\mathrm{d}x,
\\
\mathrm{II}_h:=\int_\Omega(\mathbf{u}\cdot\nabla T)S\,\mathrm{d}x-\frac{1}{2}\int_\Omega\left((\mathbf{u}_h\cdot\nabla T_h)S_h-(\mathbf{u}_h\cdot\nabla S_h)T_h\right)\mathrm{d}x.
\end{gathered}
\end{equation*}
We observe that the term $\langle\mathsf{g}_1,S-S_h\rangle$ converges to zero in view of \eqref{eq:interp_conv}. For the term $\mathrm{I}_h$, we have that $|\mathrm{I}_h|
\lesssim
\left|\int_\Omega\left(\kappa(T,C)\nabla S\cdot\nabla T-\kappa_h\nabla S\cdot\nabla T_h\right)\mathrm{d}x\right|
+\|\nabla T_h\|_{\mathbf{L}^2(\Omega)}\|\nabla(S-S_h)\|_{\mathbf{L}^2(\Omega)}$.
Indeed, the first term converges to zero because $\nabla T_h\rightharpoonup\nabla T$ and $\kappa_h\nabla S\to\kappa(T,C)\nabla S$ in $\mathbf{L}^2(\Omega)$, and the second term converges to zero by \eqref{eq:estimate_TC_discrete} and \eqref{eq:interp_conv}. For the term $\mathrm{II}_h$, since $\mathbf{u}\in\mathbf{V}(\Omega)$, an integration by parts gives $\int_\Omega(\mathbf{u}\cdot\nabla T)S\,\mathrm{d}x=\frac{1}{2}\int_\Omega((\mathbf{u}\cdot\nabla T)S-(\mathbf{u}\cdot\nabla S)T)\,\mathrm{d}x$, and thus $2\,\mathrm{II}_h$ is bounded by the terms
\[
|2\,\mathrm{II}_h|\lesssim \left|\int_\Omega\left((\mathbf{u}\cdot\nabla T)S-(\mathbf{u}_h\cdot\nabla T_h)S_h\right)\mathrm{d}x\right|+\left|\int_\Omega\left((\mathbf{u}\cdot\nabla S)T-(\mathbf{u}_h\cdot\nabla S_h)T_h\right)\mathrm{d}x\right|.
\]
To control the first and second terms in the above estimate, we use similar arguments. For the sake of brevity, we present only the first term. To do this, we invoke the convergences $\mathbf{u}_h\to\mathbf{u}$ in $\mathbf{L}^4(\Omega)$ and $\nabla T_h\rightharpoonup\nabla T$ in $\mathbf{L}^2(\Omega)$, the uniform bounds of $\{\mathbf{u}_h\}_{h>0}$ and $\{T_h\}_{h>0}$ in $\mathbf{H}_0^1(\Omega)$ and $H_0^1(\Omega)$, respectively, and \eqref{eq:interp_conv} to obtain that
\begin{multline*}
\left|\int_\Omega\left((\mathbf{u}\cdot\nabla T)S-(\mathbf{u}_h\cdot\nabla T_h)S_h\right)\mathrm{d}x\right|
\lesssim \left|\int_\Omega\left((\mathbf{u}-\mathbf{u}_h)\cdot\nabla T\right)S \mathrm{d}x\right|
\\
+
\left|\int_{\Omega}\mathbf{u}_h\cdot\nabla(T-T_h)S\mathrm{d}x\right|+
\left|\int_{\Omega}(\mathbf{u}_h\cdot\nabla T_h)(S-S_h)\mathrm{d}x\right|\to 0,\quad\text{as}\quad h\to 0.
\end{multline*}
Thus, the term $\mathrm{II}_h\to 0$. Therefore, $T$ solves the temperature equation in \eqref{eq:weak_coupled}. We can use similar arguments to show that $C$ satisfies the concentration equation in \eqref{eq:weak_coupled}.

\emph{Step 4.} \emph{Strong convergence $T_h\to T$ and $C_h\to C$ in $H_0^1(\Omega)$.} Since the skew-symmetric term is zero when $S_h=T_h$, the temperature equation in \eqref{eq:discrete_coupled} with $S_h=T_h$ gives $\int_\Omega\kappa_h|\nabla T_h|^2\,\mathrm{d}x=\langle\mathsf{g}_1,T_h\rangle$. Similarly, since $\mathbf{u}\in\mathbf{V}(\Omega)$ we have $\int_\Omega(\mathbf{u}\cdot\nabla T)T\,\mathrm{d}x=0$, and thus the temperature equation in \eqref{eq:weak_coupled} with $S=T$ gives $\int_\Omega\kappa(T,C)|\nabla T|^2\,\mathrm{d}x=\langle\mathsf{g}_1,T\rangle$. Thus, the bound $\kappa_h\geq\kappa_{-}$ yields
\[
\kappa_{-}\|\nabla(T_h-T)\|_{\mathbf{L}^2(\Omega)}^2
\leq
\langle\mathsf{g}_1,T_h\rangle
-2\int_\Omega\kappa_h\nabla T\cdot\nabla T_h\,\mathrm{d}x
+\int_\Omega\kappa_h|\nabla T|^2\,\mathrm{d}x.
\]
Since $T_h\rightharpoonup T$ in $H_0^1(\Omega)$ and $\kappa_h\nabla T\to\kappa(T,C)\nabla T$ in $\mathbf{L}^2(\Omega)$, as in the proof of \eqref{eq:conv_A_u}, the right-hand side converges to $\langle\mathsf{g}_1,T\rangle-2\int_\Omega\kappa(T,C)|\nabla T|^2\,\mathrm{d}x+\int_\Omega\kappa(T,C)|\nabla T|^2\,\mathrm{d}x=0$. Therefore, $\|\nabla(T_h-T)\|_{\mathbf{L}^2(\Omega)}\to0$ as $h\to0$. The same arguments, applied to the concentration equations, show that the term $\|\nabla(C_h-C)\|_{\mathbf{L}^2(\Omega)}\to0$ as $h\to0$.

\emph{Step 5.} \emph{Strong convergence $\mathsf{p}_h\to\mathsf{p}$ in $L^2(\Omega)$.} Let $\mathsf{q}_h\in M_h$ be the $L^2$-projection of $\mathsf{p}$ onto $M_h$, so that $\|\mathsf{p}-\mathsf{q}_h\|_{L^2(\Omega)}\to0$ as $h\to0$, according to \eqref{eq:approx_props}. If we subtract the first equations of \eqref{eq:discrete_coupled} and \eqref{eq:weak_coupled} with $\mathbf{v}=\mathbf{v}_h$, for all $\mathbf{v}_h\in\mathbf{X}_h$, we obtain
$
\int_\Omega(\mathsf{q}_h-\mathsf{p}_h)\operatorname{div}\mathbf{v}_h\,\mathrm{d}x
=
\int_\Omega(\mathsf{q}_h-\mathsf{p})\operatorname{div}\mathbf{v}_h\,\mathrm{d}x
+\langle\mathcal{A}_{T,C}(\mathbf{u})-\mathcal{A}_{T_h,C_h}(\mathbf{u}_h),\mathbf{v}_h\rangle
-\langle\mathbf{f}(T,C)-\mathbf{f}(T_h,C_h),\mathbf{v}_h\rangle$.
This identity, the discrete inf-sup condition \eqref{eq:infsup_disc}, and the bound $\|\operatorname{div}\mathbf{v}_h\|_{L^2(\Omega)}\lesssim\|\nabla\mathbf{v}_h\|_{\mathbf{L}^2(\Omega)}$ produce
\begin{multline*}
\|\mathsf{q}_h-\mathsf{p}_h\|_{L^2(\Omega)}
\lesssim
\|\mathsf{q}_h-\mathsf{p}\|_{L^2(\Omega)}
\\
+\|\mathcal{A}_{T,C}(\mathbf{u})-\mathcal{A}_{T_h,C_h}(\mathbf{u}_h)\|_{\mathbf{H}^{-1}(\Omega)}
+\|\mathbf{f}(T,C)-\mathbf{f}(T_h,C_h)\|_{\mathbf{H}^{-1}(\Omega)}.
\end{multline*}
The second term on the right-hand side converges to zero in view of \eqref{eq:conv_A}. For the last term, we use the Lipschitz continuity of $\mathbf{f}_1$ and $\mathbf{f}_2$, the Sobolev embedding $H_0^1(\Omega)\hookrightarrow L^2(\Omega)$, and the convergences $T_h\to T$ and $C_h\to C$ in $L^2(\Omega)$ to conclude that $\|\mathbf{f}(T,C)-\mathbf{f}(T_h,C_h)\|_{\mathbf{H}^{-1}(\Omega)}\leq\mathcal{C}_{2\hookrightarrow2}(C_{\mathbf{f}_1}\|T-T_h\|_{L^2(\Omega)}+C_{\mathbf{f}_2}\|C-C_h\|_{L^2(\Omega)})\to0$ as $h\to 0$. Thus, we obtain that $\|\mathsf{q}_h-\mathsf{p}_h\|_{L^2(\Omega)}\to0$ and, by the triangle inequality, $\|\mathsf{p}_h-\mathsf{p}\|_{L^2(\Omega)}\to0$ as $h\to0$. This completes the proof.
\end{proof}

\subsection{An approximation result}
\label{sec:quasi_best}
 In this section, for the finite element approximation \eqref{eq:discrete_coupled} of the coupled problem \eqref{eq:weak_coupled}, we derive a quasi-best approximation result. For this purpose, we assume suitable smallness conditions and extra regularity assumptions. To begin the analysis, we introduce the error terms $\mathbf{e}_{\mathbf{u}}:=\mathbf{u}-\mathbf{u}_h$, $e_{\mathsf{p}}:=\mathsf{p}-\mathsf{p}_h$, $e_{T}:=T-T_h$, and $e_{C}:=C-C_h$.

\begin{theorem}[quasi-best approximation]\label{th:quasi_best_3D}
In the framework of Theorem~\ref{th:existence_discrete}, let $d=3$, and let $(\mathbf{u},\mathsf{p},T,C)\in[\mathbf{W}^{1,3}(\Omega)\cap\mathbf{H}_0^1(\Omega)]\times L_0^2(\Omega)\times[W^{1,3}(\Omega)\cap H_0^1(\Omega)]\times[W^{1,3}(\Omega)\cap H_0^1(\Omega)]$ be a solution of the coupled problem \eqref{eq:weak_coupled} that satisfies the smallness conditions \eqref{eq:small_data_coupled}--\eqref{eq:small_data_coupled_2}, with $(\mathbf{u}_1,T_1,C_1)$ replaced by $(\mathbf{u},T,C)$. For every $\mathscr{T}_h\in\mathbb{T}$, let $(\mathbf{u}_h,\mathsf{p}_h,T_h,C_h)\in\mathbf{X}_h\times M_h\times Y_h\times Y_h$ be a solution of \eqref{eq:discrete_coupled}. Then,%
\begin{multline}\label{eq:quasi_best}
\|\nabla\mathbf{e}_{\mathbf{u}}\|_{\mathbf{L}^2(\Omega)}
+\|e_{\mathsf{p}}\|_{L^2(\Omega)}
+\|\nabla e_T\|_{\mathbf{L}^2(\Omega)}
+\|\nabla e_C\|_{\mathbf{L}^2(\Omega)}
\lesssim
\inf_{\mathsf{q}_h\in M_h}\|\mathsf{p}-\mathsf{q}_h\|_{L^2(\Omega)}
\\
+\inf_{\mathbf{w}_h\in\mathbf{V}_h}\|\nabla(\mathbf{u}-\mathbf{w}_h)\|_{\mathbf{L}^2(\Omega)}
+\inf_{S_h\in Y_h}\|\nabla(T-S_h)\|_{\mathbf{L}^2(\Omega)}
+\inf_{R_h\in Y_h}\|\nabla(C-R_h)\|_{\mathbf{L}^2(\Omega)},
\end{multline}
where the hidden constant does not depend on $h$.
\end{theorem}

\begin{proof}
To simplify the presentation of the proof, we introduce the notation $\kappa_h:=\kappa(T_h,C_h)$, $\mathfrak{M}:=1-\theta$, with $\theta$ as in \eqref{eq:small_data_coupled_2}, $\mathfrak{R}:=C_{\mathbf{f}}\,\mathcal{C}_{2\hookrightarrow2}^2+\varepsilon_{+}\mathcal{L}_{\nu}\mathcal{C}_{6\hookrightarrow2}\|\nabla\mathbf{u}\|_{\mathbf{L}^3(\Omega)}$, and $\mathfrak{N}:=\nu_{-}\varepsilon_{-}-\mathfrak{M}^{-1}\mathcal{C}_{4\hookrightarrow2}^2\,\Gamma(\mathsf{g}_1,\mathsf{g}_2)\,\mathfrak{R}$. We note that the conditions \eqref{eq:small_data_coupled_2} and \eqref{eq:small_data_coupled} imply that $\mathfrak{M}>0$ and $\mathfrak{N}>0$, respectively. Thus, we divide the proof into four steps.

\emph{Step 1.} \emph{Temperature and concentration.} Let $S_h\in Y_h$, and define $\eta_T:=T-S_h$ and $\zeta_T:=S_h-T_h$. Since $\mathbf{u}\in\mathbf{V}(\Omega)$, an integration by parts gives $\int_\Omega(\mathbf{u}\cdot\nabla T)S\,\mathrm{d}x=\frac{1}{2}\int_\Omega((\mathbf{u}\cdot\nabla T)S-(\mathbf{u}\cdot\nabla S)T)\,\mathrm{d}x$ for all $S\in H_0^1(\Omega)$. We now subtract the temperature equations of \eqref{eq:discrete_coupled} and \eqref{eq:weak_coupled} with the test function $\zeta_T$, and use the identity $\kappa(T,C)\nabla T-\kappa_h\nabla T_h=(\kappa(T,C)-\kappa_h)\nabla T+\kappa_h\nabla(\eta_T+\zeta_T)$. On the other hand, we write $\mathbf{u}=\mathbf{e}_{\mathbf{u}}+\mathbf{u}_h$ and $e_T=\eta_T+\zeta_T$, and note that the skew-symmetric term with $\mathbf{w}_h=\mathbf{u}_h$ and $\varphi_h=\psi_h=\zeta_T$ vanishes. Thus, we obtain that
\begin{multline*}
\int_\Omega\kappa_h|\nabla\zeta_T|^2\,\mathrm{d}x
=
-\int_\Omega\left((\kappa(T,C)-\kappa_h)\nabla T+\kappa_h\nabla\eta_T\right)\cdot\nabla\zeta_T\,\mathrm{d}x
\\
-\frac{1}{2}\int_\Omega\left((\mathbf{e}_{\mathbf{u}}\cdot\nabla T)\zeta_T-(\mathbf{e}_{\mathbf{u}}\cdot\nabla\zeta_T)T+(\mathbf{u}_h\cdot\nabla\eta_T)\zeta_T-(\mathbf{u}_h\cdot\nabla\zeta_T)\eta_T\right)\mathrm{d}x.
\end{multline*}

We now use the bounds of $\kappa_h$, the Lipschitz continuity of $\kappa$, H\"older's inequality, the Sobolev embeddings $H_0^1(\Omega)\hookrightarrow L^6(\Omega)$ and $H_0^1(\Omega)\hookrightarrow L^4(\Omega)$, and the bounds $\|\nabla T\|_{\mathbf{L}^2(\Omega)}\leq\kappa_{-}^{-1}\|\mathsf{g}_1\|_{H^{-1}(\Omega)}$ and $\|\nabla\mathbf{u}_h\|_{\mathbf{L}^2(\Omega)}\leq(\nu_{-}\varepsilon_{-})^{-1}\Lambda$.
Thus, if we divide by $\|\nabla\zeta_T\|_{\mathbf{L}^2(\Omega)}$, we obtain that
\begin{multline*}
\kappa_{-}\|\nabla\zeta_T\|_{\mathbf{L}^2(\Omega)}
\leq
\mathcal{L}_{\kappa}\mathcal{C}_{6\hookrightarrow2}\|\nabla T\|_{\mathbf{L}^3(\Omega)}\left(\|\nabla e_T\|_{\mathbf{L}^2(\Omega)}+\|\nabla e_C\|_{\mathbf{L}^2(\Omega)}\right)
\\
+\left(\kappa_{+}+\mathcal{C}_{4\hookrightarrow2}^2(\nu_{-}\varepsilon_{-})^{-1}\Lambda\right)\|\nabla\eta_T\|_{\mathbf{L}^2(\Omega)}
+\kappa_{-}^{-1}\mathcal{C}_{4\hookrightarrow2}^2\|\mathsf{g}_1\|_{H^{-1}(\Omega)}\|\nabla\mathbf{e}_{\mathbf{u}}\|_{\mathbf{L}^2(\Omega)}.
\end{multline*}
If we use the triangle inequality $\|\nabla e_T\|_{\mathbf{L}^2(\Omega)}\leq\|\nabla\eta_T\|_{\mathbf{L}^2(\Omega)}+\|\nabla\zeta_T\|_{\mathbf{L}^2(\Omega)}$, we obtain
\begin{multline*}
\|\nabla e_T\|_{\mathbf{L}^2(\Omega)}
\leq
\frac{\mathcal{L}_{\kappa}\mathcal{C}_{6\hookrightarrow2}\|\nabla T\|_{\mathbf{L}^3(\Omega)}}{\kappa_{-}}\left(\|\nabla e_T\|_{\mathbf{L}^2(\Omega)}+\|\nabla e_C\|_{\mathbf{L}^2(\Omega)}\right)
\\
+\kappa_{-}^{-1}\left(\kappa_{-}+\kappa_{+}+\mathcal{C}_{4\hookrightarrow2}^2(\nu_{-}\varepsilon_{-})^{-1}\Lambda\right)\|\nabla\eta_T\|_{\mathbf{L}^2(\Omega)}
+\kappa_{-}^{-2}\mathcal{C}_{4\hookrightarrow2}^2\|\mathsf{g}_1\|_{H^{-1}(\Omega)}\|\nabla\mathbf{e}_{\mathbf{u}}\|_{\mathbf{L}^2(\Omega)}.
\end{multline*}
On the other hand, similar arguments for the concentration equation, with $R_h\in Y_h$, $\eta_C:=C-R_h$, and $\zeta_C:=R_h-C_h$, yield a bound for $\|\nabla e_C\|_{\mathbf{L}^2(\Omega)}$ that follows from the previous one if we replace $\nabla T$, $\eta_T$, $\kappa_{\pm}$, $\mathcal{L}_{\kappa}$, and $\mathsf{g}_1$ by $\nabla C$, $\eta_C$, $\mathcal{D}_{\pm}$, $\mathcal{L}_{\mathcal{D}}$, and $\mathsf{g}_2$, respectively. On the other hand, if we add these two estimates, use the definitions of $\theta$ and $\Gamma(\mathsf{g}_1,\mathsf{g}_2)$, we can conclude that
\begin{multline}\label{eq:error_TC}
\mathfrak{M}\left(\|\nabla e_T\|_{\mathbf{L}^2(\Omega)}+\|\nabla e_C\|_{\mathbf{L}^2(\Omega)}\right)
\\
\leq
\mathfrak{K}\left(\|\nabla\eta_T\|_{\mathbf{L}^2(\Omega)}+\|\nabla\eta_C\|_{\mathbf{L}^2(\Omega)}\right)
+\mathcal{C}_{4\hookrightarrow2}^2\,\Gamma(\mathsf{g}_1,\mathsf{g}_2)\,\|\nabla\mathbf{e}_{\mathbf{u}}\|_{\mathbf{L}^2(\Omega)},
\end{multline}
where $\mathfrak{K}:=\max\{\kappa_{-}^{-1}(\kappa_{-}+\kappa_{+}+\mathcal{C}_{4\hookrightarrow2}^2(\nu_{-}\varepsilon_{-})^{-1}\Lambda),\mathcal{D}_{-}^{-1}(\mathcal{D}_{-}+\mathcal{D}_{+}+\mathcal{C}_{4\hookrightarrow2}^2(\nu_{-}\varepsilon_{-})^{-1}\Lambda)\}$.

\emph{Step 2.} \emph{Pressure.} Let $\mathsf{q}_h\in M_h$. If we subtract the first equations of \eqref{eq:discrete_coupled} and \eqref{eq:weak_coupled} with $\mathbf{v}=\mathbf{v}_h$, we obtain the following identity
\begin{equation}\label{eq:error_momentum}
\langle\mathcal{A}_{T,C}(\mathbf{u})-\mathcal{A}_{T_h,C_h}(\mathbf{u}_h),\mathbf{v}_h\rangle
-\int_\Omega e_{\mathsf{p}}\operatorname{div}\mathbf{v}_h\,\mathrm{d}x
=\langle\mathbf{f}(T,C)-\mathbf{f}(T_h,C_h),\mathbf{v}_h\rangle,
\end{equation}
for all $\mathbf{v}_h\in\mathbf{X}_h$. To continue, we invoke the Lipschitz continuity of $\nu$, $\mathbf{f}_1$, and $\mathbf{f}_2$, the bound $|\varepsilon(x,|\boldsymbol{\tau}|)\boldsymbol{\tau}|\leq\varepsilon_{+}|\boldsymbol{\tau}|$ (see \eqref{assumption_1} with $\boldsymbol{\omega}=\mathbf{0}$), H\"older's inequality, and the Sobolev embeddings $H_0^1(\Omega)\hookrightarrow L^2(\Omega)$ and $H_0^1(\Omega)\hookrightarrow L^6(\Omega)$ to obtain
\begin{multline}\label{eq:diff_Af}
\|\mathcal{A}_{T,C}(\mathbf{u})-\mathcal{A}_{T_h,C_h}(\mathbf{u})\|_{\mathbf{H}^{-1}(\Omega)}
+\|\mathbf{f}(T,C)-\mathbf{f}(T_h,C_h)\|_{\mathbf{H}^{-1}(\Omega)}
\\
\leq
\mathfrak{R}\left(\|\nabla e_T\|_{\mathbf{L}^2(\Omega)}+\|\nabla e_C\|_{\mathbf{L}^2(\Omega)}\right).
\end{multline}
In \eqref{eq:error_momentum}, we write $e_{\mathsf{p}}=(\mathsf{p}-\mathsf{q}_h)+(\mathsf{q}_h-\mathsf{p}_h)$, and add and subtract the term $\mathcal{A}_{T_h,C_h}(\mathbf{u})$. Thus, for every $\mathbf{v}_h\in\mathbf{X}_h$, we have
\begin{multline*}
\int_\Omega(\mathsf{q}_h-\mathsf{p}_h)\operatorname{div}\mathbf{v}_h\,\mathrm{d}x
=
\langle\mathcal{A}_{T,C}(\mathbf{u})-\mathcal{A}_{T_h,C_h}(\mathbf{u}),\mathbf{v}_h\rangle
+\langle\mathcal{A}_{T_h,C_h}(\mathbf{u})-\mathcal{A}_{T_h,C_h}(\mathbf{u}_h),\mathbf{v}_h\rangle
\\
-\langle\mathbf{f}(T,C)-\mathbf{f}(T_h,C_h),\mathbf{v}_h\rangle
-\int_\Omega(\mathsf{p}-\mathsf{q}_h)\operatorname{div}\mathbf{v}_h\,\mathrm{d}x.
\end{multline*}
On the other hand, we use \eqref{eq:diff_Af} in the first and third terms, \eqref{eq:prop_01} for $\mathcal{A}_{T_h,C_h}$ in the second term, and the bound $\|\operatorname{div}\mathbf{v}_h\|_{L^2(\Omega)}\leq\sqrt{d}\,\|\nabla\mathbf{v}_h\|_{\mathbf{L}^2(\Omega)}$ in the last term. Thus, the discrete inf-sup condition \eqref{eq:infsup_disc} and the triangle inequality $\|e_{\mathsf{p}}\|_{L^2(\Omega)}\leq\|\mathsf{p}-\mathsf{q}_h\|_{L^2(\Omega)}+\|\mathsf{q}_h-\mathsf{p}_h\|_{L^2(\Omega)}$ allow us to conclude the estimate
\begin{multline}\label{eq:error_p}
\tilde\beta\|e_{\mathsf{p}}\|_{L^2(\Omega)}
\leq
(\tilde\beta+\sqrt{d})\|\mathsf{p}-\mathsf{q}_h\|_{L^2(\Omega)}
+\nu_{+}\varepsilon_{+}\|\nabla\mathbf{e}_{\mathbf{u}}\|_{\mathbf{L}^2(\Omega)}
\\
+\mathfrak{R}\left(\|\nabla e_T\|_{\mathbf{L}^2(\Omega)}+\|\nabla e_C\|_{\mathbf{L}^2(\Omega)}\right).
\end{multline}

\emph{Step 3.} \emph{Velocity.} Let $\mathbf{w}_h\in\mathbf{V}_h$, and define $\eta_{\mathbf{u}}:=\mathbf{u}-\mathbf{w}_h$ and $\zeta_{\mathbf{u}}:=\mathbf{w}_h-\mathbf{u}_h$. In view of the second equation in \eqref{eq:discrete_coupled}, we have $\mathbf{u}_h\in\mathbf{V}_h$, and thus $\zeta_{\mathbf{u}}\in\mathbf{V}_h$. We now set $\mathbf{v}_h=\zeta_{\mathbf{u}}$ in \eqref{eq:error_momentum}, and note that $\int_\Omega e_{\mathsf{p}}\operatorname{div}\zeta_{\mathbf{u}}\,\mathrm{d}x=\int_\Omega(\mathsf{p}-\mathsf{q}_h)\operatorname{div}\zeta_{\mathbf{u}}\,\mathrm{d}x$, due to $\mathsf{q}_h-\mathsf{p}_h\in M_h$ and $\zeta_{\mathbf{u}}\in\mathbf{V}_h$. On the other hand, if we add and subtract the terms $\mathcal{A}_{T_h,C_h}(\mathbf{u})$ and $\mathcal{A}_{T_h,C_h}(\mathbf{w}_h)$, we can produce that
\begin{multline*}
\langle\mathcal{A}_{T_h,C_h}(\mathbf{w}_h)-\mathcal{A}_{T_h,C_h}(\mathbf{u}_h),\zeta_{\mathbf{u}}\rangle
=
\langle\mathcal{A}_{T_h,C_h}(\mathbf{w}_h)-\mathcal{A}_{T_h,C_h}(\mathbf{u}),\zeta_{\mathbf{u}}\rangle
\\
+\langle\mathcal{A}_{T_h,C_h}(\mathbf{u})-\mathcal{A}_{T,C}(\mathbf{u}),\zeta_{\mathbf{u}}\rangle
+\int_\Omega(\mathsf{p}-\mathsf{q}_h)\operatorname{div}\zeta_{\mathbf{u}}\,\mathrm{d}x
+\langle\mathbf{f}(T,C)-\mathbf{f}(T_h,C_h),\zeta_{\mathbf{u}}\rangle.
\end{multline*}
In what follows, we invoke the properties \eqref{eq:prop_02} for $\mathcal{A}_{T_h,C_h}$ on the left-hand side, and \eqref{eq:prop_01} for $\mathcal{A}_{T_h,C_h}$, \eqref{eq:diff_Af}, and the bound of $\|\operatorname{div}\zeta_{\mathbf{u}}\|_{L^2(\Omega)}$ from Step 2 on the right-hand side. Thus, we obtain that
$
\nu_{-}\varepsilon_{-}\|\nabla\zeta_{\mathbf{u}}\|_{\mathbf{L}^2(\Omega)}
\leq
\nu_{+}\varepsilon_{+}\|\nabla\eta_{\mathbf{u}}\|_{\mathbf{L}^2(\Omega)}
+\sqrt{d}\,\|\mathsf{p}-\mathsf{q}_h\|_{L^2(\Omega)}
+\mathfrak{R}\left(\|\nabla e_T\|_{\mathbf{L}^2(\Omega)}+\|\nabla e_C\|_{\mathbf{L}^2(\Omega)}\right)$.
To control the last term, we use \eqref{eq:error_TC}, divided by $\mathfrak{M}>0$, the triangle inequality $\|\nabla\mathbf{e}_{\mathbf{u}}\|_{\mathbf{L}^2(\Omega)}\leq\|\nabla\eta_{\mathbf{u}}\|_{\mathbf{L}^2(\Omega)}+\|\nabla\zeta_{\mathbf{u}}\|_{\mathbf{L}^2(\Omega)}$, and the definition of $\mathfrak{N}$ to conclude that
\begin{multline}\label{eq:error_u}
\mathfrak{N}\|\nabla\mathbf{e}_{\mathbf{u}}\|_{\mathbf{L}^2(\Omega)}
\leq
\left(\nu_{-}\varepsilon_{-}+\nu_{+}\varepsilon_{+}\right)\|\nabla\eta_{\mathbf{u}}\|_{\mathbf{L}^2(\Omega)}
+\sqrt{d}\,\|\mathsf{p}-\mathsf{q}_h\|_{L^2(\Omega)}
\\
+\mathfrak{M}^{-1}\mathfrak{R}\,\mathfrak{K}\left(\|\nabla\eta_T\|_{\mathbf{L}^2(\Omega)}+\|\nabla\eta_C\|_{\mathbf{L}^2(\Omega)}\right).
\end{multline}

\emph{Step 4.} \emph{The quasi-best estimate \eqref{eq:quasi_best}.} Finally, we divide \eqref{eq:error_u} by $\mathfrak{N}>0$, and we use the resulting bound for $\|\nabla\mathbf{e}_{\mathbf{u}}\|_{\mathbf{L}^2(\Omega)}$ in \eqref{eq:error_TC}, divided by $\mathfrak{M}>0$, and then in \eqref{eq:error_p}. Since $\mathbf{w}_h$, $\mathsf{q}_h$, $S_h$, and $R_h$ are arbitrary, we can conclude \eqref{eq:quasi_best}.
\end{proof}

In two dimensions, the regularity requirements of Theorem~\ref{th:quasi_best_3D} can be relaxed.%

\begin{theorem}[quasi-best approximation for $d=2$]\label{th:quasi_best_2D}
In the framework of Theorem~\ref{th:existence_discrete}, let $d=2$, let $\ell>2$, and let $\sigma$ be such that $1/\sigma+1/\ell=1/2$. Let $(\mathbf{u},\mathsf{p},T,C)\in[\mathbf{W}^{1,\ell}(\Omega)\cap\mathbf{H}_0^1(\Omega)]\times L_0^2(\Omega)\times[W^{1,\ell}(\Omega)\cap H_0^1(\Omega)]\times[W^{1,\ell}(\Omega)\cap H_0^1(\Omega)]$ be a solution of the coupled problem \eqref{eq:weak_coupled} that satisfies the smallness conditions \eqref{eq:small_data_coupled_2D}--\eqref{eq:small_data_coupled_2D_2}, with $(\mathbf{u}_1,T_1,C_1)$ replaced by $(\mathbf{u},T,C)$. For every $\mathscr{T}_h\in\mathbb{T}$, let $(\mathbf{u}_h,\mathsf{p}_h,T_h,C_h)\in\mathbf{X}_h\times M_h\times Y_h\times Y_h$ be a solution of \eqref{eq:discrete_coupled}. Then, the bound \eqref{eq:quasi_best} holds.
\end{theorem}

\begin{proof}
The proof is a suitable adaptation of the proof of Theorem~\ref{th:quasi_best_3D}.
\end{proof}

\subsection{A priori error estimates}
In this section, as a consequence of Theorems~\ref{th:quasi_best_3D} and \ref{th:quasi_best_2D}, we derive convergence rates for the Taylor--Hood and mini elements under high regularity assumptions on the solution as in \cite[Sections 4.2.4 and 4.2.5]{Guermond-Ern}.

\begin{theorem}[convergence rates]\label{th:rates}
Let the assumptions of Theorem~\ref{th:quasi_best_3D} hold if $d=3$, or those of Theorem~\ref{th:quasi_best_2D} if $d=2$. Let $\mathcal{E}(\mathbf{u},\mathsf{p},T,C):=\|\nabla\mathbf{e}_{\mathbf{u}}\|_{\mathbf{L}^2(\Omega)}
+\|e_{\mathsf{p}}\|_{L^2(\Omega)}
+\|\nabla e_T\|_{\mathbf{L}^2(\Omega)}
+\|\nabla e_C\|_{\mathbf{L}^2(\Omega)}$. Then, we have the following estimates, where the hidden constants do not depend on $h$.

$\checkmark$ Taylor--Hood element. If $\mathbf{u}\in\mathbf{H}^{3}(\Omega)$, $\mathsf{p}\in H^{2}(\Omega)$, and $T,C\in H^{3}(\Omega)$, then
\[
\mathcal{E}(\mathbf{u},\mathsf{p},T,C)
\lesssim
h^{2}(\|\mathbf{u}\|_{\mathbf{H}^{3}(\Omega)}+\|\mathsf{p}\|_{H^{2}(\Omega)}+\|T\|_{H^{3}(\Omega)}+\|C\|_{H^{3}(\Omega)}).
\]

$\checkmark$  Mini element. If $\mathbf{u}\in\mathbf{H}^{2}(\Omega)$, $\mathsf{p}\in H^{1}(\Omega)$, and $T,C\in H^{2}(\Omega)$, then
\[
\mathcal{E}(\mathbf{u},\mathsf{p},T,C)
\lesssim
h(\|\mathbf{u}\|_{\mathbf{H}^{2}(\Omega)}+\|\mathsf{p}\|_{H^{1}(\Omega)}+\|T\|_{H^{2}(\Omega)}+\|C\|_{H^{2}(\Omega)}).
\]
\end{theorem}

\begin{proof}
The proof consists of an application of Theorems \ref{th:quasi_best_3D} and \ref{th:quasi_best_2D}.
\end{proof}

\section{Numerical experiments}
\label{sec:numerics}
In this section, we present numerical experiments in two and three dimensions to illustrate the performance of the proposed finite element methods. The computations are performed with FreeFem++ and MUMPS. We solve the discrete problem \eqref{eq:discrete_coupled} using the scheme presented in Algorithm~\ref{alg:fixed_point} and the general framework described in Section \ref{sec:fem_spaces}.

\begin{algorithm}[t]
\caption{Fixed-point iteration}
\label{alg:fixed_point}
\textbf{Input:} Mesh $\mathscr{T}_h\in\mathbb{T}$ and initial guess $(\mathbf{u}_h^0,T_h^0,C_h^0)\in\mathbf{X}_h\times Y_h\times Y_h$.
\\
$\mathbf{1}$: For $k\geq0$, let $\eta_h^k:=\nu(T_h^k,C_h^k)\,\varepsilon(x,|\nabla\mathbf{u}_h^k|)$, $\kappa_h^k:=\kappa(T_h^k,C_h^k)$, and $\mathcal{D}_h^k:=\mathcal{D}(T_h^k,C_h^k)$. Find $(\mathbf{u}_h^{k+1},\mathsf{p}_h^{k+1})\in\mathbf{X}_h\times M_h$ such that
\begin{multline*}
\int_\Omega\eta_h^k\nabla\mathbf{u}_h^{k+1}:\nabla\mathbf{v}_h\,\mathrm{d}x-\int_\Omega\mathsf{p}_h^{k+1}\operatorname{div}\mathbf{v}_h\,\mathrm{d}x=\langle\mathbf{f}(T_h^k,C_h^k),\mathbf{v}_h\rangle,
\\
\int_\Omega\mathsf{q}_h\operatorname{div}\mathbf{u}_h^{k+1}\,\mathrm{d}x=0
\qquad\forall(\mathbf{v}_h,\mathsf{q}_h)\in\mathbf{X}_h\times M_h.
\end{multline*}
$\mathbf{2}$: Find $(T_h^{k+1},C_h^{k+1})\in Y_h\times Y_h$ such that, for all $(S_h,R_h)\in Y_h\times Y_h$,
\begin{align*}
\int_\Omega\left(\kappa_h^k\nabla T_h^{k+1}\cdot\nabla S_h+\tfrac{1}{2}\left[(\mathbf{u}_h^{k+1}\cdot\nabla T_h^{k+1})S_h-(\mathbf{u}_h^{k+1}\cdot\nabla S_h)T_h^{k+1}\right]\right)\mathrm{d}x&=\langle\mathsf{g}_1,S_h\rangle,
\\
\int_\Omega\left(\mathcal{D}_h^k\nabla C_h^{k+1}\cdot\nabla R_h+\tfrac{1}{2}\left[(\mathbf{u}_h^{k+1}\cdot\nabla C_h^{k+1})R_h-(\mathbf{u}_h^{k+1}\cdot\nabla R_h)C_h^{k+1}\right]\right)\mathrm{d}x&=\langle\mathsf{g}_2,R_h\rangle.
\end{align*}
$\mathbf{3}$: If a stopping criterion, for example \eqref{eq:numerics_picard_stop}, is satisfied, stop. Otherwise, set $k\leftarrow k+1$ and go to step $1$.
\end{algorithm}

In the following, we prove the convergence of Algorithm~\ref{alg:fixed_point} under suitable smallness conditions. In addition to the assumptions of Section~\ref{sec:intro}, we assume that $\varepsilon(x,\cdot)$ is Lipschitz continuous on $[0,\infty)$ with constant $\mathcal{L}_{\varepsilon}$, i.e., $|\varepsilon(x,s_1)-\varepsilon(x,s_2)|\leq\mathcal{L}_{\varepsilon}|s_1-s_2|$ for all $s_1,s_2\geq0$ and a.e. $x\in\Omega$.

\begin{theorem}[convergence of fixed-point strategy]\label{th:algorithm_3D}
In the framework of Theorem~\ref{th:existence_discrete}, let $d=3$, assume that $\varepsilon(x,\cdot)$ is Lipschitz continuous with constant $\mathcal{L}_{\varepsilon}$, and let $(\mathbf{u}_h,\mathsf{p}_h,T_h,C_h)\in\mathbf{X}_h\times M_h\times Y_h\times Y_h$ be a solution of \eqref{eq:discrete_coupled} that satisfies \eqref{eq:small_data_discrete_2}, with $(T_{h,1},C_{h,1})$ replaced by $(T_h,C_h)$. Let $\alpha_h:=\nu_{+}\mathcal{L}_{\varepsilon}\|\nabla\mathbf{u}_h\|_{\mathbf{L}^\infty(\Omega)}(\nu_{-}\varepsilon_{-})^{-1}$ and $\rho_h:=\theta_h+\mathcal{C}_{4\hookrightarrow2}^2\,\Gamma(\mathsf{g}_1,\mathsf{g}_2)\big(C_{\mathbf{f}}\,\mathcal{C}_{2\hookrightarrow2}^2+\varepsilon_{+}\mathcal{L}_{\nu}\mathcal{C}_{6\hookrightarrow2}\|\nabla\mathbf{u}_h\|_{\mathbf{L}^3(\Omega)}\big)(\nu_{-}\varepsilon_{-})^{-1}$, and assume that
\begin{align}\label{eq:algorithm_small}
\rho_h+\alpha_h(1-\theta_h)<1.
\end{align}
Let $\{(\mathbf{u}_h^k,\mathsf{p}_h^k,T_h^k,C_h^k)\}_{k\geq1}$ be the sequence generated by steps $1$ and $2$ of Algorithm~\ref{alg:fixed_point}. Let $\mathcal{E}_h^k:=\|\nabla(\mathbf{u}_h-\mathbf{u}_h^{k})\|_{\mathbf{L}^2(\Omega)}+\|\mathsf{p}_h-\mathsf{p}_h^{k}\|_{L^2(\Omega)}+\|\nabla(T_h-T_h^{k})\|_{\mathbf{L}^2(\Omega)}+\|\nabla(C_h-C_h^{k})\|_{\mathbf{L}^2(\Omega)}$ for $k\geq1$, and let $\mathcal{E}_h^0:=\|\nabla(\mathbf{u}_h-\mathbf{u}_h^{0})\|_{\mathbf{L}^2(\Omega)}+\|\nabla(T_h-T_h^{0})\|_{\mathbf{L}^2(\Omega)}+\|\nabla(C_h-C_h^{0})\|_{\mathbf{L}^2(\Omega)}$. Then, there exists $q_h\in[0,1)$ such that, for every $k\geq0$, we have
\begin{equation}\label{eq:algorithm_rate}
\mathcal{E}_h^{k+1}\lesssim q_h^{k}\,\mathcal{E}_h^{0},
\end{equation}
where $q_h$ and the hidden constant do not depend on $k$ or on the initial guess. In particular, the sequence $\{(\mathbf{u}_h^k,\mathsf{p}_h^k,T_h^k,C_h^k)\}_{k\geq1}$ converges to $(\mathbf{u}_h,\mathsf{p}_h,T_h,C_h)$ as $k\uparrow\infty$.
\end{theorem}

\begin{proof}
To simplify the presentation, we introduce the notation $\nu_h:=\nu(T_h,C_h)$, $\kappa_h:=\kappa(T_h,C_h)$, $\mathcal{D}_h:=\mathcal{D}(T_h,C_h)$, $\varepsilon_h:=\varepsilon(x,|\nabla\mathbf{u}_h|)$, $\nu_h^k:=\nu(T_h^k,C_h^k)$, $\varepsilon_h^k:=\varepsilon(x,|\nabla\mathbf{u}_h^k|)$, $\kappa_h^k:=\kappa(T_h^k,C_h^k)$, $\mathbf{e}_{\mathbf{u}}^{k}:=\mathbf{u}_h-\mathbf{u}_h^{k}$, $e_{\mathsf{p}}^{k}:=\mathsf{p}_h-\mathsf{p}_h^{k}$, $e_T^{k}:=T_h-T_h^{k}$, and $e_C^{k}:=C_h-C_h^{k}$. We also define $X^k:=\nu_{-}\varepsilon_{-}\|\nabla\mathbf{e}_{\mathbf{u}}^k\|_{\mathbf{L}^2(\Omega)}$, $E^k:=\|\nabla e_T^k\|_{\mathbf{L}^2(\Omega)}+\|\nabla e_C^k\|_{\mathbf{L}^2(\Omega)}$, $\gamma_h:=(\nu_{-}\varepsilon_{-})^{-1}\mathcal{C}_{4\hookrightarrow2}^2\,\Gamma(\mathsf{g}_1,\mathsf{g}_2)$, and $\mathfrak{R}_h:=C_{\mathbf{f}}\,\mathcal{C}_{2\hookrightarrow2}^2+\varepsilon_{+}\mathcal{L}_{\nu}\mathcal{C}_{6\hookrightarrow2}\|\nabla\mathbf{u}_h\|_{\mathbf{L}^3(\Omega)}$, so that $\rho_h=\theta_h+\gamma_h\mathfrak{R}_h$. We note that \eqref{eq:algorithm_small} implies $\rho_h<1$, since $\alpha_h\geq0$ and $\theta_h<1$. Since $\gamma_h\mathfrak{R}_h=\rho_h-\theta_h$, the bound $\rho_h<1$ is \eqref{eq:small_data_discrete}, and thus, by Theorem~\ref{th:uniqueness_discrete_3D}, $(\mathbf{u}_h,\mathsf{p}_h,T_h,C_h)$ is the unique solution of \eqref{eq:discrete_coupled}. We divide the proof into three steps.

\emph{Step 1.} \emph{Temperature and concentration.} We subtract the temperature equations of \eqref{eq:discrete_coupled} and step $2$ of Algorithm~\ref{alg:fixed_point} with the test function $e_T^{k+1}$, and use the identity $\kappa_h\nabla T_h-\kappa_h^k\nabla T_h^{k+1}=(\kappa_h-\kappa_h^k)\nabla T_h+\kappa_h^k\nabla e_T^{k+1}$. On the other hand, we write $\mathbf{u}_h=\mathbf{e}_{\mathbf{u}}^{k+1}+\mathbf{u}_h^{k+1}$ and $T_h=e_T^{k+1}+T_h^{k+1}$. Thus, we arrive at%
\begin{multline*}
\int_\Omega\kappa_h^k|\nabla e_T^{k+1}|^2\,\mathrm{d}x=-\int_\Omega(\kappa_h-\kappa_h^k)\nabla T_h\cdot\nabla e_T^{k+1}\,\mathrm{d}x
\\
-\frac{1}{2}\int_\Omega\left((\mathbf{e}_{\mathbf{u}}^{k+1}\cdot\nabla T_h)e_T^{k+1}-(\mathbf{e}_{\mathbf{u}}^{k+1}\cdot\nabla e_T^{k+1})T_h\right)\mathrm{d}x.
\end{multline*}
We now use the bound $\kappa_h^k\geq\kappa_{-}$, the Lipschitz continuity of $\kappa$, H\"older's inequality, the Sobolev embeddings $H_0^1(\Omega)\hookrightarrow L^6(\Omega)$ and $H_0^1(\Omega)\hookrightarrow L^4(\Omega)$, and the estimate $\|\nabla T_h\|_{\mathbf{L}^2(\Omega)}\leq\kappa_{-}^{-1}\|\mathsf{g}_1\|_{H^{-1}(\Omega)}$ from \eqref{eq:estimate_TC_discrete} to obtain
\[
\|\nabla e_T^{k+1}\|_{\mathbf{L}^2(\Omega)}\leq\frac{\mathcal{L}_{\kappa}\mathcal{C}_{6\hookrightarrow2}\|\nabla T_h\|_{\mathbf{L}^3(\Omega)}}{\kappa_{-}}E^k+\frac{\mathcal{C}_{4\hookrightarrow2}^2\|\mathsf{g}_1\|_{H^{-1}(\Omega)}}{\kappa_{-}^{2}}\|\nabla\mathbf{e}_{\mathbf{u}}^{k+1}\|_{\mathbf{L}^2(\Omega)}.
\]
Similar arguments, applied to the concentration equations with $R_h=e_C^{k+1}$ and $\mathcal{D}_h$ instead of $\kappa_h$, give a bound for $\|\nabla e_C^{k+1}\|_{\mathbf{L}^2(\Omega)}$. If we add these two estimates and use the definitions of $\theta_h$, $\Gamma(\mathsf{g}_1,\mathsf{g}_2)$, and $\gamma_h$, we can conclude that
\begin{equation}\label{eq:algorithm_TC}
E^{k+1}\leq\theta_h E^k+\gamma_h X^{k+1}.
\end{equation}

\emph{Step 2.} \emph{Velocity.} In view of the second equations of \eqref{eq:discrete_coupled} and step $1$ of Algorithm~\ref{alg:fixed_point}, we have $\mathbf{e}_{\mathbf{u}}^{k+1}\in\mathbf{V}_h$. We subtract the first equations of \eqref{eq:discrete_coupled} and step $1$, and use the identity
\[
\nu_h\varepsilon_h\nabla\mathbf{u}_h-\eta_h^k\nabla\mathbf{u}_h^{k+1}=\eta_h^k\nabla\mathbf{e}_{\mathbf{u}}^{k+1}+(\nu_h-\nu_h^k)\varepsilon_h\nabla\mathbf{u}_h+\nu_h^k(\varepsilon_h-\varepsilon_h^k)\nabla\mathbf{u}_h.
\]
Thus, we obtain the following identity
\begin{multline}\label{eq:algorithm_identity}
\int_\Omega\eta_h^k\nabla\mathbf{e}_{\mathbf{u}}^{k+1}:\nabla\mathbf{v}_h\,\mathrm{d}x-\int_\Omega e_{\mathsf{p}}^{k+1}\operatorname{div}\mathbf{v}_h\,\mathrm{d}x=\langle\mathbf{f}(T_h,C_h)-\mathbf{f}(T_h^k,C_h^k),\mathbf{v}_h\rangle
\\
-\int_\Omega(\nu_h-\nu_h^k)\varepsilon_h\nabla\mathbf{u}_h:\nabla\mathbf{v}_h\,\mathrm{d}x
-\int_\Omega\nu_h^k(\varepsilon_h-\varepsilon_h^k)\nabla\mathbf{u}_h:\nabla\mathbf{v}_h\,\mathrm{d}x
\end{multline}
for all $\mathbf{v}_h\in\mathbf{X}_h$. Let $\mathbf{v}_h=\mathbf{e}_{\mathbf{u}}^{k+1}$, and note that $\int_\Omega e_{\mathsf{p}}^{k+1}\operatorname{div}\mathbf{e}_{\mathbf{u}}^{k+1}\,\mathrm{d}x=0$, due to $e_{\mathsf{p}}^{k+1}\in M_h$. We now use the bound $\eta_h^k\geq\nu_{-}\varepsilon_{-}$, the Lipschitz continuity of $\mathbf{f}_1$, $\mathbf{f}_2$, $\nu$, and $\varepsilon(x,\cdot)$, the bound $C_{\mathbf{f}_1}\|\nabla e_T^k\|_{\mathbf{L}^2(\Omega)}+C_{\mathbf{f}_2}\|\nabla e_C^k\|_{\mathbf{L}^2(\Omega)}\leq C_{\mathbf{f}}E^k$, the bounds $\nu_h^k\leq\nu_{+}$, $|\varepsilon(x,|\boldsymbol{\tau}|)\boldsymbol{\tau}|\leq\varepsilon_{+}|\boldsymbol{\tau}|$, and $\big||\nabla\mathbf{u}_h|-|\nabla\mathbf{u}_h^k|\big|\leq|\nabla\mathbf{e}_{\mathbf{u}}^k|$, H\"older's inequality, and the Sobolev embeddings $H_0^1(\Omega)\hookrightarrow L^2(\Omega)$ and $H_0^1(\Omega)\hookrightarrow L^6(\Omega)$ to obtain
\begin{equation}\label{eq:algorithm_u}
X^{k+1}\leq\alpha_h X^k+\mathfrak{R}_h E^k.
\end{equation}

\emph{Step 3.} \emph{Convergence.} If we use \eqref{eq:algorithm_u} in \eqref{eq:algorithm_TC}, we obtain $E^{k+1}\leq\gamma_h\alpha_h X^k+\rho_h E^k$. On the other hand, if $\mu:=\mathfrak{R}_h(1-\rho_h)^{-1}$ we have that $\alpha_h(1+\gamma_h\mu)<1$. Thus, there exists $\mu_h>\mathfrak{R}_h(1-\rho_h)^{-1}$ such that $\alpha_h(1+\gamma_h\mu_h)<1$. In the following, we define $Z^k:=X^k+\mu_h E^k$ and $q_h:=\max\{\alpha_h(1+\gamma_h\mu_h),\rho_h+\mathfrak{R}_h\mu_h^{-1}\}$, and note that $q_h<1$ and that $q_h$ does not depend on $k$. These estimates yield
\[
Z^{k+1}\leq\alpha_h(1+\gamma_h\mu_h)X^k+(\mathfrak{R}_h+\mu_h\rho_h)E^k\leq q_h Z^k,\quad Z^{k}\to 0 \text{ as } k\uparrow\infty.
\]
Finally, we use the identity \eqref{eq:algorithm_identity}, the bound $\eta_h^k\leq\nu_{+}\varepsilon_{+}$, the arguments developed in Step 2, and the discrete inf-sup condition \eqref{eq:infsup_disc} to arrive at
\[
\tilde{\beta}\|e_{\mathsf{p}}^{k+1}\|_{L^2(\Omega)}\leq\nu_{+}\varepsilon_{+}\|\nabla\mathbf{e}_{\mathbf{u}}^{k+1}\|_{\mathbf{L}^2(\Omega)}+\alpha_h X^k+\mathfrak{R}_h E^k.
\]
We note that, since $X^{k+1}\leq Z^{k+1}\leq Z^k$, $E^{k+1}\leq\mu_h^{-1}Z^k$, $X^k\leq Z^k$, $E^k\leq\mu_h^{-1}Z^k$, and the bound $Z^0\leq \max\{\nu_{-}\varepsilon_{-},\mu_h\}\mathcal{E}_h^0$, we can conclude that $\mathcal{E}_h^{k+1}\lesssim Z^k\leq q_h^k Z^0\lesssim q_h^k\mathcal{E}_h^0$, which is the estimate presented in \eqref{eq:algorithm_rate}. This completes the proof.
\end{proof}

\begin{theorem}[convergence of fixed-point strategy for $d=2$]\label{th:algorithm_2D}
In the framework of Theorem~\ref{th:existence_discrete}, let $d=2$, let $\ell>2$, let $\sigma$ be such that $1/\sigma+1/\ell=1/2$, and assume that $\varepsilon(x,\cdot)$ is Lipschitz continuous with constant $\mathcal{L}_{\varepsilon}$. Let $(\mathbf{u}_h,\mathsf{p}_h,T_h,C_h)\in\mathbf{X}_h\times M_h\times Y_h\times Y_h$ be a solution of \eqref{eq:discrete_coupled} that satisfies \eqref{eq:small_data_discrete_2D_2}, with $(T_{h,1},C_{h,1})$ replaced by $(T_h,C_h)$, and \eqref{eq:algorithm_small}, with $\theta_h$ and $\rho_h$ replaced by $\theta_{h,\ell}$ and
\[
\rho_{h,\ell}:=\theta_{h,\ell}+\frac{\mathcal{C}_{4\hookrightarrow2}^2\,\Gamma(\mathsf{g}_1,\mathsf{g}_2)}{\nu_{-}\varepsilon_{-}}\left(C_{\mathbf{f}}\,\mathcal{C}_{2\hookrightarrow2}^2+\varepsilon_{+}\mathcal{L}_{\nu}\mathcal{C}_{\sigma\hookrightarrow2}\|\nabla\mathbf{u}_h\|_{\mathbf{L}^{\ell}(\Omega)}\right).
\]
Then, the conclusions of Theorem~\ref{th:algorithm_3D} hold.
\end{theorem}

\begin{proof}
The proof is a suitable adaptation of the proof of Theorem~\ref{th:algorithm_3D}.
\end{proof}

In the following, we describe the setting of our numerical experiments. In the manufactured-solution tests, the pressure mean is imposed by a global Lagrange multiplier in the velocity--pressure block system. In the applications, the natural outlet condition fixes the pressure reference. For the finest three-dimensional Taylor--Hood convergence mesh, the factorization of this block is stored on disk. We use continuous piecewise quadratic functions for the temperature and concentration with Taylor--Hood, and continuous piecewise linear functions with the mini element.

Experiments~\ref{sec:numerics_convergence} and~\ref{sec:numerics_robustness} use zero initial guesses, including $\mathsf{p}_h^0=0$. At each Picard step, we freeze the complete factor $\nu(T_h^k,C_h^k)\varepsilon(x,|\nabla\mathbf{u}_h^k|)$ and the buoyancy at the preceding iterate, and solve the linear velocity--pressure system. We then solve the two scalar equations with the new velocity and the preceding diffusivities, using the skew-symmetric convection form. In the applications, $C$ does not affect the flow and the temperature, and the scalar transport forms are specified below. For the manufactured-solution tests, with $z_h^k=(\mathbf{u}_h^k,\mathsf{p}_h^k,T_h^k,C_h^k)$, the stopping criterion is
\begin{equation}\label{eq:numerics_picard_stop}
\|z_h^{k+1}-z_h^k\|_\star/\max\{1,\|z_h^{k+1}\|_\star\}<10^{-6},
\end{equation}
where $\|(\mathbf{v},\mathsf{q},S,R)\|_\star^2
=\|\nabla\mathbf{v}\|_{\mathbf{L}^2(\Omega)}^2+\|\mathsf{q}\|_{L^2(\Omega)}^2+\|\nabla S\|_{\mathbf{L}^2(\Omega)}^2+\|\nabla R\|_{\mathbf{L}^2(\Omega)}^2$. On the other hand, when a reference solution is available, we report the absolute errors as:
$E_{\mathbf{u}}:=\|\nabla\mathbf{e}_{\mathbf{u}}\|_{\mathbf{L}^2(\Omega)}$,
$E_{\mathsf{p}}:=\|e_{\mathsf{p}}\|_{L^2(\Omega)}$,
$E_T:=\|\nabla e_T\|_{\mathbf{L}^2(\Omega)}$, and $E_C:=\|\nabla e_C\|_{\mathbf{L}^2(\Omega)}$, with $\mathbf{e}_{\mathbf{u}}$, $e_{\mathsf{p}}$, $e_T$, and $e_C$ as in Section~\ref{sec:quasi_best}.
For any of these errors, if $E_h$ and $E_{\widetilde h}$ denote its values on two consecutive meshes of sizes $h$ and $\widetilde h$, the experimental convergence rate is defined by
$r=\log(E_h/E_{\widetilde h})/\log(h/\widetilde h)$, denoted by $r_{\mathbf u}$, $r_{\mathsf p}$, $r_T$, and $r_C$, respectively.
The rates are computed from unrounded errors. We denote by $\mathrm{dof}$ the total number of degrees of freedom, including the Dirichlet degrees of freedom before elimination and, in the manufactured-solution tests, the global pressure Lagrange multiplier.
We denote by $\mathrm{it}$ the number of Picard iterations.

\subsection[Convergence in 2D and 3D]{Convergence in 2D and 3D}
\label{sec:numerics_convergence}
We study the convergence of the Taylor--Hood and mini approximations on $\Omega=(0,1)^d$, with $d\in\{2,3\}$. The meshes have $N$ subdivisions in each coordinate direction and maximum diameter $h=\sqrt{d}/N$.

The data for the model are chosen as follows: $\varepsilon(x,s)=1+0.5/\sqrt{1+s^2}$, $\nu(T,C)=1+0.20\tanh(T+C)$, $\kappa(T,C)=1+0.15\tanh(T-2C)$, $\mathcal{D}(T,C)=0.70\bigl(1+0.10\tanh(2T+C)\bigr)$,
and $\mathbf{f}(T,C)=\mathbf{f}_0+(T-C)\mathbf{e}_d$, where $\mathbf{e}_d$ is the last vector of the canonical basis of $\mathbb{R}^d$. We manufacture the data $\mathbf{f}_0$, $\mathsf{g}_1$, and $\mathsf{g}_2$ from the exact solutions below.

In two dimensions, we take $N\in\{8,16,32,64\}$ and $\psi(x,y)=0.05\sin^2(\pi x)\sin^2(\pi y)$, $\mathbf{u}=(\partial_y\psi,-\partial_x\psi)$, $\mathsf{p}(x,y)=10\cos(\pi x)\cos(\pi y)$, $T(x,y)=0.1\sin(\pi x)\sin(2\pi y)$, and $C(x,y)=0.08\sin(2\pi x)\sin(\pi y)$.
\begin{table}[tbp]
\centering\footnotesize
\setlength{\tabcolsep}{3pt}
\caption{Experiment \ref{sec:numerics_convergence}. Two-dimensional errors, experimental convergence rates, and Picard iteration counts for the Taylor--Hood and mini elements.}
\label{tab:convergence_2d}
\begin{tabular}{|rr|cc|cc|cc|cc|r|}
\hline\hline
$\mathrm{dof}$ & $h$ & $E_{\mathbf u}$ & $r_{\mathbf u}$ & $E_{\mathsf p}$ & $r_{\mathsf p}$ & $E_T$ & $r_T$ & $E_C$ & $r_C$ & $\mathrm{it}$ \\
\hline
\multicolumn{11}{|c|}{Taylor--Hood} \\
\hline
1238 & 0.17678 & {3.3577e-02} & -- & {6.6409e-02} & -- & 1.1986e-02 & -- & 9.5887e-03 & -- & 7 \\
4646 & 0.08839 & {8.1278e-03} & ${2.047}$ & {1.6205e-02} & ${2.035}$ & 3.0505e-03 & $1.974$ & 2.4404e-03 & $1.974$ & 7 \\
17990 & 0.04419 & {2.0125e-03} & ${2.014}$ & {4.0250e-03} & ${2.009}$ & 7.6621e-04 & $1.993$ & 6.1297e-04 & $1.993$ & 7 \\
70790 & 0.02210 & {5.0182e-04} & ${2.004}$ & {1.0046e-03} & ${2.002}$ & 1.9178e-04 & $1.998$ & 1.5342e-04 & $1.998$ & 7 \\
\hline
\multicolumn{11}{|c|}{Mini} \\
\hline
662 & 0.17678 & {2.1244e-01} & -- & {1.5593e-01} & -- & 1.0017e-01 & -- & 8.0135e-02 & -- & 7 \\
2470 & 0.08839 & {1.0608e-01} & ${1.002}$ & {4.7445e-02} & ${1.717}$ & 5.0876e-02 & $0.977$ & 4.0701e-02 & $0.977$ & 7 \\
9542 & 0.04419 & {5.2909e-02} & ${1.004}$ & {1.5488e-02} & ${1.615}$ & 2.5539e-02 & $0.994$ & 2.0431e-02 & $0.994$ & 7 \\
37510 & 0.02210 & {2.6408e-02} & ${1.003}$ & {5.2958e-03} & ${1.548}$ & 1.2782e-02 & $0.999$ & 1.0226e-02 & $0.999$ & 7 \\
\hline
\hline
\end{tabular}
\end{table}

Table~\ref{tab:convergence_2d} shows rates close to $\mathcal{O}(h^{2})$ for all Taylor--Hood errors. On the last refinement, the velocity and scalar rates are ${2.004}$ and $1.998$, respectively, and the pressure rate is ${2.002}$. The pressure rates also remain close to two on all refinements. For the mini element, the velocity and scalar errors approach first-order convergence, whereas the last pressure rate is $1.548$. This superconvergence for the mini element has been reported in \cite{CIONCOLINI20192432}.

In three dimensions, we take $N\in\{4,8,16,32\}$, set $\mathbf{a}=(1,2,3)^{\mathsf{T}}$, and choose the functions $\psi(x,y,z)=0.02\sin^2(\pi x)\sin^2(\pi y)\sin^2(\pi z)$, $\mathbf{u}=\nabla\times(\psi\mathbf{a})$, $\mathsf{p}(x,y,z)=10\cos(\pi x)\cos(\pi y)\cos(\pi z)$, $T(x,y,z)=0.1\sin(\pi x)\sin(2\pi y)\sin(\pi z)$, and $C(x,y,z)=0.08\sin(2\pi x)\sin(\pi y)\sin(\pi z)$.

\begin{table}[tbp]
\centering\footnotesize
\setlength{\tabcolsep}{3pt}
\caption{Experiment \ref{sec:numerics_convergence}. Three-dimensional errors, experimental convergence rates, and Picard iteration counts for the Taylor--Hood and mini elements.}
\label{tab:convergence_3d}
\begin{tabular}{|rr|cc|cc|cc|cc|r|}
\hline\hline
$\mathrm{dof}$ & $h$ & $E_{\mathbf u}$ & $r_{\mathbf u}$ & $E_{\mathsf p}$ & $r_{\mathsf p}$ & $E_T$ & $r_T$ & $E_C$ & $r_C$ & $\mathrm{it}$ \\
\hline
\multicolumn{11}{|c|}{Taylor--Hood} \\
\hline
3771 & 0.43301 & 1.5230e-01 & -- & 3.2549e-01 & -- & 4.2770e-02 & -- & 3.4216e-02 & -- & 7 \\
25295 & 0.21651 & 3.6154e-02 & $2.075$ & 7.3546e-02 & $2.146$ & 1.1780e-02 & $1.860$ & 9.4238e-03 & $1.860$ & 7 \\
184599 & 0.10825 & 8.6681e-03 & $2.060$ & 1.7806e-02 & $2.046$ & 3.0404e-03 & $1.954$ & 2.4323e-03 & $1.954$ & 7 \\
1409063 & 0.05413 & 2.1370e-03 & $2.020$ & 4.4127e-03 & $2.013$ & 7.6703e-04 & $1.987$ & 6.1363e-04 & $1.987$ & 7\\
\hline
\multicolumn{11}{|c|}{Mini} \\
\hline
1903 & 0.43301 & 4.6792e-01 & -- & 1.0163e+00 & -- & 1.6488e-01 & -- & 1.3191e-01 & -- & 6 \\
13591 & 0.21651 & 2.3593e-01 & $0.988$ & 3.2588e-01 & $1.641$ & 8.9743e-02 & $0.878$ & 7.1796e-02 & $0.878$ & 7 \\
103207 & 0.10825 & 1.1735e-01 & $1.007$ & 1.0841e-01 & $1.588$ & 4.5906e-02 & $0.967$ & 3.6725e-02 & $0.967$ & 7 \\
805447 & 0.05413 & 5.8415e-02 & $1.006$ & 3.6452e-02 & $1.572$ & 2.3087e-02 & $0.992$ & 1.8469e-02 & $0.992$ & 7 \\
\hline
\hline
\end{tabular}
\end{table}

The results in Table~\ref{tab:convergence_3d} show similar convergence rates. For Taylor--Hood, the last velocity, pressure, temperature, and concentration rates are $2.020$, $2.013$, $1.987$, and $1.987$, respectively. In particular, the pressure error decreases from $1.7806\times10^{-2}$ to $4.4127\times10^{-3}$ between $N=16$ and $N=32$. For the mini element, the velocity and scalar rates approach one, and the last pressure rate is $1.572$. These velocity and scalar rates are consistent with the approximation orders in Theorem~\ref{th:rates}. The faster pressure convergence of the mini element is an observation for the mesh families considered here.

With the tolerance in \eqref{eq:numerics_picard_stop}, Picard converges in seven iterations on every mesh, except for the coarsest three-dimensional mini mesh, which requires six. The nonlinear residuals of the flow and scalar blocks, measured in the infinity norm and divided by the maximum of one and the corresponding right-hand-side norm, are below $9.61\times10^{-8}$. Thus, the selected tolerance preserves the observed convergence rates on the meshes reported here.

As a complementary check, we retain the exact velocity, pressure, and coefficient laws, and set $T(x)=0.5\prod_{j=1}^d\sin(\pi x_j)$. We take $C(x,y)=0.4\sin(\pi x)\sin(\pi y)\bigl(1+0.30\cos(\pi x)+0.20\sin(2\pi y)\bigr)$ in 2D, and
$C(x,y,z)=0.4\sin(\pi x)\sin(\pi y)\sin(\pi z)\bigl(1+0.30\cos(\pi x)+0.20\sin(2\pi y)+0.15\cos(\pi z)\bigr)$
in 3D. The source terms are recomputed for these fields. This test increases the scalar amplitudes and introduces an asymmetric concentration profile. Table~\ref{tab:convergence_complement} recovers the same approximation orders on the two meshes considered. All runs satisfy \eqref{eq:numerics_picard_stop} in six steps in 2D and seven in 3D, with normalized block residuals below $3.18\times10^{-8}$.

\begin{table}[hbt!]
\centering\footnotesize
\caption{Experiment \ref{sec:numerics_convergence}. Complementary manufactured-solution test. Rates are computed between the two meshes of each element and dimension.}
\label{tab:convergence_complement}
\begin{tabular}{|lrr|rrrr|r|}
\hline\hline
Element & $\mathrm{dof}$ & $h$ & $r_{\mathbf u}$ & $r_{\mathsf p}$ & $r_T$ & $r_C$ & $\mathrm{it}$ \\
\hline
\multicolumn{8}{|c|}{2D} \\
\hline
Taylor--Hood & 17990 & 0.04419 & -- & -- & -- & -- & 6 \\
Taylor--Hood & 70790 & 0.02210 & $2.003$ & $2.002$ & $1.998$ & $1.998$ & 6 \\
Mini & 9542 & 0.04419 & -- & -- & -- & -- & 6 \\
Mini & 37510 & 0.02210 & $1.002$ & $1.553$ & $0.999$ & $0.999$ & 6 \\
\hline
\multicolumn{8}{|c|}{3D} \\
\hline
Taylor--Hood & 25295 & 0.21651 & -- & -- & -- & -- & 7 \\
Taylor--Hood & 184599 & 0.10825 & $2.056$ & $2.046$ & $1.954$ & $1.960$ & 7 \\
Mini & 13591 & 0.21651 & -- & -- & -- & -- & 7 \\
Mini & 103207 & 0.10825 & $1.007$ & $1.591$ & $0.967$ & $0.976$ & 7 \\
\hline
\hline
\end{tabular}
\end{table}

\subsection{Sensitivity with respect to the coefficients}
\label{sec:numerics_robustness}
In this experiment, we study the sensitivity to the rheological contrast and the effects of small viscosity and diffusivities. We reuse the original exact solutions of Experiment~\ref{sec:numerics_convergence}, and recompute the source terms for each choice of coefficients. The selected mesh sizes are $h\approx0.04419,\,0.02210$ in 2D and $h\approx0.21651,\,0.10825$ in 3D. As above, Picard starts from zero and uses the stopping criterion~\eqref{eq:numerics_picard_stop}.

For comparison between fields, we use the relative errors $\widehat{E}_{\mathbf u}={E_{\mathbf u}}/{\|\nabla\mathbf u\|_{\mathbf L^2(\Omega)}},$ $\widehat{E}_{\mathsf p}={E_{\mathsf p}}/{\|\mathsf p\|_{L^2(\Omega)}},$ $\widehat{E}_T={E_T}/{\|\nabla T\|_{\mathbf L^2(\Omega)}},$ $\widehat{E}_C={E_C}/{\|\nabla C\|_{\mathbf L^2(\Omega)}}.$ The tables report these errors in percent.

\subsubsection{Rheological contrast}\label{subsec:rheological}
The aim of this test is to assess how increasing the rheological contrast affects the approximation errors and the number of Picard iterations. We compare the Taylor--Hood and mini elements in two and three dimensions. Let us replace the constitutive law of Experiment~\ref{sec:numerics_convergence} by $\varepsilon_R(x,s)=1+(R-1)/\sqrt{1+s^2}$, $R\in\{1,1.5,10,100\}$, $s\geq0$,
and retain the other nonlinear coefficients. The parameter $R$ describes the global contrast of this law, for which $1\leq\varepsilon_R(x,s)\leq R$. Note that this parameter does not need to be equal to the spatial contrast attained by the manufactured solution.
\begin{table}[hbt!]
\centering\footnotesize
\caption{Experiment \ref{subsec:rheological}. Two-dimensional errors in percent and Picard iteration counts on the finer mesh of the rheological contrast study.}
\label{tab:robustness_rheology_2d}
\begin{tabular}{|rrr|rrrr|r|}
\hline\hline
$R$ & $\mathrm{dof}$ & $h$ & $100\widehat{E}_{\mathbf u}$ & $100\widehat{E}_{\mathsf p}$ & $100\widehat{E}_T$ & $100\widehat{E}_C$ & $\mathrm{it}$ \\
\hline
\multicolumn{8}{|c|}{Taylor--Hood} \\
\hline
1 & 70790 & 0.02210 & 0.0720 & 0.0201 & 0.0546 & 0.0546 & 4 \\
1.5 & 70790 & 0.02210 & 0.0719 & 0.0201 & 0.0546 & 0.0546 & 7 \\
10 & 70790 & 0.02210 & 0.0718 & 0.0201 & 0.0546 & 0.0546 & 12 \\
100 & 70790 & 0.02210 & 0.0718 & 0.0201 & 0.0546 & 0.0546 & 15 \\
\hline
\multicolumn{8}{|c|}{Mini} \\
\hline
1 & 37510 & 0.02210 & 3.7846 & 0.0749 & 3.6392 & 3.6392 & 4 \\
1.5 & 37510 & 0.02210 & 3.7841 & 0.1059 & 3.6392 & 3.6392 & 7 \\
10 & 37510 & 0.02210 & 3.7846 & 0.6458 & 3.6392 & 3.6392 & 12 \\
100 & 37510 & 0.02210 & 3.7850 & 6.3799 & 3.6392 & 3.6392 & 15 \\
\hline
\hline
\end{tabular}
\end{table}

\begin{table}[hbt!]
\centering\footnotesize
\caption{Experiment \ref{subsec:rheological}. Three-dimensional errors in percent and Picard iteration counts on the finer mesh of the rheological contrast study.}
\label{tab:robustness_rheology_3d}
\begin{tabular}{|rrr|rrrr|r|}
\hline\hline
$R$ & $\mathrm{dof}$ & $h$ & $100\widehat{E}_{\mathbf u}$ & $100\widehat{E}_{\mathsf p}$ & $100\widehat{E}_T$ & $100\widehat{E}_C$ & $\mathrm{it}$ \\
\hline
\multicolumn{8}{|c|}{Taylor--Hood} \\
\hline
1 & 184599 & 0.10825 & 1.2513 & 0.5036 & 1.1175 & 1.1175 & 4 \\
1.5 & 184599 & 0.10825 & 1.2121 & 0.5036 & 1.1175 & 1.1175 & 7 \\
10 & 184599 & 0.10825 & 1.1731 & 0.5045 & 1.1175 & 1.1175 & 14 \\
100 & 184599 & 0.10825 & 1.1720 & 0.5850 & 1.1175 & 1.1175 & 18 \\
\hline
\multicolumn{8}{|c|}{Mini} \\
\hline
1 & 103207 & 0.10825 & 16.4253 & 2.1957 & 16.8728 & 16.8728 & 4 \\
1.5 & 103207 & 0.10825 & 16.4104 & 3.0662 & 16.8728 & 16.8729 & 7 \\
10 & 103207 & 0.10825 & 16.4255 & 18.0096 & 16.8728 & 16.8729 & 14 \\
100 & 103207 & 0.10825 & 16.4429 & $\mathbf{176.3720}$ & 16.8728 & 16.8729 & 18 \\
\hline
\hline
\end{tabular}
\end{table}

Tables~\ref{tab:robustness_rheology_2d} and~\ref{tab:robustness_rheology_3d} report the results on the finer mesh of each dimension. All 32 runs meet the stopping criterion. On these finer meshes, the iteration count increases from four at $R=1$ to fifteen in 2D and eighteen in 3D at $R=100$, for both elements. The velocity and scalar errors remain of the same magnitude as $R$ varies. The rates computed from the two meshes remain close to two for Taylor--Hood and one for mini in these three fields.

\begin{figure}[hbt!]
\centering
\includegraphics[width=0.8\linewidth]{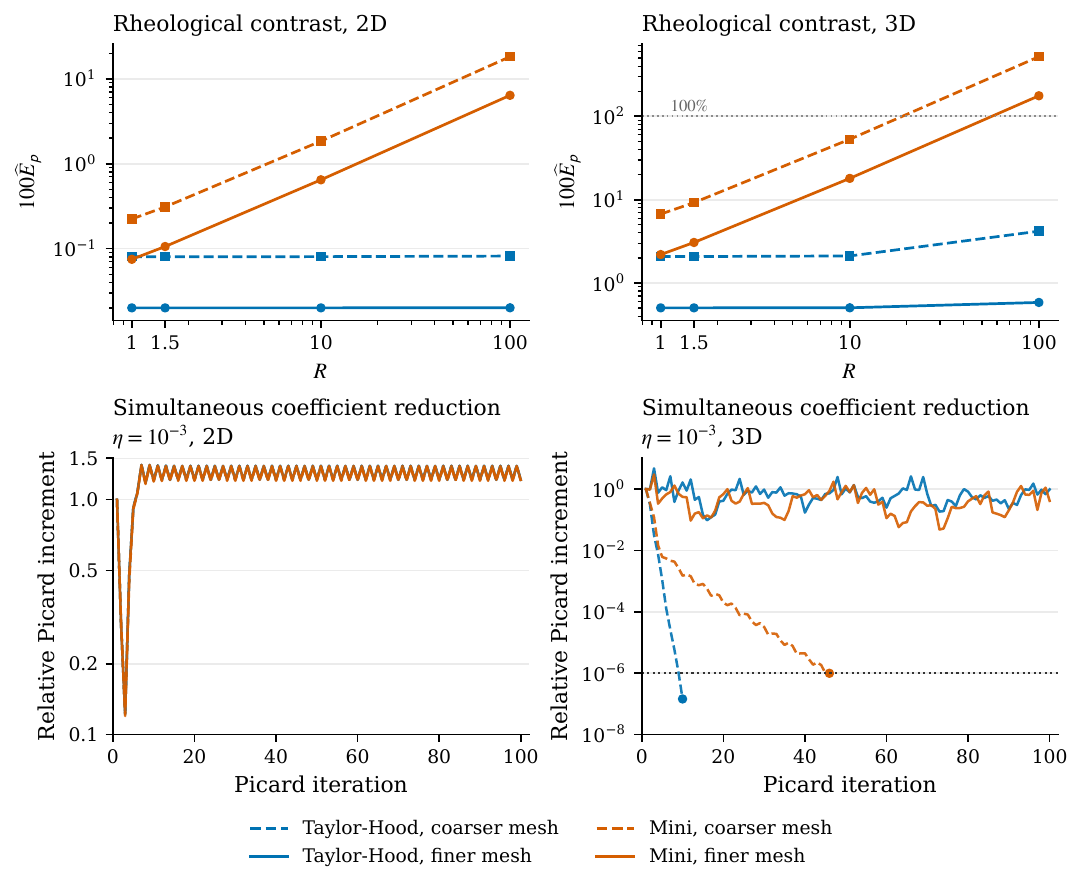}
\caption{Experiment \ref{sec:numerics_robustness}. Relative pressure errors for the rheological contrast study (top), and Picard histories for simultaneous coefficient reduction with $\eta=10^{-3}$ (bottom). Dashed and solid curves refer to the coarser and finer meshes. The dotted line in the bottom right panel marks the stopping tolerance $10^{-6}$, below the range shown in the bottom left panel.}
\label{fig:robustness_coefficients}
\end{figure}

The pressure error of the mini element is more sensitive to $R$. In 3D, its relative error increases from $2.20\%$ at $R=1$ to $176.37\%$ at $R=100$, whereas the corresponding Taylor--Hood error increases from $0.50\%$ to $0.59\%$. For mini at $R=100$, refinement reduces the pressure error from $519.09\%$ to $176.37\%$, as shown in Figure~\ref{fig:robustness_coefficients}. On the finer mesh, $E_{\mathsf p}=6.2357$ and $\|\mathsf p\|_{L^2(\Omega)}=3.5355$.

\begin{figure}[hbt!]
\centering
\begin{minipage}{\linewidth}
\centering
\includegraphics[width=0.8\linewidth]{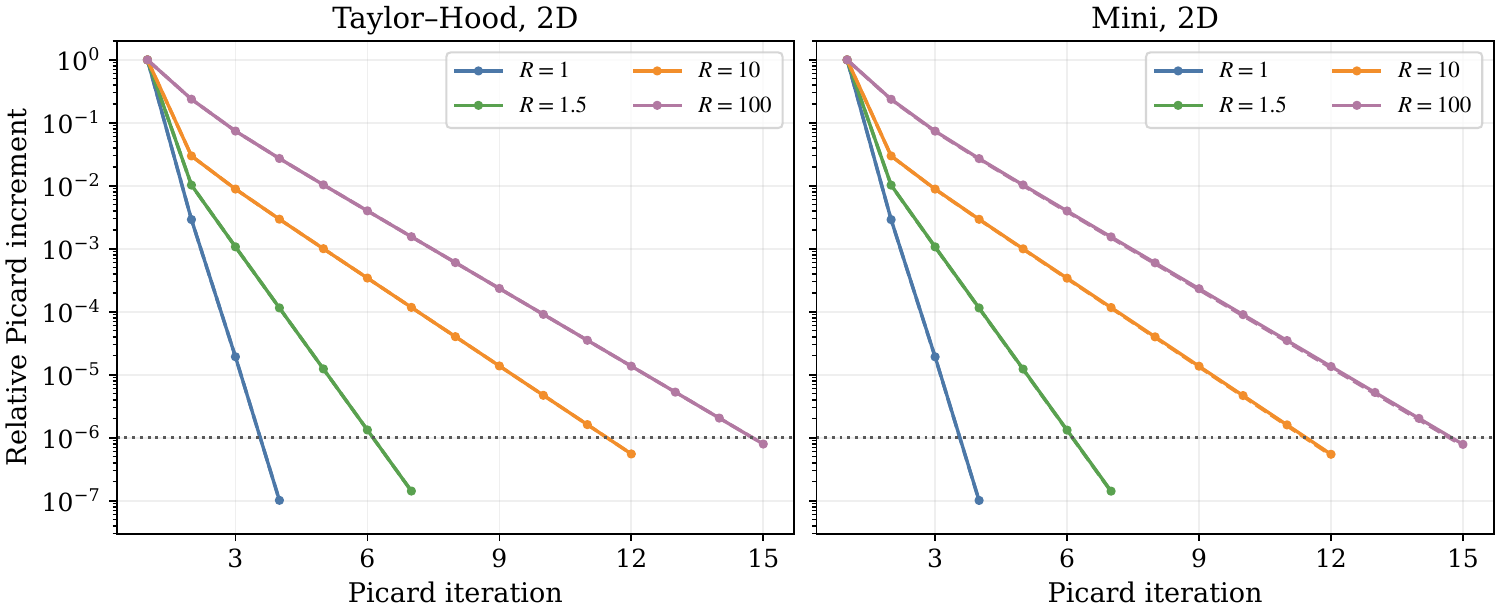}
\end{minipage}\par\smallskip
\begin{minipage}{\linewidth}
\centering
\includegraphics[width=0.8\linewidth]{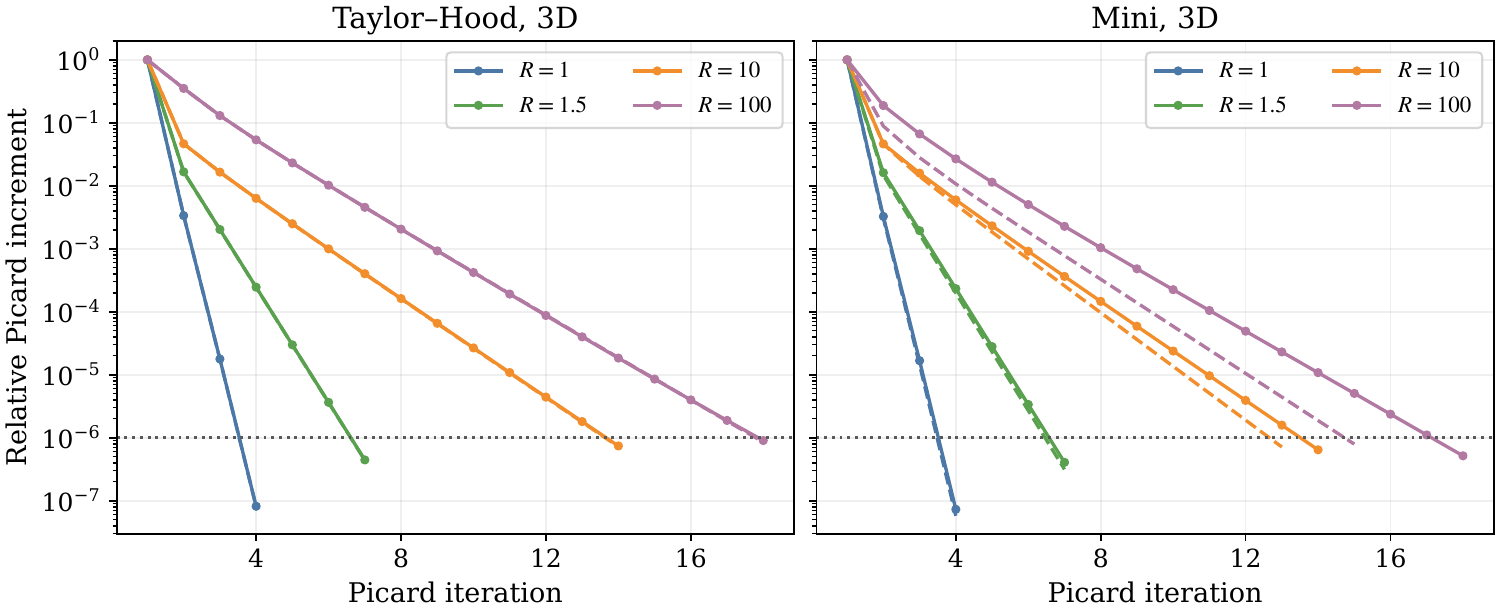}
\end{minipage}
\caption{Experiment \ref{subsec:rheological}. Relative Picard increments in 2D (top) and 3D (bottom), for Taylor--Hood (left) and mini (right). Dashed and solid curves refer to the coarser and finer meshes, with circles on the latter. The dotted line marks $10^{-6}$. The two mesh curves nearly coincide in 2D.}
\label{fig:robustness_picard_histories}
\end{figure}

Figure~\ref{fig:robustness_picard_histories} shows the Picard histories on both meshes. The curves are close under refinement, whereas increasing $R$ slows the decay of the increments. To quantify this behavior, we take the median of the last four available ratios $\|z_h^{k+1}-z_h^k\|_\star/\|z_h^k-z_h^{k-1}\|_\star$. On the finer 2D Taylor--Hood mesh, the resulting values are $0.107$, $0.342$, and $0.387$ for $R=1.5$, $10$, and $100$, respectively. The sufficient constants in Theorems~\ref{th:algorithm_3D} and~\ref{th:algorithm_2D} have not been evaluated. In particular, the four-step histories at $R=1$ are too short to identify an asymptotic factor.

The pressure results favor Taylor--Hood for the larger contrasts considered here, while the mini element retains first-order convergence for velocity and the scalars.

\subsubsection{Small viscosity and diffusivities}\label{subsec:viscosity}
Now, we aim to study the effects on the scheme when we choose small values for the viscosity and diffusivities. To do this, we fix $R=1.5$ and multiply the complete coefficient laws by positive factors $\nu_a(T,C)=a\,\nu(T,C),\;\kappa_b(T,C)=b\,\kappa(T,C),\;\mathcal D_c(T,C)=c\,\mathcal D(T,C).$ For $\eta\in\{1,10^{-1},10^{-2},10^{-3}\}$, we consider $(a,b,c)\in
\{(\eta,1,1),\ (1,\eta,1),\ (1,1,\eta),\ (\eta,\eta,\eta)\}.$ We fix a maximum of 100 Picard iterations. The case $\eta=1$ is shared by the four sweeps, giving 104 distinct discrete problems. Reducing only $\kappa$ or only $\mathcal D$ requires seven iterations in every tested case. Table~\ref{tab:robustness_small_coefficients} gives the counts for viscosity reduction and simultaneous coefficient reduction.

\begin{table}[hbt!]
\centering\footnotesize
\caption{Experiment \ref{subsec:viscosity}. Picard iteration counts for reduced coefficients. Each pair refers to the coarser and finer meshes, respectively. NC indicates that the stopping criterion is not met within 100 steps. An asterisk marks a converged run with at least one relative error above $100\%$.}
\label{tab:robustness_small_coefficients}
\begin{tabular}{|r|cc|cc|}
\hline\hline
 & \multicolumn{2}{c|}{2D} & \multicolumn{2}{c|}{3D} \\
$\eta$ & Taylor--Hood & Mini & Taylor--Hood & Mini \\
\hline
$h$ & \multicolumn{2}{c|}{0.04419 / 0.02210} & \multicolumn{2}{c|}{0.21651 / 0.10825} \\
$\mathrm{dof}$ & 17990 / 70790 & 9542 / 37510 & 25295 / 184599 & 13591 / 103207 \\
\hline
\multicolumn{5}{|c|}{$(a,b,c)=(\eta,1,1)$} \\
\hline
$1$ & $7 / 7$ & $7 / 7$ & $7 / 7$ & $7 / 7$ \\
$0.1$ & $7 / 7$ & $7 / 7$ & $7 / 7$ & $7 / 7$ \\
$0.01$ & $7 / 7$ & $8 / 7$ & $8^{*} / 7$ & $8^{*} / 7$ \\
$0.001$ & $8 / 8$ & $8^{*} / 8$ & $6^{*} / 8^{*}$ & $6^{*} / 7^{*}$ \\
\hline
\multicolumn{5}{|c|}{$(a,b,c)=(\eta,\eta,\eta)$} \\
\hline
$1$ & $7 / 7$ & $7 / 7$ & $7 / 7$ & $7 / 7$ \\
$0.1$ & $7 / 7$ & $7 / 7$ & $7 / 7$ & $7 / 7$ \\
$0.01$ & $9 / 9$ & $9 / 9$ & $9^{*} / 9$ & $8^{*} / 9$ \\
$0.001$ & $\mathrm{NC} / \mathrm{NC}$ & $\mathrm{NC} / \mathrm{NC}$ & $10^{*} / \mathrm{NC}$ & $46^{*} / \mathrm{NC}$ \\
\hline
\hline
\end{tabular}
\end{table}

Viscosity reduction alone requires six to eight iterations. Under simultaneous coefficient reduction, the counts increase to eight or nine at $\eta=10^{-2}$. At $\eta=10^{-3}$, under simultaneous coefficient reduction, all four 2D runs and both finer 3D runs reach the iteration limit. On the coarser 3D mesh, Taylor--Hood and mini meet the stopping criterion after 10 and 46 iterations, respectively, but their relative velocity errors are $3164.74\%$ and $4332.17\%$. Overall, 98 of the 104 problems meet the stopping criterion, and 11 of these have at least one relative error above $100\%$.

To examine the large velocity errors under viscosity reduction, we repeat the tests at $\eta=10^{-3}$ with $\mathsf p_\eta=\eta\mathsf p$, retaining the exact velocity and scalars and recomputing the source terms. Table~\ref{tab:pressure_control} shows that the velocity error decreases by factors between $20.4$ and $611.4$ when only viscosity is reduced. All eight control runs converge in nine steps. This sensitivity to the pressure amplitude is consistent with pressure-gradient contamination of the discrete velocity at small viscosity. We report absolute pressure errors because $\|\mathsf p_\eta\|_{L^2(\Omega)}$ is also reduced by a factor of $1000$.

\begin{table}[hbt!]
\centering\footnotesize
\caption{Experiment \ref{subsec:viscosity}. Pressure control for $(a,b,c)=(10^{-3},1,1)$. Superscripts $o$ and $s$ denote the original and scaled pressures, respectively. Iteration counts are given as original/scaled.}
\label{tab:pressure_control}
\begin{tabular}{|lrr|rr|cc|c|}
\hline\hline
Element & $\mathrm{dof}$ & $h$ & $100\widehat{E}_{\mathbf u}^{o}$ & $100\widehat{E}_{\mathbf u}^{s}$ & $E_{\mathsf p}^{o}$ & $E_{\mathsf p}^{s}$ & $\mathrm{it}$ \\
\hline
\multicolumn{8}{|c|}{2D} \\
\hline
Taylor--Hood & 17990 & 0.04419 & 32.446 & 0.288 & 4.025e-03 & 4.025e-06 & 8/9 \\
Taylor--Hood & 70790 & 0.02210 & 4.070 & 0.072 & 1.005e-03 & 1.005e-06 & 8/9 \\
Mini & 9542 & 0.04419 & 359.544 & 7.795 & 4.032e-03 & 4.872e-05 & 8/9 \\
Mini & 37510 & 0.02210 & 77.607 & 3.812 & 1.005e-03 & 1.274e-05 & 8/9 \\
\hline
\multicolumn{8}{|c|}{3D} \\
\hline
Taylor--Hood & 25295 & 0.21651 & 3119.262 & 5.102 & 7.351e-02 & 7.472e-05 & 6/9 \\
Taylor--Hood & 184599 & 0.10825 & 387.036 & 1.213 & 1.781e-02 & 1.783e-05 & 8/9 \\
Mini & 13591 & 0.21651 & 4320.666 & 49.831 & 7.601e-02 & 9.202e-04 & 6/9 \\
Mini & 103207 & 0.10825 & 1019.941 & 19.419 & 1.792e-02 & 2.379e-04 & 7/9 \\
\hline
\hline
\end{tabular}
\end{table}

The same pressure scaling does not restore convergence under simultaneous coefficient reduction: all eight control runs reach the 100-step limit. Thus, the viscosity-only control does not explain the failures when all three coefficients are reduced. These computations also do not evaluate the sufficient smallness conditions of the fixed-point theorems.

The joint stopping criterion can mask larger updates in individual fields. Of the 98 converged original runs, 60 have at least one separately normalized field increment above $10^{-6}$. For the coarser 3D mini run with simultaneous reduction at $\eta=10^{-3}$, the joint criterion is first met at step 46, but the concentration update is $1.09\times10^{-5}$. Continuing until all four separately normalized increments are below $10^{-6}$ requires 58 steps. The count of 46 in Table~\ref{tab:robustness_small_coefficients} corresponds to the original joint criterion.

Simultaneous coefficient reduction is the most demanding case for the undamped Picard iteration in this study. The pressure control and the continued iteration show why the stopping criterion, the accuracy of each field, and the coefficient regime must be assessed separately.

\subsection{Applications}
\label{sec:numerics_applications}
The aim of this experiment is to assess two configurations for combined cooling and solute removal. In 2D, we compare uniform and staged cooling to examine their effect on the temperature distribution, dissolved \(\mathrm{CO}_2\) removal, and pressure drop. In 3D, we study sequential cooling and dissolved \(\mathrm{NH}_3\) removal in two connected tanks, focusing on the thermal uniformity in the membrane tank and the fraction of solute removed.
\subsubsection{Staged cooling in a hollow-fiber oxygenator}
\label{sec:numerics_staged_cooling}
In this test, we consider a two-dimensional channel of length $1200\,\mu\mathrm{m}$ and height $H=200\,\mu\mathrm{m}$, with six square obstacles of side $60\,\mu\mathrm{m}$ representing the fiber cross sections. The centers of $S_1,\ldots,S_6$ are spaced $120\,\mu\mathrm{m}$ apart in the streamwise direction, starting at $x=200\,\mu\mathrm{m}$, and alternate between $y=70$ and $130\,\mu\mathrm{m}$. We use the equations in \eqref{eq:model} with $\mathbf{f}=\mathbf{0}$, $\mathsf{g}_1=\mathsf{g}_2=0$, and mixed boundary conditions. At the inlet, we prescribe a parabolic velocity profile with mean $U=1\,\mathrm{mm}/\mathrm{s}$, a temperature of $37\,^{\circ}\mathrm{C}$, and a normalized dissolved $\mathrm{CO}_2$ concentration $C=1$. We impose no slip on the walls and fibers and zero natural traction at the outlet for the full-gradient operator. Fibers $S_1,S_3,S_5$ are perfect sinks ($C=0$), and $S_2,S_4,S_6$ are held at $36$, $34.5$, and $33\,^{\circ}\mathrm{C}$, respectively. For each scalar, the diffusive flux vanishes on the remaining boundary.

We scale lengths by $H$, velocities by $U$, and pressure by $\mu_{37}U/H$, where $\mu_{37}=0.056\,\mathrm{Pa}\,\mathrm{s}$. We write $T=(\vartheta-33)/4$ for the dimensionless temperature, where $\vartheta$ is measured in degrees Celsius, and set $\widehat T=\min\{1,\max\{0,T\}\}$. With $a_1=4\log(1.02)$, $b_1=4/306.15$, and $\varepsilon_\infty=0.00345/0.056$, we use
\begin{align*}
\nu(T,C)&=\exp(a_1(1-\widehat T)),&
\varepsilon(x,s)&=\varepsilon_\infty+(1-\varepsilon_\infty)(1+(16.565s)^2)^{(0.357-1)/2},\\
\kappa(T,C)&=\frac{1.3\times10^{-7}}{UH},&
\mathcal{D}(T,C)&=\frac{1.5\times10^{-9}}{UH}\frac{1+b_1\widehat T}{1+b_1}\exp(a_1(\widehat T-1)),
\end{align*}
with $U$ and $H$ in SI units. We use Taylor--Hood elements and continuous piecewise quadratics for both scalars on a mesh graded near the fiber corners. We start Picard with the parabolic velocity, zero pressure, and $T=1$. The largest of the relative $L^2(\Omega)$ updates in $\nabla\mathbf{u}$, $\mathsf{p}$, and $\nabla T$, with each component normalized separately in the same way as in \eqref{eq:numerics_picard_stop}, falls below $10^{-6}$ in 12 steps for both cooling configurations. The scalar transport is assembled in conservative weak form, $-\int_\Omega \varphi_h\mathbf{u}_h\cdot\nabla\psi_h\,\mathrm{d}x+\int_{\Gamma_{\rm in}\cup\Gamma_{\rm out}}(\mathbf{u}_h\cdot\mathbf n)\varphi_h\psi_h\,\mathrm{d}s$, including the open-boundary fluxes. Here, $\Gamma_{\rm in}$ and $\Gamma_{\rm out}$ are the inlet and outlet, and $\mathbf{n}$ is the outward unit normal. The outlet condition fixes the pressure reference, and we solve for $C$ after convergence of the velocity and temperature.

\begin{figure}[hbt!]
\centering
\begin{minipage}{0.97\linewidth}
\centering
\includegraphics[width=0.95\linewidth]{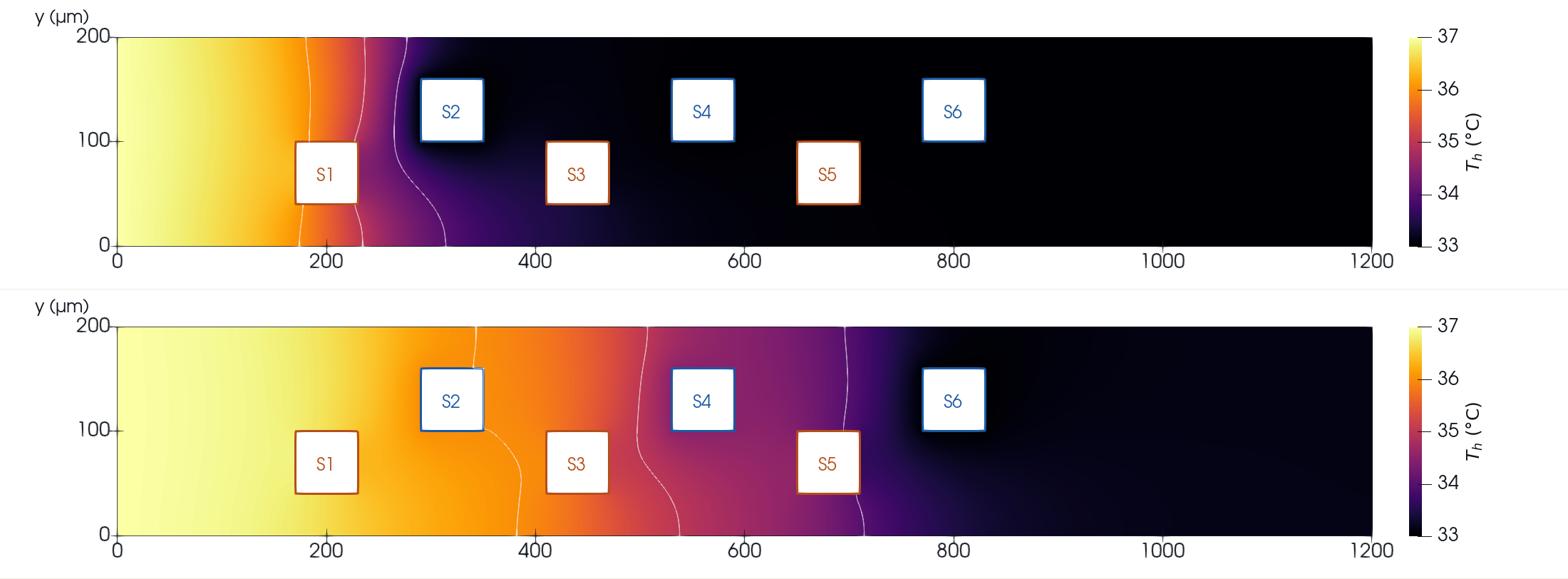}
\end{minipage}\par\nointerlineskip
\begin{minipage}{0.97\linewidth}
\centering
\includegraphics[width=0.85\linewidth]{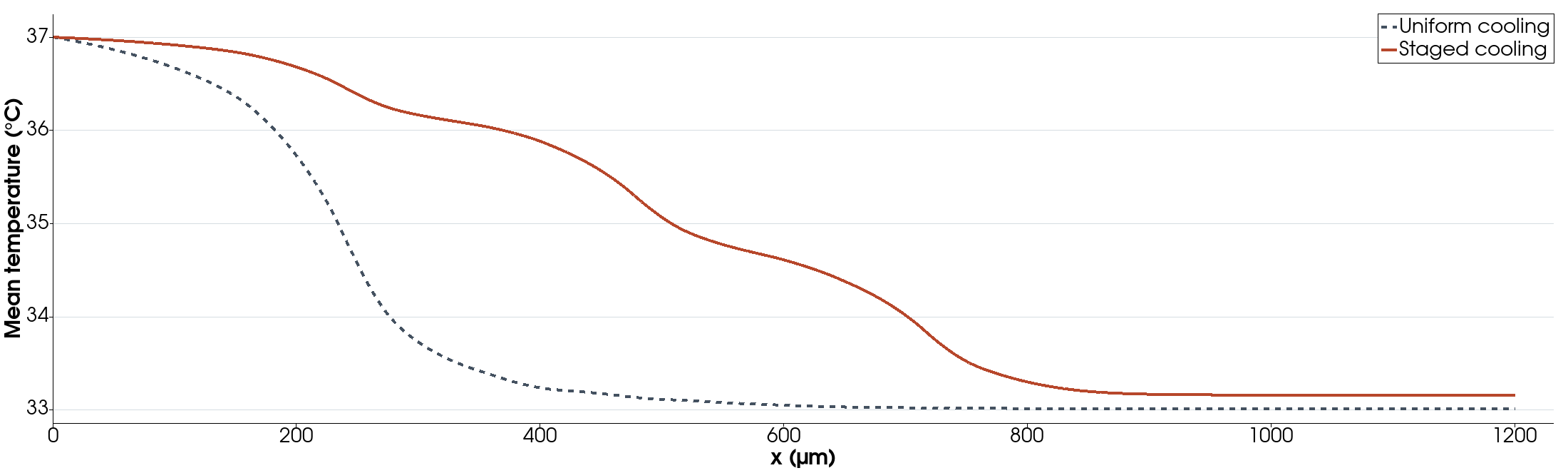}
\end{minipage}
\caption{Experiment \ref{sec:numerics_staged_cooling}. Temperature for uniform cooling (top) and staged cooling (middle), and flow-weighted mean temperature along the channel (bottom).}
\label{fig:experiment3_staged_cooling}
\end{figure}

\begin{figure}[t!]
\centering
\includegraphics[width=0.92\linewidth]{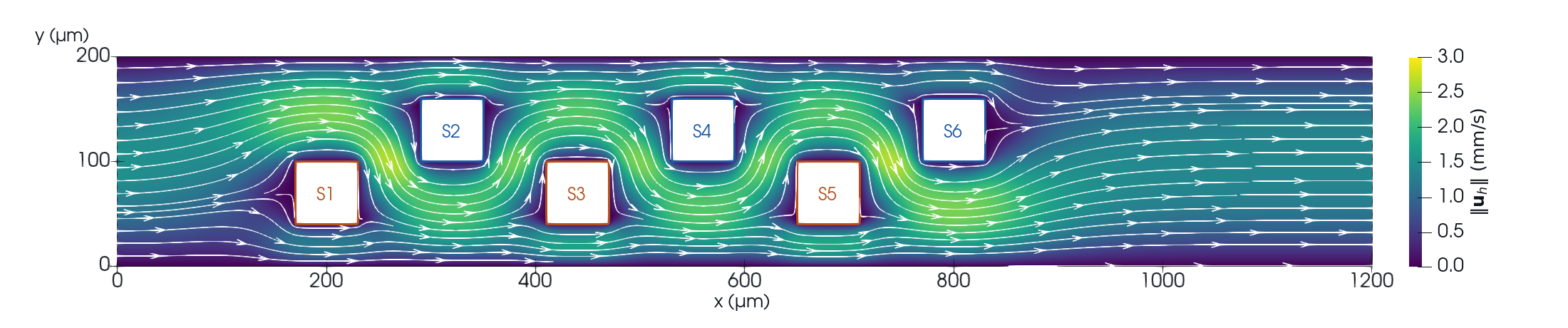}\par
\includegraphics[width=0.92\linewidth]{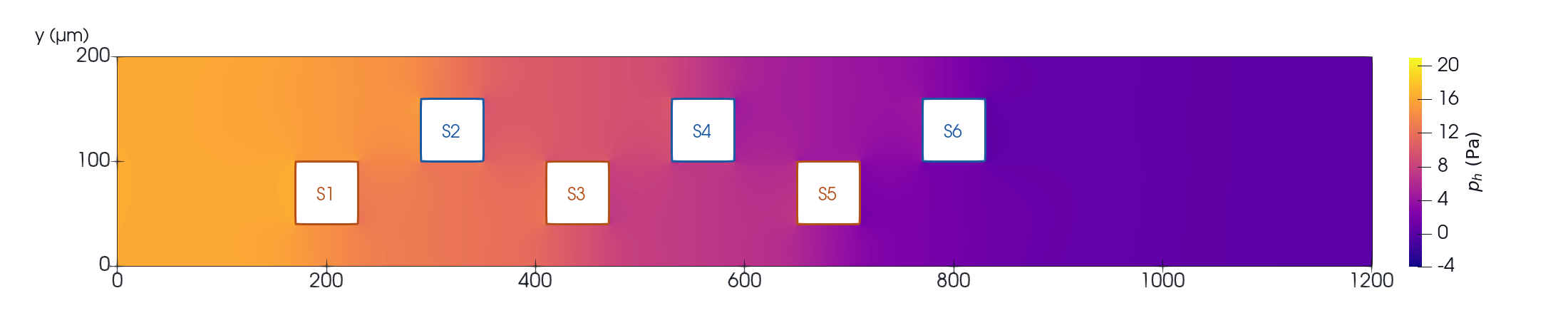}\par
\includegraphics[width=0.92\linewidth]{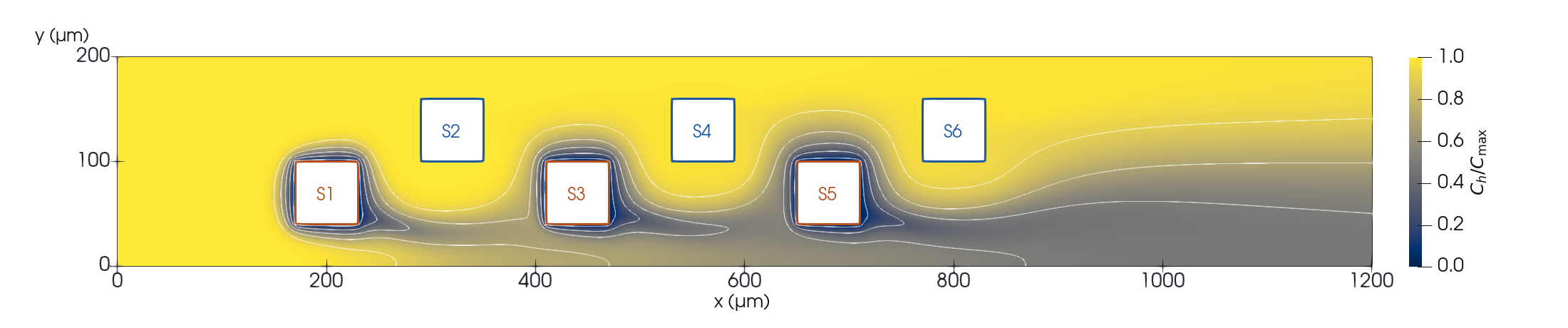}
\caption{Experiment \ref{sec:numerics_staged_cooling}. Speed and streamlines (top), pressure relative to the outlet mean (middle), and normalized dissolved $\mathrm{CO}_2$ concentration (bottom) for staged cooling.}
\label{fig:experiment3_velocity_concentration}
\end{figure}

Figure~\ref{fig:experiment3_staged_cooling} compares staged cooling with uniform cooling at $33\,^{\circ}\mathrm{C}$ on all three fibers. At $x=400\,\mu\mathrm{m}$, the mean temperature, weighted by the streamwise velocity, has undergone $94\%$ of its inlet-to-outlet decrease under uniform cooling and $29\%$ under staged cooling. Staging thus spreads the temperature drop along the channel, while the outlet means ($33.005$ and $33.154\,^{\circ}\mathrm{C}$) differ by only $0.15\,^{\circ}\mathrm{C}$.\par
In Figure~\ref{fig:experiment3_velocity_concentration}, the streamlines bend alternately above and below the fibers, with a peak speed of $2.65\,\mathrm{mm}/\mathrm{s}$ in the gaps. The pressure falls mainly across the fiber array and varies little downstream. The difference between the inlet and outlet pressure means is $16.20\,\mathrm{Pa}$ for staged cooling and $16.68\,\mathrm{Pa}$ for uniform cooling.

The concentration decreases near the sinks $S_1,S_3,S_5$, with depleted regions extending downstream along the lower half of the channel. The sinks remove $28.78\%$ of the incoming dissolved $\mathrm{CO}_2$ flux under staged cooling, compared with $28.18\%$ under uniform cooling. These fractions are computed from the discrete boundary reactions, with the diffusive contribution included in the inlet flux.

For this configuration, prescribing different fiber temperatures provides control over where cooling occurs along the channel, with modest changes in dissolved $\mathrm{CO}_2$ removal and the pressure required to maintain the prescribed flow rate.

\subsubsection{Cooling and solute removal in two connected tanks}
\label{sec:numerics_tanks_3d}
Let us consider now two cubic tanks, $A$ and $B$, of side $10\,\mathrm{cm}$, joined by a cylindrical pipe of length $5\,\mathrm{cm}$ and radius $1\,\mathrm{cm}$. The pipe axis lies $2\,\mathrm{cm}$ above the tank floor. The inlet and outlet pipes, located at height $8\,\mathrm{cm}$, have radius $0.8\,\mathrm{cm}$ and length $3\,\mathrm{cm}$. Tank $A$ contains three horizontal cooling rings of major radius $3\,\mathrm{cm}$ and tube radius $0.4\,\mathrm{cm}$, centered at heights $3$, $5$, and $7\,\mathrm{cm}$. Tank $B$ contains five cylindrical membranes of radius $0.4\,\mathrm{cm}$, parallel to the $y$-axis and spanning the tank.

We prescribe a parabolic inlet velocity with mean $U=2\,\mathrm{mm}/\mathrm{s}$, a temperature of $55\,^{\circ}\mathrm{C}$, and a normalized dissolved $\mathrm{NH}_3$ concentration $C=1$. The resulting flow rate is $24.02\,\mathrm{mL}/\mathrm{min}$. The rings are held at $35\,^{\circ}\mathrm{C}$, and the membranes are perfect sinks with $C=0$. We impose no slip on the walls, rings, and membranes, zero natural traction for the full-gradient operator at the outlet, and zero diffusive flux for each scalar on the remaining boundary.

We scale lengths by $H=0.1\,\mathrm{m}$, velocities by $U$, and pressure by $\mu_\star U/H$, with $\mu_\star=0.3\,\mathrm{Pa}\,\mathrm{s}$. We set $T=(\vartheta-35)/20$, with $\vartheta$ measured in degrees Celsius, and use $\widehat T=\min\{1,\max\{0,T\}\}$. The dimensionless coefficients are
\begin{align*}
\nu(T,C)&=\exp(a_2(1-\widehat T)),&
\varepsilon(x,s)&=0.02+0.98(1+(50s)^2)^{-1/4},\\
\kappa(T,C)&=\frac{1}{400},&
\mathcal{D}(T,C)&=\frac{1-b_2(1-\widehat T)}{2000}\exp(-0.41(1-\widehat T)),
\end{align*}
where $a_2=17000\cdot20/(8.314462618\cdot308.15^2)\approx0.430646$ and $b_2=20/328.15$. We take $\mathbf{f}(T,C)=0.5T\mathbf{e}_3$, where $\mathbf{e}_3$ is the vertical unit vector, and $\mathsf{g}_1=\mathsf{g}_2=0$.

We use Taylor--Hood elements and quadratic scalar spaces on a mesh of $382\,328$ tetrahedra, with $2\,850\,919$ total degrees of freedom. The scalar equations use the advective term $\int_\Omega(\mathbf{u}_h\cdot\nabla\varphi_h)\psi_h\,\mathrm{d}x$. We start with zero interior velocity and pressure and $T=1$ in the interior, imposing the prescribed boundary values. With the full buoyancy force from the first step and no relaxation, Picard meets the blockwise tolerance $10^{-6}$ used in 2D after 13 steps. The resulting fields differ by less than $4.2\times10^{-7}$ in the respective relative norms from those obtained in five steps with precomputed initial fields.

\begin{figure}[hbt!]
\centering
\begin{minipage}[t]{0.45\linewidth}
\centering
\includegraphics[width=\linewidth,trim=100 110 40 205,clip]{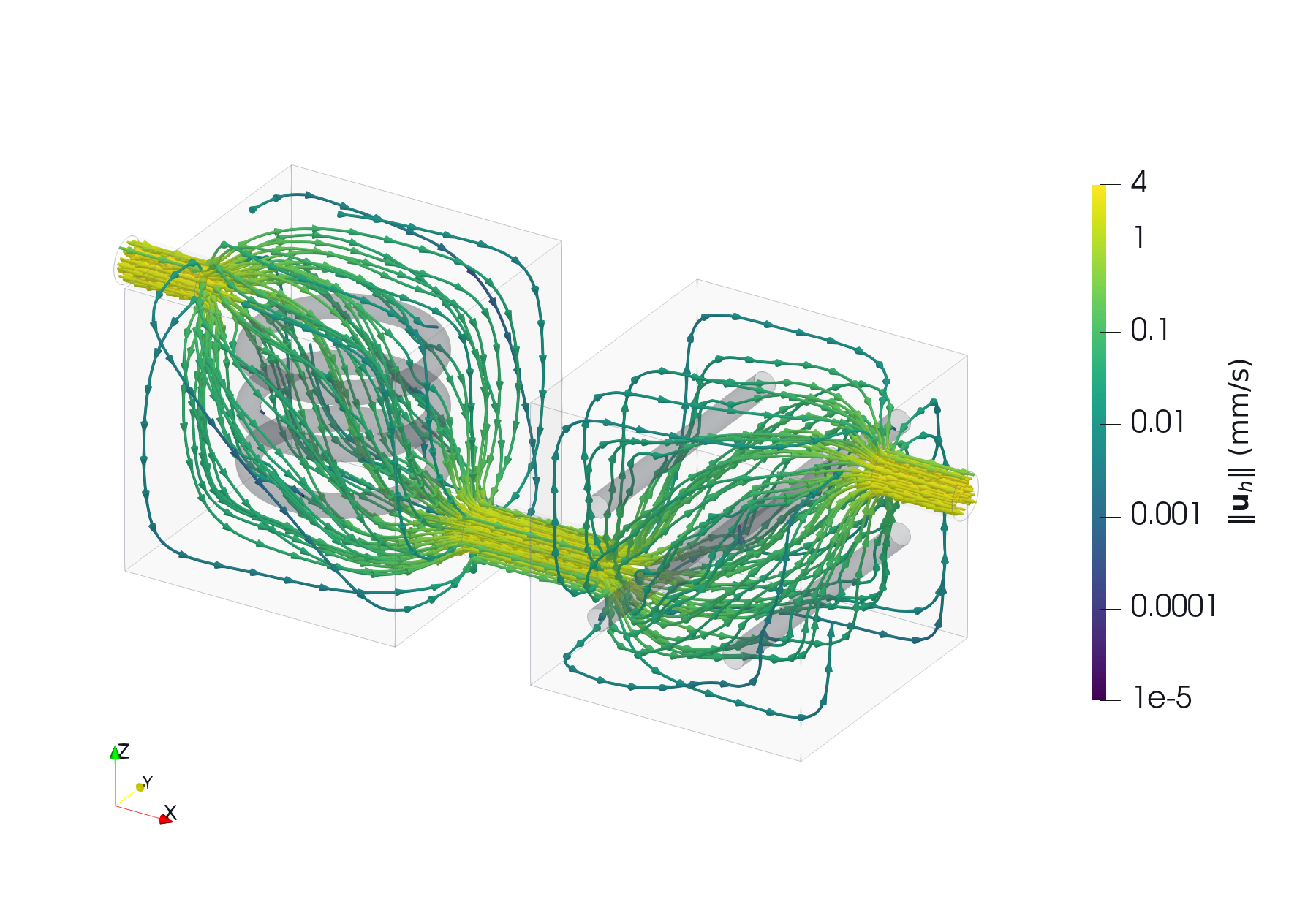}
\end{minipage}\hfill
\begin{minipage}[t]{0.45\linewidth}
\centering
\includegraphics[width=\linewidth,trim=100 110 40 205,clip]{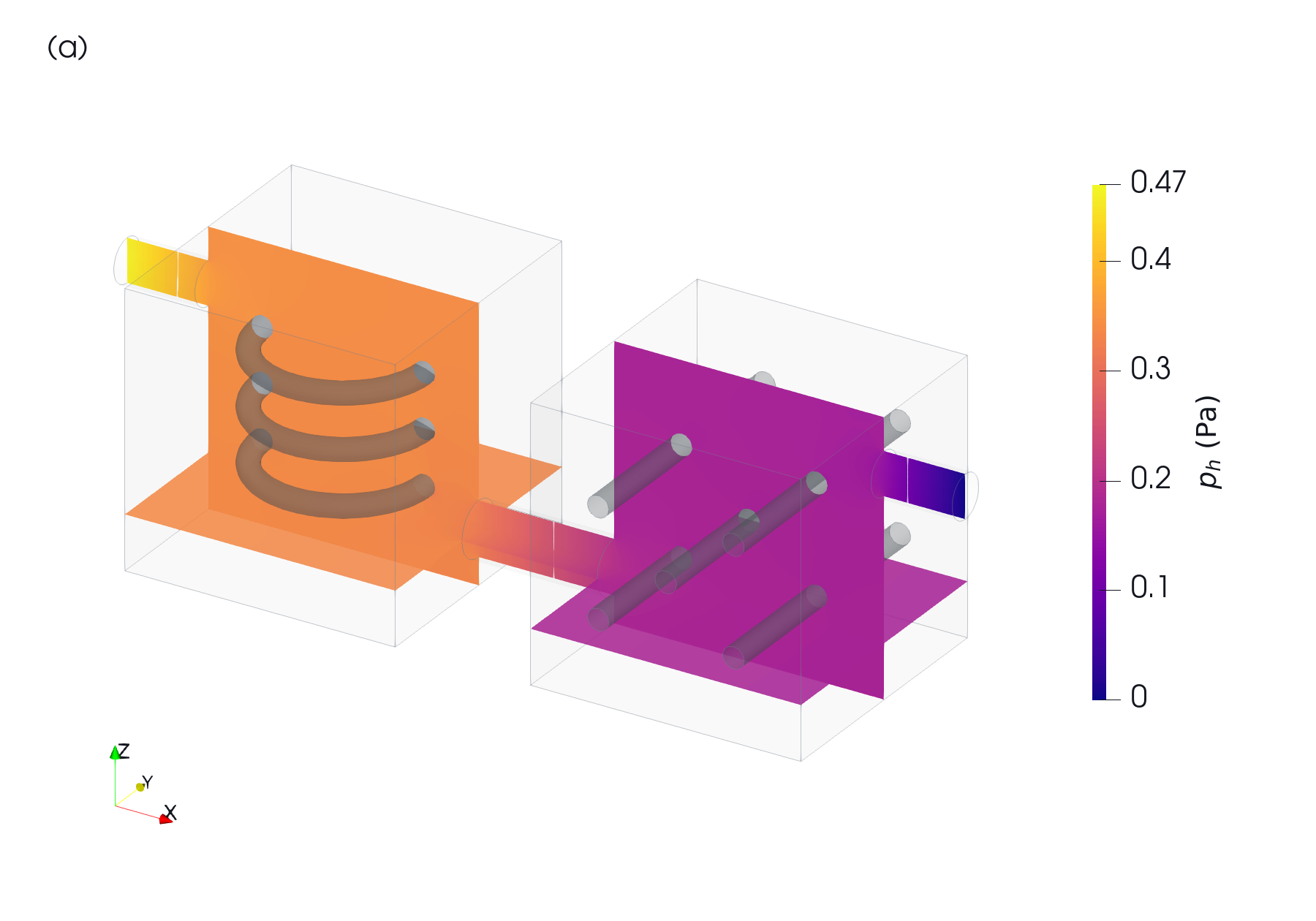}
\end{minipage}\par\smallskip
\begin{minipage}[t]{0.45\linewidth}
\centering
\includegraphics[width=\linewidth,trim=100 110 40 205,clip]{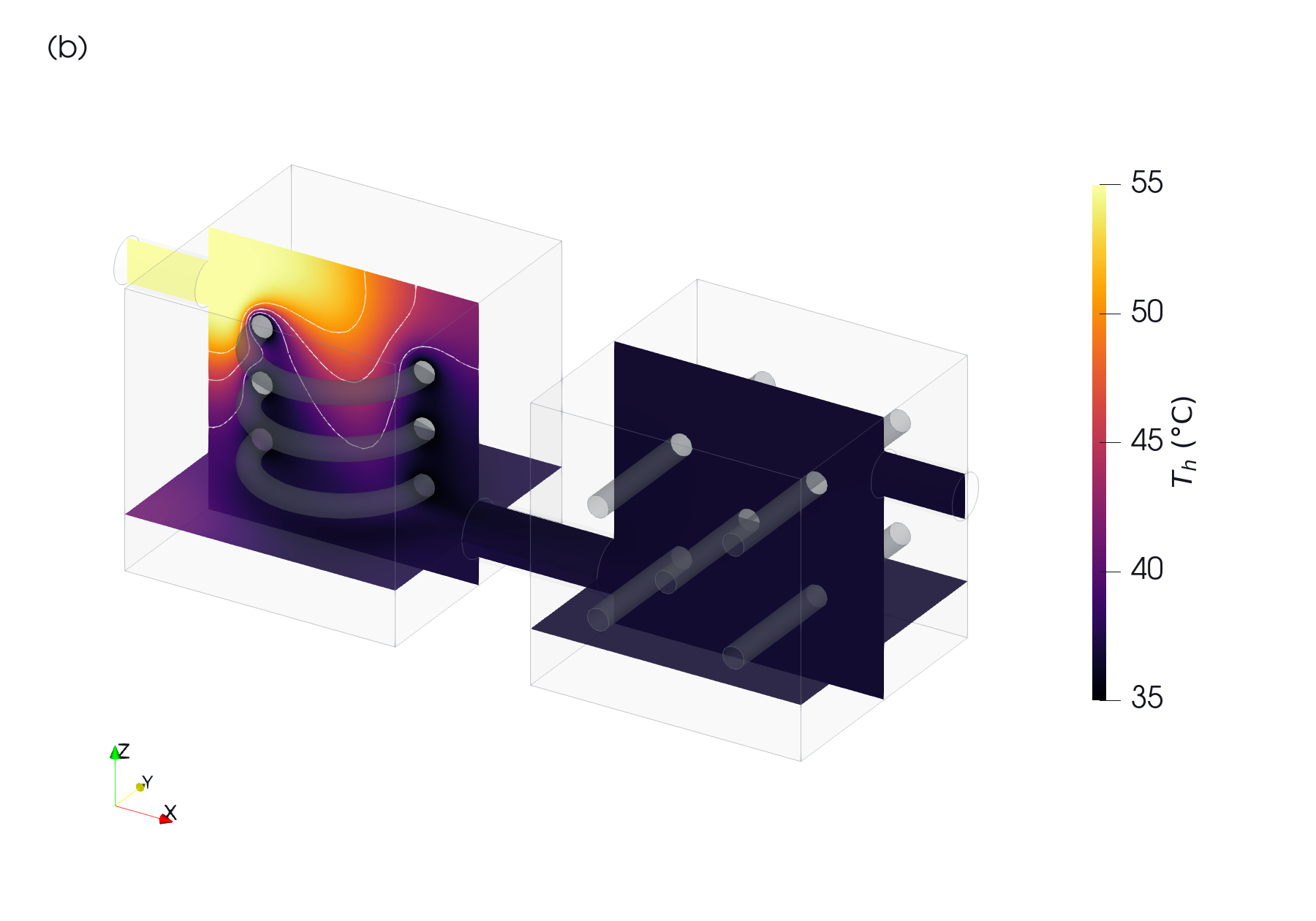}
\end{minipage}\hfill
\begin{minipage}[t]{0.45\linewidth}
\centering
\includegraphics[width=\linewidth,trim=100 110 40 205,clip]{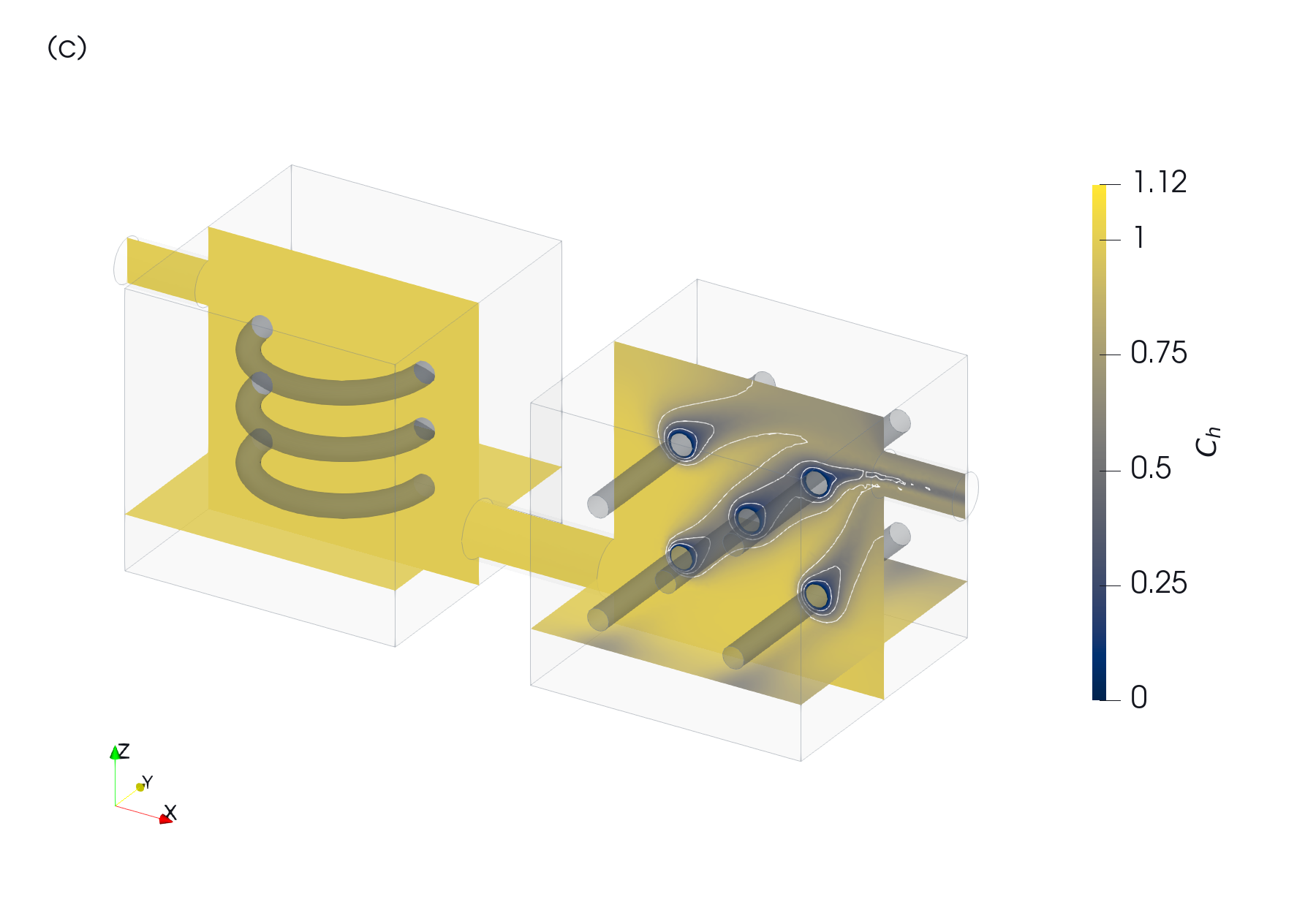}
\end{minipage}
\caption{Experiment \ref{sec:numerics_tanks_3d}. Streamlines colored by speed (top left), pressure (top right), temperature (bottom left), and concentration (bottom right). The cuts are at $y=5\,\mathrm{cm}$ and $z=2\,\mathrm{cm}$. The speed scale is logarithmic.}
\label{fig:experiment3_tanks_3d}
\end{figure}

Figure~\ref{fig:experiment3_tanks_3d} shows the faster flow in the narrow pipes and the slower motion around the internal surfaces of both tanks. The maximum nodal speed is $3.98\,\mathrm{mm}/\mathrm{s}$. The pressure is nearly uniform within each tank and falls mainly along the pipes. Its mean inlet-to-outlet drop is $0.459\,\mathrm{Pa}$.

Most of the cooling occurs in $A$, where the inlet stream meets the cold rings. The flow-weighted outlet temperature is $36.624\,^{\circ}\mathrm{C}$. Tank $B$ is almost isothermal, with a volume mean of $36.628\,^{\circ}\mathrm{C}$ and a standard deviation of $0.042\,^{\circ}\mathrm{C}$. The computed temperature ranges from $35$ to $55\,^{\circ}\mathrm{C}$ within numerical precision.

The concentration remains close to the inlet value in $A$ and decreases around the membranes in $B$. The fraction of incoming solute flux removed, computed from the discrete boundary reactions as in the two-dimensional case, is $41.08\%$. The quadratic concentration reaches a maximum of $1.029$, a local numerical overshoot of $2.91\%$ above the inlet value.

To examine this overshoot, we refine the concentration mesh locally near the membranes, keeping the quadratic velocity and temperature fields fixed by exact transfer to the nested mesh. Table~\ref{tab:concentration_refinement} shows that the overshoot decreases to $0.73\%$, while the extraction fraction changes by $0.0364$ percentage points. The relative $L^2(\Omega)$ change in concentration is $4.81\times10^{-4}$. The maxima are evaluated over the quadratic polynomials on each tetrahedron, including their faces and edges.

\begin{table}[hbt!]
\centering\footnotesize
\caption{Experiment \ref{sec:numerics_tanks_3d}. Local refinement of the passive concentration problem with fixed velocity and temperature. The maximum mesh diameter remains unchanged because refinement is local. Here, $\mathrm{dof}_C$ is the number of concentration degrees of freedom.}
\label{tab:concentration_refinement}
\begin{tabular}{|lrrr|rr|}
\hline\hline
Mesh & Tetrahedra & $\mathrm{dof}_C$ & $h$ & $\max C_h$ & Extraction (\%) \\
\hline
Original & 382328 & 555124 & 0.12372 & 1.029078 & 41.0847 \\
Refined & 751629 & 1074687 & 0.12372 & 1.007310 & 41.0483 \\
\hline\hline
\end{tabular}
\end{table}

For these data, the tanks act as successive cooling and solute-removal stages. The membranes in $B$ receive fluid at an almost uniform temperature after cooling in $A$.

\section*{Conflict of interest}
The authors report no conflicts of interest.

\section*{Data availability}
The datasets generated and/or analyzed during this study are available from the corresponding author upon reasonable and justified request.

\bibliographystyle{siamplain}
\bibliography{biblio}

@article {MR4610451,
    AUTHOR = {Grasselli, Maurizio and Parolini, Nicola and Poiatti, Andrea and Verani, Marco},
     TITLE = {Non-isothermal non-{N}ewtonian fluids: the stationary case},
   JOURNAL = {Math. Models Methods Appl. Sci.},
  FJOURNAL = {Mathematical Models and Methods in Applied Sciences},
    VOLUME = {33},
      YEAR = {2023},
    NUMBER = {9},
     PAGES = {1747--1801},
      ISSN = {0218-2025,1793-6314},
   MRCLASS = {76A05 (35Q30 76M10)},
  MRNUMBER = {4610451},
       DOI = {10.1142/S0218202523500410},
       URL = {https://doi.org/10.1142/S0218202523500410},
}

@article {MR3262601,
    AUTHOR = {Bul\'{\i}\v{c}ek, Miroslav and Pust\v{e}jovsk\'{a}, Petra},
     TITLE = {Existence analysis for a model describing flow of an incompressible chemically reacting non-{N}ewtonian fluid},
   JOURNAL = {SIAM J. Math. Anal.},
  FJOURNAL = {SIAM Journal on Mathematical Analysis},
    VOLUME = {46},
      YEAR = {2014},
    NUMBER = {5},
     PAGES = {3223--3240},
      ISSN = {0036-1410,1095-7154},
   MRCLASS = {35Q35 (35A01 76A05 80A32)},
  MRNUMBER = {3262601},
       DOI = {10.1137/130927589},
       URL = {https://doi.org/10.1137/130927589},
}

@article {MR3023245,
    AUTHOR = {Bul\'{\i}\v{c}ek, Miroslav and Pust\v{e}jovsk\'{a}, Petra},
     TITLE = {On existence analysis of steady flows of generalized {N}ewtonian fluids with concentration dependent power-law index},
   JOURNAL = {J. Math. Anal. Appl.},
  FJOURNAL = {Journal of Mathematical Analysis and Applications},
    VOLUME = {402},
      YEAR = {2013},
    NUMBER = {1},
     PAGES = {157--166},
      ISSN = {0022-247X,1096-0813},
   MRCLASS = {76D05 (76Z05 92C35)},
  MRNUMBER = {3023245},
MRREVIEWER = {Jeffrey\ J.\ Heys},
       DOI = {10.1016/j.jmaa.2012.12.066},
       URL = {https://doi.org/10.1016/j.jmaa.2012.12.066},
}

@article {MR3834434,
    AUTHOR = {Ko, Seungchan and Pust\v{e}jovsk\'{a}, Petra and S\"{u}li, Endre},
     TITLE = {Finite element approximation of an incompressible chemically reacting non-{N}ewtonian fluid},
   JOURNAL = {ESAIM Math. Model. Numer. Anal.},
  FJOURNAL = {ESAIM. Mathematical Modelling and Numerical Analysis},
    VOLUME = {52},
      YEAR = {2018},
    NUMBER = {2},
     PAGES = {509--541},
      ISSN = {2822-7840,2804-7214},
   MRCLASS = {65N30 (76A05 76V05 80A32)},
  MRNUMBER = {3834434},
MRREVIEWER = {V\'{\i}t\ Dolej\v{s}\'{\i}},
       DOI = {10.1051/m2an/2017043},
       URL = {https://doi.org/10.1051/m2an/2017043},
}

@article {MR3904139,
    AUTHOR = {Ko, Seungchan and S\"{u}li, Endre},
     TITLE = {Finite element approximation of steady flows of generalized {N}ewtonian fluids with concentration-dependent power-law index},
   JOURNAL = {Math. Comp.},
  FJOURNAL = {Mathematics of Computation},
    VOLUME = {88},
      YEAR = {2019},
    NUMBER = {317},
     PAGES = {1061--1090},
      ISSN = {0025-5718,1088-6842},
   MRCLASS = {65N30 (76A05 76M10)},
  MRNUMBER = {3904139},
MRREVIEWER = {Kamel\ Al-Khaled},
       DOI = {10.1090/mcom/3379},
       URL = {https://doi.org/10.1090/mcom/3379},
}

@article {MR4344862,
    AUTHOR = {Caucao, Sergio and Gatica, Gabriel N. and Ortega, Juan P.},
     TITLE = {A fully-mixed formulation in {B}anach spaces for the coupling of the steady {B}rinkman-{F}orchheimer and double-diffusion equations},
   JOURNAL = {ESAIM Math. Model. Numer. Anal.},
  FJOURNAL = {ESAIM. Mathematical Modelling and Numerical Analysis},
    VOLUME = {55},
      YEAR = {2021},
    NUMBER = {6},
     PAGES = {2725--2758},
      ISSN = {2822-7840,2804-7214},
   MRCLASS = {65N30 (35Q79 65N12 65N15 76D07 76R05)},
  MRNUMBER = {4344862},
       DOI = {10.1051/m2an/2021072},
       URL = {https://doi.org/10.1051/m2an/2021072},
}

@article {MR4659442,
    AUTHOR = {Carrasco, Sergio and Caucao, Sergio and Gatica, Gabriel N.},
     TITLE = {New mixed finite element methods for the coupled convective {B}rinkman-{F}orchheimer and double-diffusion equations},
   JOURNAL = {J. Sci. Comput.},
  FJOURNAL = {Journal of Scientific Computing},
    VOLUME = {97},
      YEAR = {2023},
    NUMBER = {3},
     PAGES = {Paper No. 61, 49},
      ISSN = {0885-7474,1573-7691},
   MRCLASS = {65N30 (35Q79 65N12 65N15 80A19)},
  MRNUMBER = {4659442},
       DOI = {10.1007/s10915-023-02371-7},
       URL = {https://doi.org/10.1007/s10915-023-02371-7},
}

@article {MR4404048,
    AUTHOR = {\'{A}lvarez, Mario and Colmenares, Eligio and Sequeira, Fil\'{a}nder A.},
     TITLE = {Analysis of a semi-augmented mixed finite element method for double-diffusive natural convection in porous media},
   JOURNAL = {Comput. Math. Appl.},
  FJOURNAL = {Computers \& Mathematics with Applications. An International Journal},
    VOLUME = {114},
      YEAR = {2022},
     PAGES = {112--131},
      ISSN = {0898-1221,1873-7668},
   MRCLASS = {65N30 (65N15 76S05)},
  MRNUMBER = {4404048},
       DOI = {10.1016/j.camwa.2022.03.032},
       URL = {https://doi.org/10.1016/j.camwa.2022.03.032},
}

@article {MR1069652,
    AUTHOR = {Baranger, Jacques and Najib, Khalid},
     TITLE = {Analyse num\'{e}rique des \'{e}coulements quasi-newtoniens dont la viscosit\'{e} ob\'{e}it \`a la loi puissance ou la loi de {C}arreau},
   JOURNAL = {Numer. Math.},
  FJOURNAL = {Numerische Mathematik},
    VOLUME = {58},
      YEAR = {1990},
    NUMBER = {1},
     PAGES = {35--49},
      ISSN = {0029-599X,0945-3245},
   MRCLASS = {65N30 (76A05 76D05)},
  MRNUMBER = {1069652},
MRREVIEWER = {Teodor\ Potra},
       DOI = {10.1007/BF01385609},
       URL = {https://doi.org/10.1007/BF01385609},
}

@article {MR1211613,
    AUTHOR = {Sandri, D.},
     TITLE = {Sur l'approximation num\'{e}rique des \'{e}coulements quasi-newtoniens dont la viscosit\'{e} suit la loi puissance ou la loi de {C}arreau},
   JOURNAL = {RAIRO Mod\'{e}l. Math. Anal. Num\'{e}r.},
  FJOURNAL = {RAIRO Mod\'{e}lisation Math\'{e}matique et Analyse Num\'{e}rique},
    VOLUME = {27},
      YEAR = {1993},
    NUMBER = {2},
     PAGES = {131--155},
      ISSN = {0764-583X},
   MRCLASS = {76M10 (65N30 76A05)},
  MRNUMBER = {1211613},
MRREVIEWER = {C.\ Wayne\ Mastin},
       DOI = {10.1051/m2an/1993270201311},
       URL = {https://doi.org/10.1051/m2an/1993270201311},
}

@article {MR1213411,
    AUTHOR = {Barrett, John W. and Liu, W. B.},
     TITLE = {Finite element error analysis of a quasi-{N}ewtonian flow obeying the {C}arreau or power law},
   JOURNAL = {Numer. Math.},
  FJOURNAL = {Numerische Mathematik},
    VOLUME = {64},
      YEAR = {1993},
    NUMBER = {4},
     PAGES = {433--453},
      ISSN = {0029-599X,0945-3245},
   MRCLASS = {65N30 (65N15 76A05 76M10)},
  MRNUMBER = {1213411},
MRREVIEWER = {A.\ Jami},
       DOI = {10.1007/BF01388698},
       URL = {https://doi.org/10.1007/BF01388698},
}

@article {MR1301740,
    AUTHOR = {Barrett, John W. and Liu, W. B.},
     TITLE = {Quasi-norm error bounds for the finite element approximation of a non-{N}ewtonian flow},
   JOURNAL = {Numer. Math.},
  FJOURNAL = {Numerische Mathematik},
    VOLUME = {68},
      YEAR = {1994},
    NUMBER = {4},
     PAGES = {437--456},
      ISSN = {0029-599X,0945-3245},
   MRCLASS = {65N30 (65N15 76A05 76M10)},
  MRNUMBER = {1301740},
MRREVIEWER = {G.\ W.\ Hedstrom},
       DOI = {10.1007/s002110050071},
       URL = {https://doi.org/10.1007/s002110050071},
}

@article {MR1348583,
    AUTHOR = {Baranger, Jacques and Mikeli\'{c}, Andro},
     TITLE = {Stationary solutions to a quasi-{N}ewtonian flow with viscous heating},
   JOURNAL = {Math. Models Methods Appl. Sci.},
  FJOURNAL = {Mathematical Models and Methods in Applied Sciences},
    VOLUME = {5},
      YEAR = {1995},
    NUMBER = {6},
     PAGES = {725--738},
      ISSN = {0218-2025,1793-6314},
   MRCLASS = {76A05 (35Q35)},
  MRNUMBER = {1348583},
MRREVIEWER = {Niko\ Sauer},
       DOI = {10.1142/S0218202595000401},
       URL = {https://doi.org/10.1142/S0218202595000401},
}

@article {MR1636360,
    AUTHOR = {Wardi, S.},
     TITLE = {A convergence result for an iterative method for the equations of a stationary quasi-{N}ewtonian flow with temperature dependent viscosity},
   JOURNAL = {RAIRO Mod\'{e}l. Math. Anal. Num\'{e}r.},
  FJOURNAL = {RAIRO Mod\'{e}lisation Math\'{e}matique et Analyse Num\'{e}rique},
    VOLUME = {32},
      YEAR = {1998},
    NUMBER = {4},
     PAGES = {391--404},
      ISSN = {0764-583X},
   MRCLASS = {76A05 (65M99 76M99)},
  MRNUMBER = {1636360},
MRREVIEWER = {Bruce\ A.\ Finlayson},
       DOI = {10.1051/m2an/1998320403911},
       URL = {https://doi.org/10.1051/m2an/1998320403911},
}

@article {MR2320257,
    AUTHOR = {Zhu, Jiang and Loula, Abimael F. D. and Karam F., Jos\'{e} and Guerreiro, Jo\~{a}o N. C.},
     TITLE = {Finite element analysis of a coupled thermally dependent viscosity flow problem},
   JOURNAL = {Comput. Appl. Math.},
  FJOURNAL = {Computational \& Applied Mathematics},
    VOLUME = {26},
      YEAR = {2007},
    NUMBER = {1},
     PAGES = {45--66},
      ISSN = {1807-0302},
   MRCLASS = {65N30 (76M10)},
  MRNUMBER = {2320257},
MRREVIEWER = {Riccardo\ Sacco},
       DOI = {10.1590/S0101-82052007000100003},
       URL = {https://doi.org/10.1590/S0101-82052007000100003},
}

@article {MR2794183,
    AUTHOR = {Zhu, Jiang and Yu, Xijun},
     TITLE = {Thermally coupled quasi-{N}ewtonian flows: analysis and computation},
   JOURNAL = {J. Comput. Appl. Math.},
  FJOURNAL = {Journal of Computational and Applied Mathematics},
    VOLUME = {235},
      YEAR = {2011},
    NUMBER = {13},
     PAGES = {3898--3909},
      ISSN = {0377-0427,1879-1778},
   MRCLASS = {76D05 (35A01 35A02 35Q35 65N30 76A05)},
  MRNUMBER = {2794183},
MRREVIEWER = {Georgios\ C.\ Georgiou},
       DOI = {10.1016/j.cam.2011.01.036},
       URL = {https://doi.org/10.1016/j.cam.2011.01.036},
}

@article {MR3070180,
    AUTHOR = {Zhang, Jiansong and Zhu, Jiang and Yu, Xijun and Loula, Abimael F. D. and Bevilacqua, Luiz},
     TITLE = {Mixed finite element analysis of thermally coupled quasi-{N}ewtonian flows},
   JOURNAL = {Int. J. Numer. Anal. Model. Ser. B},
  FJOURNAL = {International Journal of Numerical Analysis and Modeling. Series B},
    VOLUME = {4},
      YEAR = {2013},
    NUMBER = {1},
     PAGES = {35--49},
      ISSN = {1923-2942,1923-2950},
   MRCLASS = {65N30 (65N15 76A05 76M10 80A20)},
  MRNUMBER = {3070180},
MRREVIEWER = {Abdallah\ Bradji},
}

@article {MR4379972,
    AUTHOR = {Farrell, Patrick and Gazca Orozco, Pablo Alexei and S\"{u}li, Endre},
     TITLE = {Finite element approximation and preconditioning for anisothermal flow of implicitly-constituted non-{N}ewtonian fluids},
   JOURNAL = {Math. Comp.},
  FJOURNAL = {Mathematics of Computation},
    VOLUME = {91},
      YEAR = {2022},
    NUMBER = {334},
     PAGES = {659--697},
      ISSN = {0025-5718,1088-6842},
   MRCLASS = {65N30 (65F08 65N22 76A05)},
  MRNUMBER = {4379972},
       DOI = {10.1090/mcom/3703},
       URL = {https://doi.org/10.1090/mcom/3703},
}

@article {MR2515781,
    AUTHOR = {Bul\'{\i}\v{c}ek, M. and M\'{a}lek, J. and Rajagopal, K. R.},
     TITLE = {Mathematical analysis of unsteady flows of fluids with pressure, shear-rate, and temperature dependent material moduli that slip at solid boundaries},
   JOURNAL = {SIAM J. Math. Anal.},
  FJOURNAL = {SIAM Journal on Mathematical Analysis},
    VOLUME = {41},
      YEAR = {2009},
    NUMBER = {2},
     PAGES = {665--707},
      ISSN = {0036-1410,1095-7154},
   MRCLASS = {76D03 (35Q35 35Q74 76A05)},
  MRNUMBER = {2515781},
MRREVIEWER = {Ilia\ S.\ Mogilevski\u{\i}},
       DOI = {10.1137/07069540X},
       URL = {https://doi.org/10.1137/07069540X},
}

@book {CiarletBook,
    AUTHOR = {Ciarlet, Philippe G.},
     TITLE = {The finite element method for elliptic problems},
    SERIES = {Classics in Applied Mathematics},
    VOLUME = {40},
      NOTE = {Reprint of the 1978 original},
 PUBLISHER = {Society for Industrial and Applied Mathematics (SIAM), Philadelphia, PA},
      YEAR = {2002},
     PAGES = {xxviii+530},
      ISBN = {0-89871-514-8},
   MRCLASS = {65N30 (65-02 65N15 74G15 74S05)},
  MRNUMBER = {1930132},
       DOI = {10.1137/1.9780898719208},
}

@book {MR3014456,
    AUTHOR = {Roub\'{\i}\v{c}ek, Tom\'{a}\v{s}},
     TITLE = {Nonlinear partial differential equations with applications},
    SERIES = {International Series of Numerical Mathematics},
    VOLUME = {153},
   EDITION = {Second},
 PUBLISHER = {Birkh\"{a}user/Springer Basel AG, Basel},
      YEAR = {2013},
     PAGES = {xx+476},
      ISBN = {978-3-0348-0512-4; 978-3-0348-0513-1},
   MRCLASS = {35-02 (35J60 35K55 35Q30)},
  MRNUMBER = {3014456},
       DOI = {10.1007/978-3-0348-0513-1},
}

@book {MR0603444,
    AUTHOR = {Temam, Roger},
     TITLE = {Navier-{S}tokes equations},
    SERIES = {Studies in Mathematics and its Applications},
    VOLUME = {2},
   EDITION = {Revised},
      NOTE = {Theory and numerical analysis, With an appendix by F. Thomasset},
 PUBLISHER = {North-Holland Publishing Co., Amsterdam-New York},
      YEAR = {1979},
     PAGES = {x+519},
      ISBN = {0-444-85307-3},
   MRCLASS = {35Q10 (49D99 65P05 76D05)},
  MRNUMBER = {603444},
}

@article{CIONCOLINI20192432,
author = {Andrea Cioncolini and Daniele Boffi},
title = {The MINI mixed finite element for the Stokes problem: An experimental investigation},
journal = {Computers \& Mathematics with Applications},
volume = {77},
number = {9},
pages = {2432-2446},
year = {2019},
issn = {0898-1221},
doi = {10.1016/j.camwa.2018.12.028},
}

@article{LIU,
  author  = {Liu, Jie},
  title   = {Open and traction boundary conditions for the incompressible {Navier--Stokes} equations},
  journal = {Journal of Computational Physics},
  volume  = {228},
  number  = {19},
  pages   = {7250--7267},
  year    = {2009},
  doi     = {10.1016/j.jcp.2009.06.021}
}

@book {MR3097958,
    AUTHOR = {Boffi, Daniele and Brezzi, Franco and Fortin, Michel},
     TITLE = {Mixed finite element methods and applications},
    SERIES = {Springer Series in Computational Mathematics},
    VOLUME = {44},
 PUBLISHER = {Springer, Heidelberg},
      YEAR = {2013},
     PAGES = {xiv+685},
      ISBN = {978-3-642-36518-8; 978-3-642-36519-5},
   MRCLASS = {65-02 (65M60 65N30)},
  MRNUMBER = {3097958},
MRREVIEWER = {Beny\ Neta},
       DOI = {10.1007/978-3-642-36519-5},
}

@book {MR2373954,
    AUTHOR = {Brenner, Susanne C. and Scott, L. Ridgway},
     TITLE = {The mathematical theory of finite element methods},
    SERIES = {Texts in Applied Mathematics},
    VOLUME = {15},
   EDITION = {Third},
 PUBLISHER = {Springer, New York},
      YEAR = {2008},
     PAGES = {xviii+397},
      ISBN = {978-0-387-75933-3},
   MRCLASS = {65-01 (65-02)},
  MRNUMBER = {2373954},
       DOI = {10.1007/978-0-387-75934-0},
}

@book {Guermond-Ern,
    AUTHOR = {Ern, Alexandre and Guermond, Jean-Luc},
     TITLE = {Theory and practice of finite elements},
    SERIES = {Applied Mathematical Sciences},
    VOLUME = {159},
 PUBLISHER = {Springer-Verlag, New York},
      YEAR = {2004},
     PAGES = {xiv+524},
      ISBN = {0-387-20574-8},
   MRCLASS = {65-02 (65M60 65N30 74S05 76M10 78M10)},
  MRNUMBER = {2050138},
MRREVIEWER = {R. S. Anderssen},
       DOI = {10.1007/978-1-4757-4355-5},
}

@book {MR3136903,
    AUTHOR = {Ciarlet, Philippe G.},
     TITLE = {Linear and nonlinear functional analysis with applications},
 PUBLISHER = {Society for Industrial and Applied Mathematics, Philadelphia,
              PA},
      YEAR = {2013},
     PAGES = {xiv+832},
      ISBN = {978-1-611972-58-0},
   MRCLASS = {46-02},
  MRNUMBER = {3136903},
MRREVIEWER = {Bhavana Deshpande},
}

@book {MR2424078,
    AUTHOR = {Adams, Robert A. and Fournier, John J. F.},
     TITLE = {Sobolev spaces},
    SERIES = {Pure and Applied Mathematics (Amsterdam)},
    VOLUME = {140},
   EDITION = {Second},
 PUBLISHER = {Elsevier/Academic Press, Amsterdam},
      YEAR = {2003},
     PAGES = {xiv+305},
      ISBN = {0-12-044143-8},
   MRCLASS = {46E35 (46-01 46-02 46B70 46Exx)},
  MRNUMBER = {2424078},
}

@book {MR851383,
    AUTHOR = {Girault, Vivette and Raviart, Pierre-Arnaud},
     TITLE = {Finite element methods for {N}avier-{S}tokes equations},
    SERIES = {Springer Series in Computational Mathematics},
    VOLUME = {5},
      NOTE = {Theory and algorithms},
 PUBLISHER = {Springer-Verlag, Berlin},
      YEAR = {1986},
     PAGES = {x+374},
      ISBN = {3-540-15796-4},
   MRCLASS = {65N30 (65-02 76-08)},
  MRNUMBER = {851383},
MRREVIEWER = {Max D. Gunzburger},
       DOI = {10.1007/978-3-642-61623-5},
}

@article {MR3232446,
    AUTHOR = {Oyarz\'{u}a, Ricardo and Qin, Tong and Sch\"{o}tzau, Dominik},
     TITLE = {An exactly divergence-free finite element method for a
              generalized {B}oussinesq problem},
   JOURNAL = {IMA J. Numer. Anal.},
  FJOURNAL = {IMA Journal of Numerical Analysis},
    VOLUME = {34},
      YEAR = {2014},
    NUMBER = {3},
     PAGES = {1104--1135},
      ISSN = {0272-4979},
   MRCLASS = {65N30 (65N12 76Dxx 86A05)},
  MRNUMBER = {3232446},
MRREVIEWER = {John A. Loustau},
       DOI = {10.1093/imanum/drt043},
}

@article {MR4337455,
    AUTHOR = {Sayah, Toni and Semaan, Georges and Triki, Faouzi},
     TITLE = {Finite element methods for the {D}arcy-{F}orchheimer problem
              coupled with the convection-diffusion-reaction problem},
   JOURNAL = {ESAIM Math. Model. Numer. Anal.},
  FJOURNAL = {ESAIM. Mathematical Modelling and Numerical Analysis},
    VOLUME = {55},
      YEAR = {2021},
    NUMBER = {6},
     PAGES = {2643--2678},
      ISSN = {0764-583X},
   MRCLASS = {65N30 (35K05 74S05 76S05)},
  MRNUMBER = {4337455},
       DOI = {10.1051/m2an/2021066},
       URL = {https://doi.org/10.1051/m2an/2021066},
}

@article {MR3119721,
    AUTHOR = {Congreve, Scott and Houston, Paul and S\"{u}li, Endre and
              Wihler, Thomas P.},
     TITLE = {Discontinuous {G}alerkin finite element approximation of
              quasilinear elliptic boundary value problems {II}: strongly
              monotone quasi-{N}ewtonian flows},
   JOURNAL = {IMA J. Numer. Anal.},
  FJOURNAL = {IMA Journal of Numerical Analysis},
    VOLUME = {33},
      YEAR = {2013},
    NUMBER = {4},
     PAGES = {1386--1415},
      ISSN = {0272-4979,1464-3642},
   MRCLASS = {65N30 (35Q35 65N15 65N50)},
  MRNUMBER = {3119721},
MRREVIEWER = {J\'{a}n\ Lov\'{\i}\v{s}ek},
       DOI = {10.1093/imanum/drs046},
       URL = {https://doi.org/10.1093/imanum/drs046},
}

@article {MR4911804,
    AUTHOR = {Campa\~na, Gilberto and Mu\~noz, Pablo and Ot\'arola, Enrique},
     TITLE = {Finite element approximation for a convective
              {B}rinkman-{F}orchheimer problem coupled with a heat equation},
   JOURNAL = {J. Sci. Comput.},
  FJOURNAL = {Journal of Scientific Computing},
    VOLUME = {104},
      YEAR = {2025},
    NUMBER = {1},
     PAGES = {Paper No. 21, 33},
      ISSN = {0885-7474,1573-7691},
   MRCLASS = {65N30 (35A01 35A02 35Q35)},
  MRNUMBER = {4911804},
       DOI = {10.1007/s10915-025-02932-y},
       URL = {https://doi.org/10.1007/s10915-025-02932-y},
}

@article {MR1370148,
    AUTHOR = {Lorca, Sebasti\'{a}n A. and Boldrini, Jos\'{e} Luiz},
     TITLE = {Stationary solutions for generalized {B}oussinesq models},
   JOURNAL = {J. Differential Equations},
  FJOURNAL = {Journal of Differential Equations},
    VOLUME = {124},
      YEAR = {1996},
    NUMBER = {2},
     PAGES = {389--406},
      ISSN = {0022-0396,1090-2732},
   MRCLASS = {35Q30 (35B45 35B65 76D05)},
  MRNUMBER = {1370148},
MRREVIEWER = {Nikolay\ G.\ Kuznetsov},
       DOI = {10.1006/jdeq.1996.0016},
       URL = {https://doi.org/10.1006/jdeq.1996.0016},
}

@article {MR3423229,
    AUTHOR = {Alvarez, Mario and Gatica, Gabriel N. and Ruiz-Baier, Ricardo},
     TITLE = {An augmented mixed-primal finite element method for a coupled
              flow-transport problem},
   JOURNAL = {ESAIM Math. Model. Numer. Anal.},
  FJOURNAL = {ESAIM. Mathematical Modelling and Numerical Analysis},
    VOLUME = {49},
      YEAR = {2015},
    NUMBER = {5},
     PAGES = {1399--1427},
      ISSN = {2822-7840,2804-7214},
   MRCLASS = {65N30 (35Q35 65N12 65N15 76D07 76R05)},
  MRNUMBER = {3423229},
       DOI = {10.1051/m2an/2015015},
       URL = {https://doi.org/10.1051/m2an/2015015},
}

@article {MR3959546,
    AUTHOR = {B\"{u}rger, Raimund and M\'{e}ndez, Paul E. and Ruiz-Baier,
              Ricardo},
     TITLE = {On {$H(\rm{div})$}-conforming methods for double-diffusion
              equations in porous media},
   JOURNAL = {SIAM J. Numer. Anal.},
  FJOURNAL = {SIAM Journal on Numerical Analysis},
    VOLUME = {57},
      YEAR = {2019},
    NUMBER = {3},
     PAGES = {1318--1343},
      ISSN = {0036-1429,1095-7170},
   MRCLASS = {65N30 (76R50 76S05)},
  MRNUMBER = {3959546},
       DOI = {10.1137/18M1196108},
       URL = {https://doi.org/10.1137/18M1196108},
}

@article {MR4659334,
    AUTHOR = {Allendes, Alejandro and Campa\~{n}a, Gilberto and Ot\'{a}rola,
              Enrique},
     TITLE = {Numerical discretization of a {D}arcy-{F}orchheimer problem
              coupled with a singular heat equation},
   JOURNAL = {SIAM J. Sci. Comput.},
  FJOURNAL = {SIAM Journal on Scientific Computing},
    VOLUME = {45},
      YEAR = {2023},
    NUMBER = {5},
     PAGES = {A2755--A2780},
      ISSN = {1064-8275,1095-7197},
   MRCLASS = {65N15 (35R06 65N50 76S05)},
  MRNUMBER = {4659334},
       DOI = {10.1137/22M1536340},
       URL = {https://doi.org/10.1137/22M1536340},
}

@article {MR5061536,
    AUTHOR = {Alekseev, G. V. and Soboleva, O. V.},
     TITLE = {Solvability of a boundary value problem for the stationary
              heat and mass transfer equations with variable coefficients},
      NOTE = {Translation of Differ. Uravn. {\bf 61} (2025), no. 12,
              1603--1619},
   JOURNAL = {Differ. Equ.},
  FJOURNAL = {Differential Equations},
    VOLUME = {61},
      YEAR = {2025},
    NUMBER = {12},
     PAGES = {1853--1868},
      ISSN = {0012-2661,1608-3083},
   MRCLASS = {35K05 (35A01 35A02)},
  MRNUMBER = {5061536},
       DOI = {10.1134/s001226612512002x},
       URL = {https://doi.org/10.1134/s001226612512002x},
}

@article {TaylorHood1973,
    AUTHOR = {Taylor, C. and Hood, P.},
     TITLE = {A numerical solution of the {N}avier-{S}tokes equations using
              the finite element technique},
   JOURNAL = {Internat. J. Comput. \& Fluids},
  FJOURNAL = {An International Journal. Computers \& Fluids},
    VOLUME = {1},
      YEAR = {1973},
    NUMBER = {1},
     PAGES = {73--100},
      ISSN = {0045-7930},
   MRCLASS = {76.65},
  MRNUMBER = {339677},
MRREVIEWER = {Nima\ Geffen},
       DOI = {10.1016/0045-7930(73)90027-3},
}

@article {MR0799997,
    AUTHOR = {Arnold, D. N. and Brezzi, F. and Fortin, M.},
     TITLE = {A stable finite element for the {S}tokes equations},
   JOURNAL = {Calcolo},
  FJOURNAL = {Calcolo. A Quarterly on Numerical Analysis and Theory of Computation},
    VOLUME = {21},
      YEAR = {1984},
    NUMBER = {4},
     PAGES = {337--344},
      ISSN = {0008-0624},
   MRCLASS = {65N30 (76-08 76D07)},
  MRNUMBER = {799997},
MRREVIEWER = {Gheorghe\ Micula},
       DOI = {10.1007/BF02576171},
}

@article {MR1011446,
    AUTHOR = {Scott, L. Ridgway and Zhang, Shangyou},
     TITLE = {Finite element interpolation of nonsmooth functions satisfying
              boundary conditions},
   JOURNAL = {Math. Comp.},
  FJOURNAL = {Mathematics of Computation},
    VOLUME = {54},
      YEAR = {1990},
    NUMBER = {190},
     PAGES = {483--493},
      ISSN = {0025-5718,1088-6842},
   MRCLASS = {65D05 (65N30)},
  MRNUMBER = {1011446},
MRREVIEWER = {Qian\ Li},
       DOI = {10.2307/2008497},
       URL = {https://doi.org/10.2307/2008497},
}

@article {Clement1975,
    AUTHOR = {Cl\'{e}ment, Ph.},
     TITLE = {Approximation by finite element functions using local
              regularization},
   JOURNAL = {Rev. Fran\c{c}aise Automat. Informat. Recherche
              Op\'{e}rationnelle S\'{e}r. Rouge Anal. Num\'{e}r.},
  FJOURNAL = {Revue Fran\c{c}aise d'Automatique, Informatique et Recherche
              Op\'{e}rationnelle S\'{e}rie Rouge. Analyse Num\'{e}rique},
    VOLUME = {9},
      YEAR = {1975},
    NUMBER = {R-2},
     PAGES = {77--84},
      ISSN = {0397-9342,2777-3493},
   MRCLASS = {65N15},
  MRNUMBER = {400739},
MRREVIEWER = {Stephen\ Hilbert},
}

@book {MR0548867,
    AUTHOR = {Girault, V. and Raviart, P.-A.},
     TITLE = {Finite element approximation of the {N}avier-{S}tokes
              equations},
    SERIES = {Lecture Notes in Mathematics},
    VOLUME = {749},
 PUBLISHER = {Springer-Verlag, Berlin-New York},
      YEAR = {1979},
     PAGES = {vii+200},
      ISBN = {3-540-09557-8},
   MRCLASS = {65N30 (35Q10 76D05)},
  MRNUMBER = {548867},
MRREVIEWER = {W.\ Velte},
}
\end{document}